\documentclass{article}
\usepackage[utf8]{inputenc}
\usepackage{hyperref}

\title{Comparing Classes of Estimators: When Does Gradient Descent Beat Ridge Regression in Linear Models? }

\usepackage{times}
\usepackage{graphicx}
\usepackage{subcaption}
\usepackage{amsfonts}
\usepackage{amsmath,amssymb}
\usepackage{bbm}
\usepackage{xcolor}
\usepackage{pgfplots}

\usepackage[normalem]{ulem}

\usepackage{geometry}

\usepackage{enumitem}

\newcommand{\R}{\mathbb{R}}

\newcommand{\E}{\mathbf{E}}
\usepackage{algorithm,algpseudocode}
\newcommand{\indep}{\perp \!\!\! \perp}

\DeclareMathOperator{\argmin}{argmin}

\newtheorem{theorem}{Theorem}
\newtheorem{definition}{Definition}
\newtheorem{lemma}{Lemma}
\newtheorem{corollary}{Corollary}

\newtheorem{proposition}{Proposition}

\DeclareMathOperator{\trace}{Tr}

\definecolor{ed}{RGB}{225,0,0}

\newcommand{\bitem}{\begin{itemize}}
\newcommand{\eitem}{\end{itemize}}
\newcommand{\benum}{\begin{enumerate}}
\newcommand{\eenum}{\end{enumerate}}
\newcommand{\beq}{\begin{equation}}
\newcommand{\eeq}{\end{equation}}
\newcommand{\beqs}{\begin{equation*}}
\newcommand{\eeqs}{\end{equation*}}
\newcommand{\ep}{\varepsilon}
\newcommand{\hbeta}{\widehat{\beta}}

\usepackage{amsfonts} %mathbb
\newcommand{\oldep}{\epsilon}
\author{
Dominic Richards,\footnote{Amazon.} \footnote{Research performed whilst at the University of Oxford.}
\,
Edgar Dobriban,\footnote{Department of Statistics and Data Science, and Department of Computer and Information Science, University of Pennsylvania. E-mail: \texttt{dobriban@wharton.upenn.edu}.}
\,
and Patrick Rebeschini\footnote{Department of Statistics, Oxford University. E-mail: \texttt{patrick.rebeschini@stats.ox.ac.uk}.}
}

\begin{document}
\maketitle

\begin{abstract}
Methods for learning from data depend on various types of tuning parameters, such as 
penalization strength or step size.
Since performance can depend strongly on these parameters, it is important to 
compare \emph{classes} of estimators---by considering prescribed \emph{finite sets} of tuning parameters---
not just particularly tuned methods. 
In this work, we investigate classes of methods via the relative performance of the \emph{best method in the class}. 
We consider the central problem of linear regression---with a random isotropic ground truth---and investigate the estimation performance of two fundamental methods, 
gradient descent 
and ridge regression. 
We 
unveil the following phenomena.
(1) For general designs, constant stepsize gradient descent  outperforms ridge regression when the eigenvalues of the empirical data covariance matrix decay slowly,  as a power law with exponent less than unity. 
If instead the eigenvalues decay quickly, as a power law with exponent greater than unity or exponentially, we show that 
ridge regression outperforms gradient descent. 
(2) For orthogonal designs, we compute the exact minimax optimal \emph{class} of estimators (achieving min-max-\emph{min} optimality), 
showing it is equivalent to gradient descent with 
decaying learning rate. 
We find the sub-optimality of ridge regression and gradient descent with constant step size.
Our results highlight that statistical performance can  depend strongly
on tuning parameters. 
In particular, while optimally tuned ridge regression is the best estimator in our setting, 
it can be outperformed by gradient descent by an \emph{arbitrary/unbounded} amount
when both methods are only tuned over finitely many regularization parameters.

\end{abstract}
% \tableofcontents

\section{Introduction}

\sloppy Least squares regression is a workhorse in modern statistics, machine learning, engineering, and signal processing, often lying at the heart of more complex algorithms and methodologies. 
It is well known that various forms of shrinkage or regularization often improve estimation performance in the presence of noise \cite{stein1956inadmissibility,james1961estimation,stein1981estimation}. 
The optimal level of shrinkage required can depend upon \emph{unknown} properties of the problem, such as the signal-to-noise ratio. As a result, in practice the regularization strength is routinely chosen through a two-step procedure: first, a \emph{finite} class, or collection, of candidate models (estimators) is proposed by considering a range of regularization values. 
Then a candidate model is chosen from the class through a criterion---such as prediction or estimation error---which itself is estimated, say, through cross-validation using held-out data \cite{allen1974relationship,stone1977asymptotics,geisser1975predictive}. 
This leads to a natural statistical question: How should \emph{classes} of candidate models be chosen and compared? 

We investigate this fundamental question when classes of candidate models are compared through the \emph{best-in-class model}. Studying the best candidate model in a class allows us to understand the behavior and sensitivity of model classes with respect to their tuning parameters at a level of precision that has not been investigated before. 
Indeed, the standard approach is to compare estimators that are optimally tuned with respect to hyperparameters, say by balancing bias and variance. 
In contrast, our approach allows us to
evaluate model classes when their hyperparameters are restricted to a finite collection/grid; as they must be in practice. 

We instantiate this approach for two popular classes of estimators, the iterates of gradient descent (GD), and ridge regression with a varying regularization strength. Each of these classes is extremely important and widely used. 
Gradient descent is the canonical first-order optimization algorithm, and its variants are currently the method of choice for training state-of-the art machine learning models, including deep neural networks \cite{nesterov2013introductory,bubeck2015convex,krizhevsky2012imagenet,brown2020language,lecun2015deep,schmidhuber2015deep}. Ridge regression, and more broadly $\ell_2$ regularization, has been possibly the most fundamental approach to regularization in inverse problems and statistical estimation from the 1940s onwards \cite{tikhonov1943stability,stein1956inadmissibility,james1961estimation,hoerl1970ridge,stein1981estimation}. For modern deep neural networks, it is linked to the highly popular method of weight decay \cite{hanson1988comparing,loshchilov2017decoupled}.

\begin{figure}
    \centering
    \includegraphics[width=0.49\textwidth]{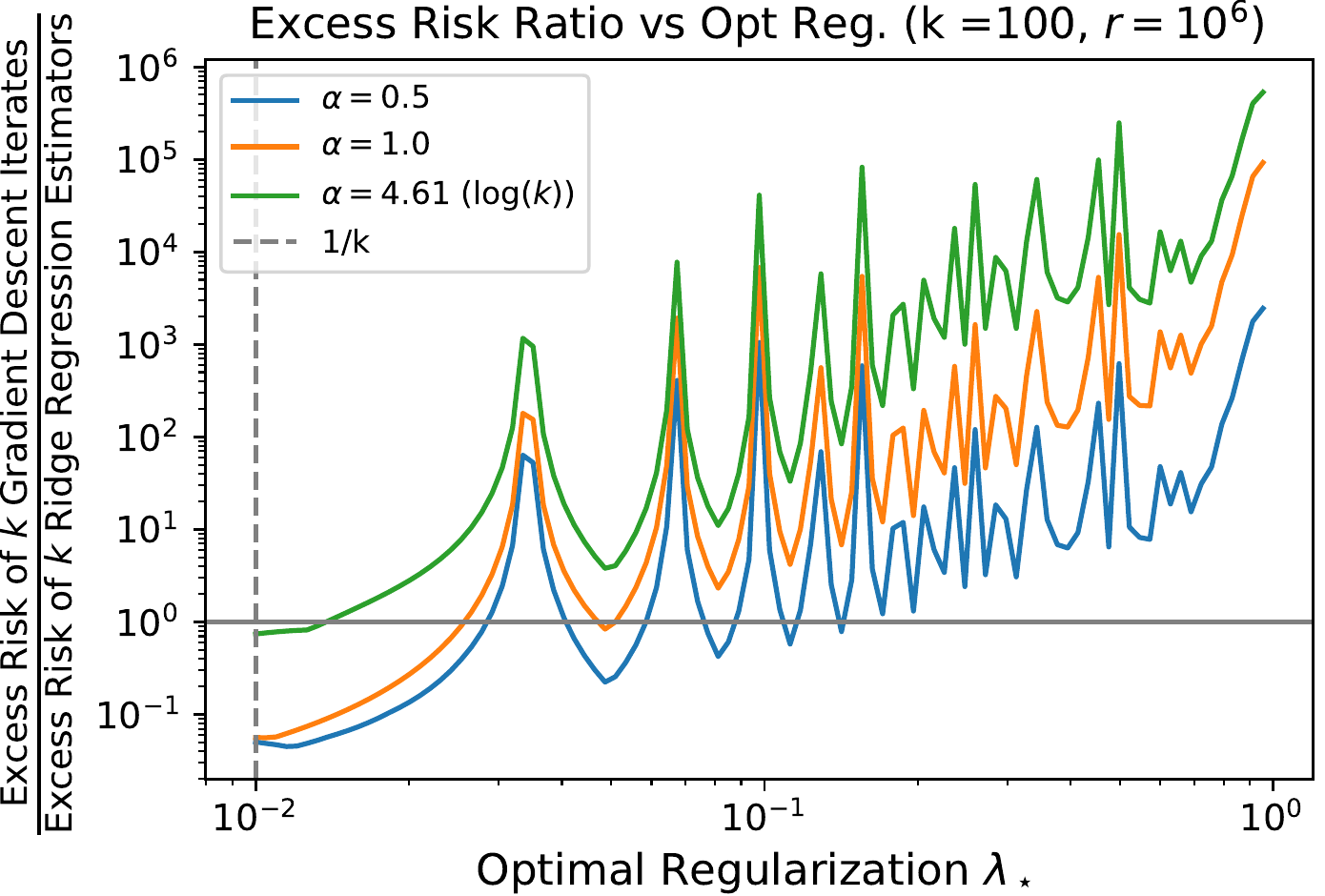}
    \includegraphics[width=0.49\textwidth]{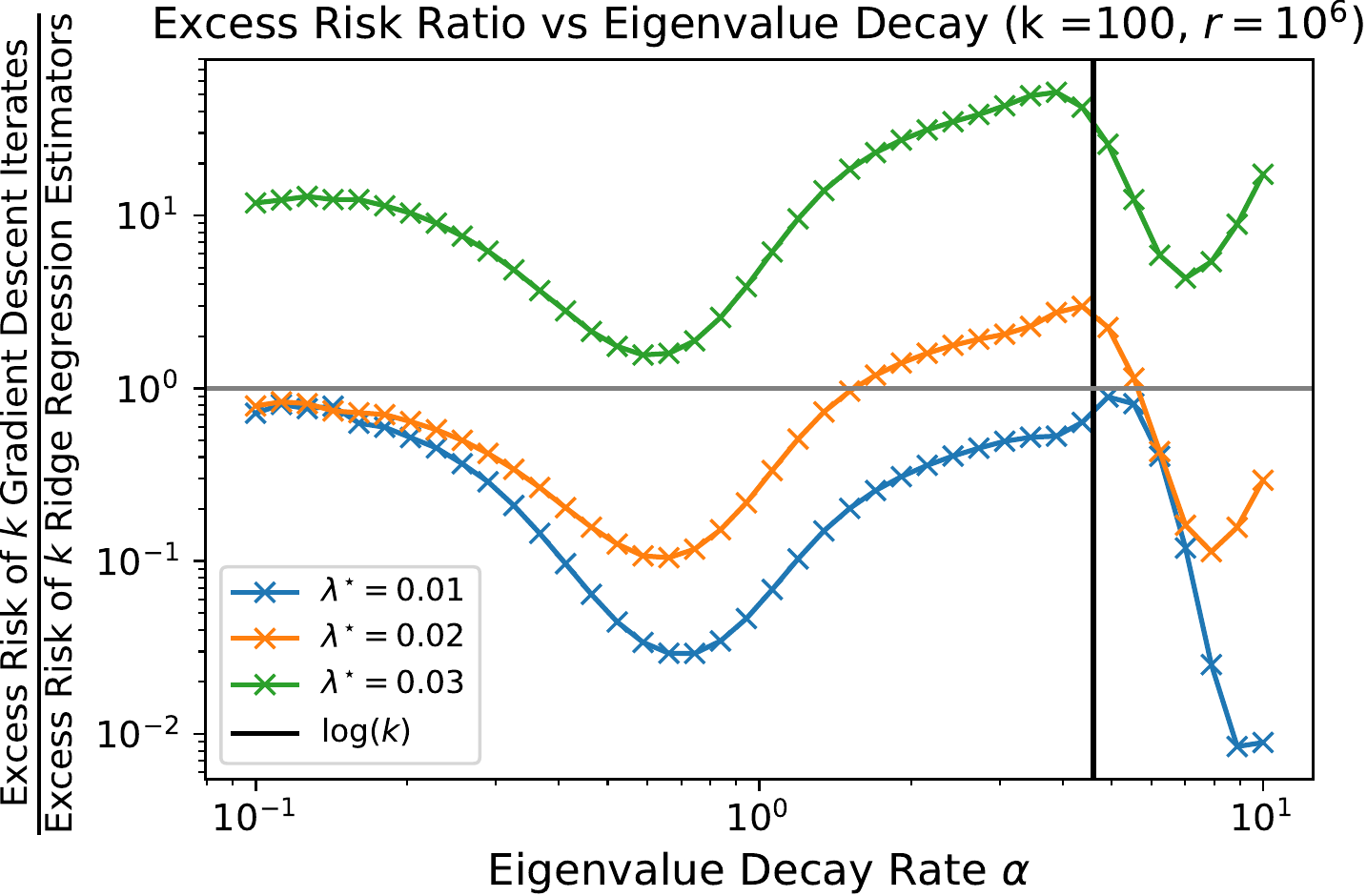}
    \caption{Ratio of the minimum excess risk achieved by $k$ iterates of gradient descent (with constant step size $\eta = O(1/k)$) and the minimum excess risk achieved by $k$ ridge regression estimators (with ridge regularization parameter spread uniformly in the interval $[0,1]$). Plot is as a function of the optimal ridge regularization parameter $\lambda_{\star} > 0$ (\textit{left}) and the exponent $\alpha \geq 0$ controlling the polynomial decay of the eigenvalues  $s_i = i^{-\alpha}$ for $i \ge 1$, of the empirical covariance matrix of the data (\textit{right}). 
    Dimension and sample size $d = n = 10^{6}$ with $k = 100$. 
    } 
    \label{fig:Compared}
\end{figure}

We develop an extensive set of theoretical results characterizing their relative performance for least squares regression,
in the central problem of linear regression---with a random isotropic ground truth parameter.
It is well known that \emph{optimally tuned} ridge regression is a linear Bayes method---and thus optimal---in this setting; henceforth, let $\lambda_\star > 0$ be the optimal choice of the ridge penalty term. 
In particular, the class of \emph{all} ridge regression estimators is at least as good as, and possibly much better than,  the class of \emph{all} gradient descent estimators.
One may thus expect that this is the case even for finite classes, when one considers a finite set of tuning parameters. 
However, we show that this is far from true. 
Along the way, we unveil several properties that do not seem to have been previously noticed in the literature.
These are summarised within our two main contributions:
\begin{enumerate}
    \item \textbf{Gradient Descent vs. Ridge Regression.}  
    We perform a detailed analysis of the ratio of excess risks (\emph{relative sub-optimality}, c.f.\ Definition \ref{res}) between classes of gradient descent estimators (with constant step size, i.e., with a step size independent of iteration time and only depending on the time horizon) and ridge regression estimators (with choices of regularization parameters discussed later; including uniform discretization). 
    Throughout, by \emph{excess risk of a class of estimators} we mean the best (i.e., minimal) excess risk achieved by the estimators in the class. 
    This turns out to depend strongly on the structure of the data, and on the rate of decay of the eigenvalues of the empirical covariance matrix of the data.
    \\
    \textbf{Fast Eigenvalue Decay.} When eigenvalues of the empirical covariance matrix decay fast--- according to a power law with exponent greater than unity, or exponentially fast---then \emph{ridge regression outperforms gradient descent}.
    This mirrors the behavior when all estimators are considered; and thus is the expected behavior.
    Precisely, there are problem sequences with a sample of size $n$, in dimension $d$, and number of estimators $k$,
    such that
    \footnote{Formal result in Corollary \ref{RiskLimit:Lower:Formal} Appendix \ref{sec:sub-opt-limit}.} 
    \begin{align*}
        \lim_{\substack{ k \rightarrow \infty \\ \lambda_{\star} \rightarrow 0}} \lim_{\substack{ d,n \rightarrow \infty}} 
        \frac{\text{Excess Risk of $k$ Gradient Descent Iterates}}{\text{Excess Risk of $k$ Ridge Regression Estimators}}
        = \infty.
    \end{align*}
    However, surprisingly, the performance of ridge regression relative to gradient descent behaves differently when the eigenvalues \emph{decay slowly} (Figure \ref{fig:Compared} and Theorem \ref{thm:Informal}).
    \\
    
    \textbf{Slow Eigenvalue Decay \& High SNR.} When eigenvalues of the empirical covariance matrix of the data decay according to a power law with exponent less than unity, and with a high signal-to-noise ratio (SNR), \emph{gradient descent outperforms ridge regression}. 
        In particular, with sample of size $n$, dimension $d$ and number of estimators $k$
        \footnote{Formal result in Corollary \ref{RiskLimit:Formal} Appendix \ref{sec:sub-opt-limit}.} 
        \begin{align*}
            \liminf_{\substack{ k \rightarrow \infty \\ \lambda_{\star}\rightarrow 0}} \lim_{\substack{ d,n \rightarrow \infty}} 
            \frac{\text{Excess Risk of $k$ Gradient Descent Iterates}}{\text{Excess Risk of $k$ Ridge Regression Estimators}}
            = 0.
        \end{align*} 
    The fact that the excess risk ratio tends to zero for appropriate sequences of problems 
    highlights that gradient descent can outperform ridge regression by an \emph{arbitrary factor}. 
    More generally, 
    Theorem \ref{thm:Informal} characterizes the excess risk ratio for finite $n,d,k,\lambda_{\star}$.
    \\
    
     Overall, the ratio of excess risks can behave in subtle ways, see Figure \ref{fig:Compared}. Precisely, we find that the excess risk ratio is highly oscillatory as a function of the optimal ridge regularization parameter (for other parameters fixed), 
     and has several regions of monotonicity as a function of the eigenvalue decay rate $\alpha$ (we identify three), see the discussion after Theorem \ref{thm:Informal}. 
     To control the ratio, we develop new  upper and lower bounds on the excess risk of individual classes of gradient descent and ridge regression estimators
    (Section \ref{b-ind}).
    
    \item \textbf{Min-Max-Min Optimality, Exact Computations, and Phase Transitions.} For orthogonal designs, i.e., when the empirical covariance matrix of the data is a scaled identity matrix, 
    we find the \emph{exact minimax optimal classes} of estimators over a range of signal-to-noise ratios  (Theorem \ref{mmxgrid}). We refer to these as \emph{min-max-min} optimal due to the formulation we consider, c.f.\ equation \eqref{mmxrisk}, which evaluates the minimum over classes, the maximum over the problems,  and  the minimum (best-in-class) performance. 
    Perhaps surprisingly, 
    this class is achieved \emph{neither}
    by ridge regression with a uniform discretization grid (or other common discretizations such as log-uniform) 
    \emph{nor} by gradient descent with a constant step size.
    However, it can be represented as gradient descent with a step size decreasing polynomially in the iteration counter (Section \ref{sec:app:compgrid}). 
    
    For orthogonal designs, we find the \emph{exact} sub-optimality of ridge regression and gradient descent (Section \ref{sec:SuboptCompGrid}) compared to the optimum. 
    We find that the sub-optimality has a ``switching" behavior as a function of the number $k$ of estimators in the class, behaving as either $1/k$ or $1/k^2$ for certain ranges of $k$; this matches numerical results (Figure 
    \ref{fig:MinimaxVsTheory}) and 
    results for general designs, c.f.\ equation \eqref{equ:uniform_vs_log_grid}.
\end{enumerate}

For further details on Figure \ref{fig:Compared}, the left plot illustrates the highly oscillatory behavior of the excess risk ratio (relative sub-optimality) of the two methods. This is partly due to the discreteness of the grid of tuning parameters: when a ridge regularization grid element falls close to the optimal regularizer, the best-in-class method performs very well. The right plot illustrates the non-trivial behavior of the same performance metric as a function of the geometry of the data/design matrix, as captured by the decay rate of the eigenvalues of the empirical covariance matrix. 
Our results identify several regions of monotonicity 
of these curves.

The code to reproduce our experimental results is available at \url{https://github.com/DominicRichards/ComparingGradientDescentRidge}.

\subsection{Related Literature}
\label{sec:lit}
From the early works of \cite{tikhonov1943stability,stein1956inadmissibility,james1961estimation,hoerl1970ridge,stein1981estimation}, a vast literature investigating the performance of shrinkage methods in the context of mean estimation, inverse problems, and least squares regression has developed. 
The PhD thesis by Thisted \cite{thisted1977}, in particular, studies ridge regression, making connections with minimax estimation and empirical Bayes methods, including a thorough historical review. \cite{thisted1977} argues for the importance of making a distinction between ridge regression estimators for a fixed regularization, for adaptively chosen regularization, and for classes of models with varying regularization strengths. Our work makes steps towards rigorously formalizing such a distinction, focusing on the performance of classes of models defined by varying the regularization strength. General references on ridge regression include \cite{gruber2017improving,van2015lecture}.

We  highlight works related to inverse problems and kernel regression, such as \cite{landweber1951iteration,engl1996regularization,bissantz2007convergence,caponnetto2007optimal,yao2007early,bauer2007regularization,raskutti2014early,rosasco2015learning,blanchard2018optimal,pagliana2019implicit,lin2020optimal}, where the shrinkage performance of iterative methods like gradient descent has gained much attention. Here statistical performance is encoded by an \emph{optimal rate} at which the estimation or prediction error decays with respect to the sample size \cite{caponnetto2007optimal,blanchard2018optimal}. The rate depends on both the regularity of the problem (e.g., smoothness of the optimal regression function), and the effective dimension (often referred to as the capacity) of the kernel. The guarantees for various estimators---e.g., ridge regression and early stopped gradient descent---aim to show that the optimal rate is achieved in a worst case sense. However, the bounds usually do not capture the dependence on the discretization, which is a key component of our results.

The concurrent work \cite{zou2021benefits} compares the sample inflation, in a Bahadur efficiency sense \cite{bahadur1967rates,bahadur1971some}, of single-pass stochastic gradient descent (SGD) to ridge regression i.e., the additional sample size required by SGD in order to have comparable or better performance than ridge regression, when each is optimally tuned. Notably, they observe (Corollary 5.1 within) when the population eigenvalues decay with a power law with exponent less than unity, then the sample inflation is mild (logarithmic). This finding is connected to our results (Theorem \ref{thm:Informal}), in that we find that gradient descent can outperform ridge regression when the empirical eigenvalues decay at this rate. However, our perspective and results are different and complementary. For instance, we consider a fixed-design setting instead of a random-design setting, and focus on the sensitivity to tuning parameters with respect to the underlying unknown signal-to-noise level, including the effects of discretization, instead of looking at optimal tuning.        

The two  concurrent works \cite{sheng2022accelerated,wang2022implicit} use an error decomposition similar to Lemma \ref{lem:ErrorDecomp} to explore the excess risk of spectral filters arising from continuous time accelerated gradient descent. While our work is linked to excess risk (see the relative sub-optimality definition \ref{res} and Section \ref{b-ind-RR}), we study finite classes of estimators. 
The continuous time limit could perhaps be associated to an infinitely sized class, but this is distinct from our work. 
A natural future direction is to explore 
using our problem dependent class bounds for ridge regression or gradient descent, to compare with accelerated gradient descent. For extended related literature see supplementary material \ref{sec:app:lit}.

\subsection{Setup} 
\label{sec:Setup}
Consider a standard noisy least squares regression setting with $n$ data points of dimension $d$. We observe a design matrix $X \in \mathbb{R}^{n \times d}$ (also referred to as a data matrix, inputs, features, or covariates) and  response $Y \in \mathbb{R}^{n}$ (also referred to as outcomes or outputs) such that
\begin{align}
\label{equ:Response}
    Y = X \beta_{\star} + \sigma \oldep,
\end{align}
where the ground truth regression parameter $\beta_{\star} \in \mathbb{R}^{d}$,  noise $\oldep \in \mathbb{R}^{n}$, and noise level $\sigma > 0$ are unknown. 
We assume that the noise and ground truth $(\oldep,\beta_{\star})$ are random and satisfy the following distributional assumptions: 
\begin{align}
\label{equ:NoiseSignal:ass}
    \E[\oldep] = 0, \quad\quad \E[\oldep \oldep^{\top}] = I_n,
    \quad\quad \E[\beta_{\star} \beta_{\star}^{\top}]  = \frac{\psi}{d} I_d,
    \quad\quad \oldep \indep \beta_{\star}.
\end{align}
The noise $\oldep$ and ground truth $\beta_{\star}$ are mutually independent, each with zero mean  and an isotropic covariance  matrix. The scalar $\psi >0$ controls the expected squared norm of $\beta_{\star}$, as $\E \|\beta_{\star}\|_2^2 = \psi$, and we will refer to it as the \emph{signal strength}. 
Given $X$ and $Y$, we focus on estimating the associated ground truth $\beta_{\star}$ with squared $\ell_2$-norm loss, 
for estimators $\widehat{\beta}$
\begin{align*}
    L_{\beta_{\star}}(\widehat{\beta}) := \|\widehat{\beta} - \beta_{\star}\|_2^2.
\end{align*}

We denote the eigenvector-eigenvalue decomposition of the un-centered empirical covariance matrix as $n^{-1}X^{\top}X $ $=$ $\sum_{i=1}^{r} s_i v_i v_i^{\top}$, where the rank is denoted $r \leq \min(n,d)$, the non-zero eigenvalues  $\{s_i\}_{i=1}^{r}$, and associated eigenvectors $\{v_i\}_{i=1}^{r}$. 
For a measurable function $\Phi:\mathbb{R}_{\geq 0} \rightarrow \mathbb{R}_{\geq 0}$, where $\mathbb{R}_{\geq 0} = [0,\infty)$, we then denote the application of a function to a matrix by applying it to the eigenspectrum, i.e., $\Phi(n^{-1}X^{\top} X) = \sum_{i=1}^{r} \Phi(s_i) v_i v_i^{\top}$.\footnote{Sometimes we will refer to $\Phi$ as a \emph{spectral shrinker}.}
Given this, we will study estimators 
$\widehat{\beta}_{\Phi}:\R^{n(p+1)}\to \R^{p}$
of the form
\begin{align}\label{hbphi}
    \widehat{\beta}_{\Phi} = \widehat{\beta}_{\Phi}(X,Y) = \Phi\Big(\frac{X^{\top} X}{n}\Big) \frac{X^{\top}Y}{n}.
\end{align}
We now introduce two estimators studied in this work, beginning with ridge regression. 
\begin{definition}[Ridge Regression]
\label{def:RR}
Define ridge regression with penalty $\lambda > 0$ as \footnote{The expression $\widehat{\beta}_{\lambda}$ is a slight abuse of notation  compared to $\widehat{\beta}_{\Phi_\lambda}$, and is only used for convenience.} 
\begin{align*}
    \widehat{\beta}_{\lambda} := \Big( \frac{X^{\top}X}{n} + \lambda I_d \Big)^{-1} \frac{X^{\top} Y}{n}
    = \argmin_{\beta \in \mathbb{R}^{d}} \Big\{ \frac{1}{n} \|X \beta - Y \|_2^2 + \lambda \|\beta\|_2^2 \Big\}.
\end{align*}
This corresponds to the spectral shrinker 
$\Phi_\lambda$  in \eqref{hbphi}, such that for any $u\ge 0$,
$\Phi_\lambda(u) = 1/(u+\lambda)$.  
\end{definition}
\sloppy The second estimator arises from the iterates of gradient descent applied to the squared loss  $\frac{1}{2n} \|X \beta - Y\|_2^2$ with a fixed stepsize.
\begin{definition}[Gradient Descent]
\label{def:GD}
Initialized at $\widehat{\beta}_{\eta,0} = 0$, define the iterates of gradient descent recursively for stepsize $\eta \geq 0$ and iteration $t \geq 0$ as
\begin{align*}
    \widehat{\beta}_{\eta,t+1}
    : = 
    \widehat{\beta}_{\eta,t} 
    - 
    \frac{\eta }{n} X^{\top} (X \widehat{\beta}_{\eta,t} - Y)
    = 
    \sum_{\ell=0}^{t} \eta \left(I_d - \eta \frac{X^{\top}X}{n}\right)^{\ell} \frac{X^{\top} Y}{n}.
\end{align*}
This corresponds to the spectral shrinker
$\Phi$  in \eqref{hbphi}, such that for any $u\ge 0$,
$\Phi(u) = \eta \sum_{\ell=0}^{t-1}(1-\eta u)^{\ell}$.
\end{definition}

The optimal choice of regularization parameter for the ridge regression estimator is  $\lambda_{\star} =  \sigma^2 d/(\psi n)$, as shown in the following known lemma, which states that ridge regression with the appropriate regularization parameter is a \emph{linear Bayes} estimator, see \cite{hartigan1969linear,rao1975simultaneous, gruber2017improving}.
\begin{lemma}
\label{lem:ErrorDecomp}
Consider data generated according to \eqref{equ:Response} and \eqref{equ:NoiseSignal:ass}.  Then for any $\Phi: \mathbb{R}_{\geq 0} \rightarrow \mathbb{R}_{\geq 0}$ we have  $\E_{\beta_{\star},\oldep}[L_{\beta_{\star}}(\widehat{\beta}_{\Phi})] \geq \E_{\beta_{\star},\oldep}[L_{\beta_{\star}}(\widehat{\beta}_{\lambda_{\star}})]$.
\end{lemma}
Motivated by this, we define the expected excess risk (or, sub-optimality) of an estimator $\hbeta$ by  
\beq\label{subopt}
\mathcal{E}(\widehat{\beta})=
\E_{\beta_{\star},\oldep}[L_{\beta_{\star}}(\widehat{\beta})-
L_{\beta_{\star}}(\widehat{\beta}_{\lambda_\star})].
\eeq

\subsubsection{Comparing classes via the relative sub-optimality}

It is clear from Lemma \ref{lem:ErrorDecomp} that when the noise variance $\sigma^2$ and signal strength $\psi$ is known, it is not possible to outperform ridge regression in expectation. However, in practice $\sigma^2$ and $\psi$ are usually not known, and may be hard to estimate. This leaves open the choice of a ``best" estimator.

One approach is to use hyperparameter search to choose the best regularization parameter from a class. In general, one can consider a class $\mathcal{C} = \{\widehat{\beta}^i\}_{i=1}^{k}$ of $k \ge 1$ estimators $\widehat{\beta}^i$, for $i=1,\dots,k$. For instance, for ridge regression we can choose a  grid $\Gamma = \{0,1/(k-1),2/(k-1),\dots,1\}$ of size $k$ over the unit interval $[0,1]$, and set $\mathcal{C} = \{\widehat{\beta}_\lambda \}_{\lambda \in \Gamma}$.  

We now introduce the measure by which we compare \emph{two  classes} of models. This definition applies to any estimation problem, and any classes of estimators, not just linear regression with quadratic loss and the specific estimators that we focus on in this paper. 
\begin{definition}[Relative Sub-optimality]\label{res}
For two classes of estimators $\mathcal{C}_1 = \{\widehat{\beta}^{u}\}_{u \in U_1}$ and $\mathcal{C}_2 = \{ \widehat{\beta}^{u}\}_{u \in U_{2}}$ define their \emph{relative sub-optimality} using the sub-optimality $\mathcal{E}$ from \eqref{subopt} 
as\footnote{We assume $\min_{\widehat{\beta} \in \mathcal{C}_2} \mathcal{E}(\widehat{\beta})>0$.}
\begin{align}\label{rs}
    \mathcal{S}(\mathcal{C}_1,\mathcal{C}_2)
     := 
     \frac{\min_{\widehat{\beta} \in \mathcal{C}_1} \mathcal{E}(\widehat{\beta})} 
    {\min_{\widehat{\beta} \in \mathcal{C}_2} \mathcal{E}(\widehat{\beta})}
    =
    \frac{ 
    \min_{\widehat{\beta} \in \mathcal{C}_1 } 
    \E_{\beta_{\star},\oldep}[L_{\beta_{\star}}(\widehat{\beta})]
    - 
    \E_{\beta_{\star},\oldep}[L_{\beta_{\star}}(\widehat{\beta}_{\lambda_{\star}})]}{
    \min_{\hbeta \in \mathcal{C}_2 } 
    \E_{\beta_{\star},\oldep}[L_{\beta_{\star}}(\widehat{\beta})]
    - 
    \E_{\beta_{\star},\oldep}[L_{\beta_{\star}}(\widehat{\beta}_{\lambda_{\star}})]
    }.
\end{align}
\end{definition}
The relative sub-optimality $\mathcal{S}(\cdot,\cdot)$ is directly related to the excess risk of each class when compared to optimally tuned ridge regression. Indeed, if  we have $m_1 \geq 0, m_2>0$ such that
\begin{align*}
    \min_{\widehat{\beta} \in \mathcal{C}_{i}} \E_{\beta_{\star},\oldep}[L_{\beta_{\star}}(\widehat{\beta})] = (1+m_i) \E_{\beta_{\star},\oldep}[L_{\beta_{\star}}(\widehat{\beta}_{\lambda_{\star}})]
\end{align*}
then precisely $\mathcal{S}(\mathcal{C}_1,\mathcal{C}_2) = m_1/m_2$.

\section{Gradient Descent versus Ridge Regression}
\label{gd-vs-ridge}
In this section we compare the classes of models produced by the regularization path of gradient descent to ridge regression with a grid of penalization parameters. 
For simplicity of notation, in this section we assume without loss of generality that $\psi=1$. This can always be achieved by rescaling the problem.
We assume that there is a \emph{known} lower bound 
$\lambda_{\min} > 0$
on the optimal ridge regularization parameter, 
such that $\lambda_{\star} \in [\lambda_{\min},1]$.
We will specify $\lambda_{\min}$ in each result to follow. 
This is equivalent to having a known lower bound $\sigma_{\min} > 0$ on the variance of the noise,  so that $\sigma \ge \sigma_{\min}$, with the relationship 
\beq\label{lst}
\lambda_{\star} = \frac{d}{n} \sigma^2 \geq \frac{d}{n} \sigma_{\min}^2 = \lambda_{\min}.
\eeq
Given this, we define the following two classes of models, for $k\ge 2$\footnote{The class $\mathcal{C}^{\text{GD}}(\eta,k)$ is defined for $k\ge1$, but we only consider $k\ge 2$, because we want to have an equal number of GD and ridge estimators.}
\begin{align}
    \mathcal{C}^{\text{GD}}(\eta,k) 
    & := 
    \{ \widehat{\beta}_{\eta,t} : 1 \leq t \leq k\},\label{cl}\\
    \mathcal{C}^{\text{Ridge}}(\lambda_{\min},k) 
    & := 
    \left\{ \widehat{\beta}_{\lambda} : \lambda \in \{\lambda_{\min},\lambda_{\min}+\delta,\lambda_{\min} + 2\delta,\dots,1\}, \delta = \frac{1-\lambda_{\min}}{k-1} \right\},
    \label{cl2}
\end{align}
The first class $\mathcal{C}^{\text{GD}}(\eta,k) $ contains the first $k$ iterates of gradient descent.
Meanwhile, $\mathcal{C}^{\text{Ridge}}(\lambda_{\min},k)$ is the class of ridge regression estimators with regularization parameters chosen from \emph{uniformly} discretizing the interval $[\lambda_{\min},1]$, with window $\delta>0$.

Since the two sets are of the same size $k$, 
it is natural to scale $\eta \!= \! (k \lambda_{\min})^{-1}$. In this case, if $\lambda_{\min} = \Theta( \lambda_{\star})$, then we have $k = O((\eta \lambda_{\star})^{-1})$, which aligns with the scaling within prior work on inverse problems and linear regression \cite{RosascoIterateAveraging,ali2019continuous,lin2020optimal}.\footnote{Intuitively, this scaling ensures the shrinkage of gradient descent  in the direction associated to the smallest eigenvalue 
is \emph{at least as strong} as for optimally tuned ridge regression. 
That is, for an eigenvalue $s > 0$, let $k^{\star}(s)$ be the theoretical number of gradient descent iterations---not necessarily an integer---that matches optimally tuned ridge regression, so $(1-\eta s)^{k^{\star}(s)} = \lambda_{\star}/(\lambda_\star + s)$. 
Noting that $k^{\star}(s)$ is a decreasing function in $s$, the maximum number of gradient descent iterations (to match ridge regression) occurs when $s \rightarrow 0$ in which case  $k^{\star}(s) = \log(1+s/\lambda_{\star})/\log(1/(1-\eta s)) \rightarrow 1/(\eta \lambda_{\star})$.}
Now, when $\lambda_{\min}$ is of a smaller order, one should either perform additional iterations if the step size is fixed, or choose a larger step size if the number of iterations is fixed, to ensure gradient descent includes sufficiently regularized models.

We now give a theorem which describes the main insights in this section.  
We recall that $A \gtrsim B$ if there is a constant $C$ independent of $d,n,\sigma,k$ such that $A \ge CB$.\footnote{When this statement refers to an assumption, the constant needs to be sufficiently large.} 
We write $A\lesssim B$ when $B \gtrsim A$, and 
$A\simeq B$ when $A\lesssim B$ and $B\lesssim A$.
For a real number $a$
and a set of real numbers $\Gamma = \{ b_1,b_2,\dots,b_k\}$ let
$\text{Dist}(a,\Gamma) = \min_{b \in \Gamma} |a-b|$.
We denote the distance scaled by a real number  $ c > 0$ as
\beq\label{dist}
\text{Dist}_c(a,\Gamma) = \text{Dist}(a,\Gamma)/c = \frac{\min_{b \in \Gamma} |a-b|}{c}.
\eeq
For a number $k > 1$ of classes, let $\delta = (1-\lambda_{\min})/(k-1)$ be size of the window of a grid, and define the grid 
$$\Gamma = \{\lambda_{\min},\lambda_{\min}+\delta,\lambda_{\min} + 2\delta,\dots,1\}$$
of ridge regularization parameters. 
Further, let  $\eta =1/(k\lambda_{\min})$ the learning rate for gradient descent. 
To avoid a degenerate problem 
we assume that the optimal ridge regularization parameter does not intersect the discretization, i.e., $\lambda_{\star} \not\in \Gamma = \{\lambda_{\min},\lambda_{\min}+\delta,\lambda_{\min} + 2\delta,\dots,1\}$. Similarly, we assume $(\eta \lambda_\star)^{-1} \not\in \{0,1,2,3,\dots,k-1\}$, which is associated to the ``optimal'' number of gradient descent iterations. 
\begin{theorem}
\label{thm:Informal}
In the random-effects linear model \eqref{equ:Response},
\eqref{equ:NoiseSignal:ass} in the above setting,
assume 
that the lower bound $\sigma_{\min}$ on the noise level $\sigma$ obeys $0 < \sigma_{\min} \leq \sigma /\sqrt{2}$
 and
the number $k$ of classes is large enough that, for the optimal regularization parameter $\lambda_{\star}$ and its known lower bound $\lambda_{\min}$ from \eqref{lst}
$$k \gtrsim \lambda_{\star}^{-1}\max\left\{1,\frac{1}{\lambda_{\min}} \log\left(1+\frac{1}{\lambda_{\star}}\right)\right\}.$$ 
We then have the following results, under two sets of conditions for the non-zero eigenvalues  $\{s_i\}_{i=1}^{r}$, of  $n^{-1}X^{\top}X$.\footnote{See the related work section for a motivation of these conditions.} 
\begin{itemize}[leftmargin=*]
    \item \textbf{Slow Power Law Decay}: Suppose $s_i = i^{-\alpha}$, $i=1,\ldots, r$, for $\alpha \in (0,1)$. Then there is a function  $r_{\alpha,\lambda_{\star},k,\lambda_{\min}}  > 1$ 
    of $\alpha,\lambda_{\star},k,\lambda_{\min}$
    such that if the rank $r \geq r_{\alpha,\lambda_{\star},k,\lambda_{\min}}$, the relative sub-optimality \eqref{rs} of the class of gradient descent and ridge regression estimators from \eqref{cl}, \eqref{cl2}
    has the order, with the scaled distance from \eqref{dist} 
    \begin{align*}
        \mathcal{S}(\mathcal{C}^{\text{GD}}(\eta,k),\mathcal{C}^{\text{Ridge}}(\lambda_{\min},k)) 
        \simeq 
        \frac{1}{\mathrm{Dist}_{\delta}(\lambda_{\star},\Gamma)^2}
        \left(\frac{d}{n}\right)^2 \left(\frac{\sigma^{2}}{\sigma_{\mathrm{min}}}\right)^4.
    \end{align*}
    \item \textbf{Fast Power Law or Exponential Decay}: Suppose $s_i = i^{-\alpha}$, $i=1,\ldots, r$, for 
    $\log(k) \gtrsim \alpha > 1$ or $s_i = \exp(-\rho(i-1))$, $i=1,\ldots, r$, for $\log(k) \gtrsim \rho > 0$. Further, suppose that $s_r \leq \lambda_{\star}$. Then the relative sub-optimality \eqref{rs} of the class of gradient descent and ridge regression estimators from \eqref{cl}, \eqref{cl2} is lower bounded as
    \begin{align*}
        \mathcal{S}(\mathcal{C}^{\text{GD}}(\eta,k),\mathcal{C}^{\text{Ridge}}(\lambda_{\min},k)) 
        \gtrsim 
        \frac{1}{\mathrm{Dist}_{\delta}(\lambda_{\star},\Gamma)^2} 
        \left(\frac{d}{n}\right)^2 \left(k \sigma^{2}\right)^2.
    \end{align*}
\end{itemize}
\end{theorem}

The above theorem provides insights into the relative sub-optimality of gradient descent and ridge regression under different assumptions of the eigenvalue decay. 
The bounds are scaled by $\text{Dist}_{\delta}(\lambda_{\star},\Gamma)^{-2}$ (recall \eqref{dist}), which encodes the distance between the optimal regularization $\lambda_{\star}$ and the ridge regression discretization grid $\Gamma$. 
In particular, if $\lambda_\star = \lambda_{\min} + (j+\ep)\delta$ for an integer $k-2 \geq j \geq 0$ and $\ep \in (0,1)$ then $\text{Dist}_{\delta}(\lambda_{\star},\Gamma) = \min\{1-\ep,\ep\}$. 
Therefore, if $\lambda_{\star}$ lands in the middle of two discretization points---so $\ep = 1/2$---we find $\text{Dist}_{\delta}(\lambda_{\star},\Gamma)  = 1/2$.
Alternatively,  if the optimal regularization $\lambda_\star$ nears a point within the discretization $\Gamma$---so $\ep \rightarrow 0$ or $1$---we then get  $\text{Dist}_{\delta}(\lambda_{\star},\Gamma) \rightarrow 0$, reflecting that the denominator of $\mathcal{S}(\mathcal{C}^{\text{GD}}(\eta,k) ,\mathcal{C}^{\text{Ridge}}(\lambda_{\min},k) )$ vanishes. In this case ridge regression becomes ``infinitely more accurate" than gradient descent; which is the expected behavior.

Theorem \ref{thm:Informal} then considers two main types of problem structures. 
In the first case, when the eigenvalues decay slowly (at a power law rate of $s_i=i^{-\alpha}$ for $0 < \alpha < 1$), the ratio is on the order of $\left(d/n\right)^2 \left(\sigma^{2}/\sigma_{\text{min}}\right)^4$ provided the rank $r=\min(n,d)$ is sufficiently large with respect to the problem parameters $(\alpha,\lambda_{\star})$ as well as the model parameters $(k,\lambda_{\min})$. When $d \simeq n$ and $\sigma_{\min} = \Theta(\sigma)$, the relative sub-optimality is on the order of $O(\sigma^4)$, and thus, perhaps surprisingly, gradient descent performs better then ridge regression in low noise settings. 
In particular, if $\sigma\to0$, gradient descent can outperform ridge regression by an \emph{arbitrary factor}. 

In general, the smaller the \emph{a priori} lower bound $\sigma_{\min}$, the less efficient gradient descent becomes. The reason is that it ``wastes'' models in this case, because the number of models $k$ is fixed, and thus the step size $\eta = O(1/(k\sigma_{\min}^2))$ is large. 
Therefore gradient descent takes large steps with most of the iterations being close to the convergence point, i.e., the least squares estimator.
In particular, the bound is of a constant  order if $\sigma_{\min} \simeq \sigma^2$. 
Similar results hold in the ``low-dimensional'' setting where $d = n^{q}$ for some $q \in [0,1)$, details of which are provided in the results to follow.

The second case of Theorem \ref{thm:Informal} considers eigenvalues that decay \emph{more quickly}, namely, with a power law so that $s_i = i^{-\alpha}$ for $\log(k) \gtrsim \alpha > 1$, or exponentially quickly. In either of these cases we see that the relative sub-optimality grows at least with the \emph{square of the number of models $k$ for any rank $r$}. This implies that ridge regression is the better choice as number of models grows; which aligns with our expectations based on considering all ridge estimators. 

When considering the monotonicity of the relative sub-optimality $\mathcal{S}(\mathcal{C}^{\text{GD}}(\eta,k),\mathcal{C}^{\text{Ridge}}(\lambda_{\min},k))$ as a function of the eigenvalue polynomial decay rate $\alpha \geq 0$ i.e., $s_i = i^{-\alpha}$, we note an alignment between 
the empirical evaluation in Figure \ref{fig:Compared} and our theoretical results. Precisely, we highlight three regimes: \textbf{(1)} a slow decay $\alpha \in (0,1)$, which is favorable to gradient descent (Theorem \ref{thm:Informal}); \textbf{(2)}, a faster decay $\log(k) \gtrsim \alpha \geq 1$, which is favorable to ridge regression (Theorem \ref{thm:mainGDvsRidge_Upper:poly}) ; \textbf{(3)}, a very fast decay $\alpha \gtrsim \log(k)$, which we conjecture is favorable to gradient descent. 
The cases of very small or large $\alpha$ are consistent with our findings from the orthogonal design setting in Section \ref{mmx-ortho-sec}, which intuitively corresponds to $\alpha\to0$ or $\alpha\to\infty$.

\textbf{Computational Questions.} Within the results above, we fix the class sizes to be equal, so $|\mathcal{C}^{\text{GD}}(\eta,k)| = |\mathcal{C}^{\text{Ridge}}(\lambda_{\min},k)| = k$.
Alternatively, the estimator's computational cost (run-time or memory) can be encoded in the class sizes. 
More generally, we can consider $|\mathcal{C}^{\text{GD}}(\eta,k_{GD})| = k_{\text{GD}}$ and $|\mathcal{C}^{\text{Ridge}}(\lambda_{\min},k_{\text{Ridge}})| = k_{\text{Ridge}}$
and provide bounds 
in terms of relative size $c = k_{\text{GD}}/k_{\text{Ridge}}$. When the eigenvalues decaying slowly, the relative sub-optimality will then be scaled by a multiplicative factor \footnote{
Looking to the dominant second term in  equation \eqref{equ:ratio:sketch}, with $\eta = 1/(k_{\text{GD}}\lambda_{\star})$ and $\delta = O(1/k_{\text{Ridge}})$ we note it is scaled by $\eta^2/\delta^2 = 1/(c^2 \lambda_{\star}^2)$.  
} of $1/c^2$ in this case. As gradient descent is typically cheaper to compute than ridge regression, we are lead to consider $c \geq 1$.
Thus, accounting for computational cost typically only further improves the performance of gradient descent relative to ridge regression.  

\textbf{Supplementary Results and Discussion}
Additional results and discussion are included within Section \ref{sec:app:additional_results_discussion} of the supplementary material. Precisely, Theorem  \ref{thm:log-ridgeVsGD} in Section \ref{sec:app:log_grid} demonstrates gradient descent still outperforms ridge regression tuned over a logarithmic grid when eigenvalues decay slowly, mirroring Theorem \ref{thm:Informal}.  Section \ref{sec:app:add_discussion} provides discussion related to the framework investigated within this work as well as directions for future research, and 
Section \ref{sec:ProofSketch} contains a proof sketch for Theorem \ref{thm:Informal}.

The following subsection, Section \ref{sec:experiments}, provides an empirical investigation into the upper and lower bounds from Theorem \ref{thm:Informal} alongside an extension to the random design setting.  
The following sections then provide more detailed results for the relative sub-optimality, for various assumptions on the eigenvalues. 
Section \ref{gd-vs-ridge-slowdecay} provides the results for a slow power law eigenvalue decay, while Section  \ref{gd-vs-ridge-expdecay} provides the results for exponential and fast power law decays. Finally, Section \ref{b-ind} provides bounds for individual classes of estimators.

\subsection{Experiments}
\label{sec:experiments}
We study how the relative sub-optimality 
$\mathcal{S}(\mathcal{C}^{\text{GD}}(\eta,k),\mathcal{C}^{\text{Ridge}}(\lambda_{\min},k))$
between gradient descent and ridge regression  scales with the noise standard deviation $\sigma$ in the \emph{slow power law decay} regime. To do so, we consider two experiment settings. The first evaluates the upper and lower bounds on the relative sub-optimality $\mathcal{S}\left( \mathcal{C}^{\text{GD}}(\eta,k), \mathcal{C}^{\text{Ridge}}(\lambda_{\min},k)\right)$ when $s_i = i^{-\alpha}$ for $i=1,\dots,r$. The second considers the randomness introduced from generating data according to \eqref{equ:Response} and \eqref{equ:NoiseSignal:ass} with the eigenvalue decay at the population level, specifically when each row in the design matrix $X$ is drawn from a multivariate normal distribution with a population covariance matrix $\text{Diag}(s_1,s_2,\dots,s_d)$ where $s_i = i^{-\alpha}$. This slightly differs from our theoretical analysis, which considers a \emph{fixed design} setting; and the goal is to evaluate our results in a somewhat broader range of settings. Each of these experiments is presented in its own paragraph. 

\begin{figure}[!h]
    \centering
    \includegraphics[width=0.5\textwidth]{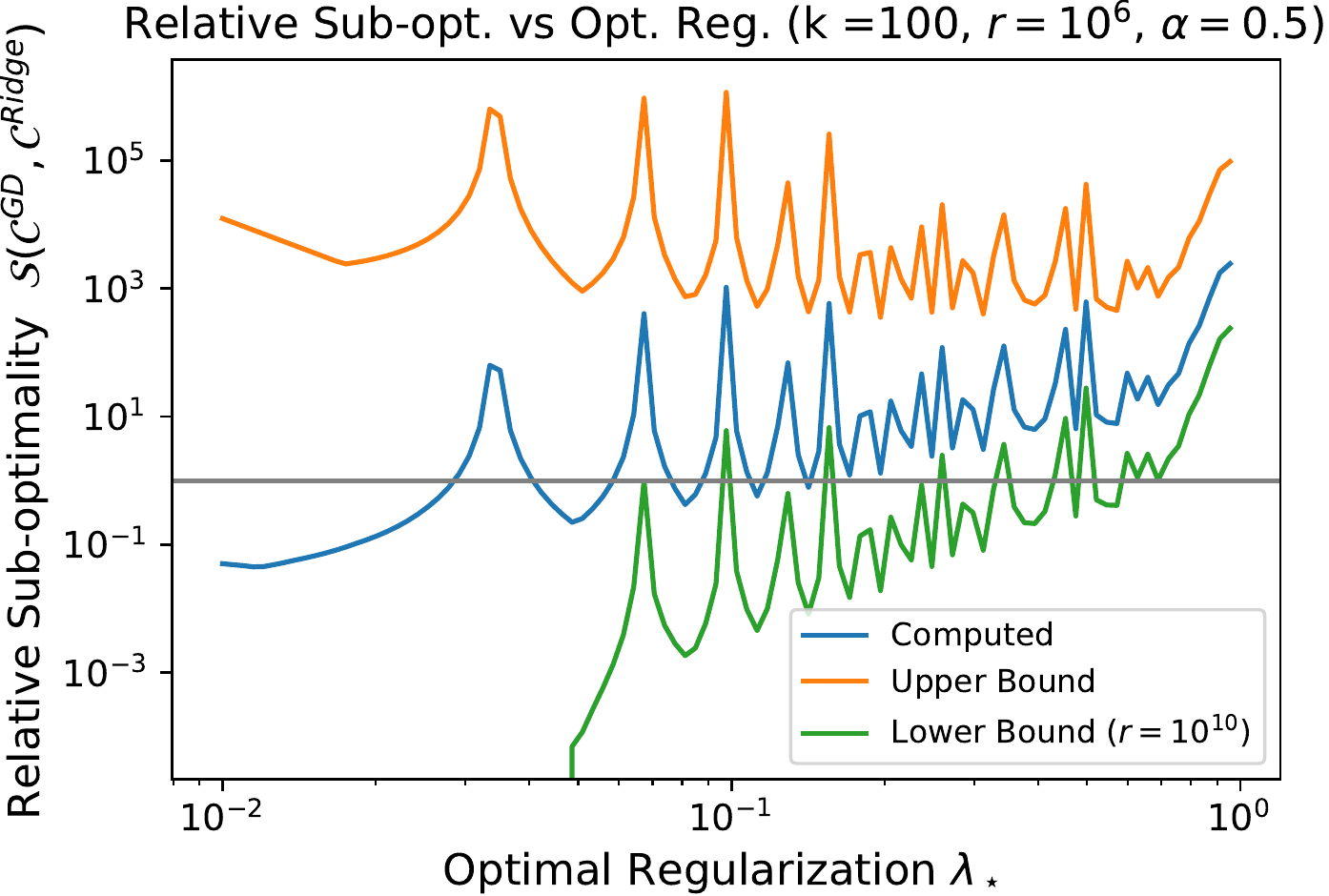}
    \caption{Plot of relative sub-optimality $\mathcal{S}(\mathcal{C}^{\text{GD}}(\eta,k),\mathcal{C}^{\text{Ridge}}(\lambda_{\min},k))$ against the optimal regularization $\lambda_\star$ on a $\log-\log$ scale.
    Problem parameters: $\lambda_{\min} = 0.71245212 \cdot \lambda_\star$ (the factor of $0.71245212$ ensures $1/\eta \lambda_{\min}$ is not an integer), $k=100$, $\eta = 1/(k \lambda_{\min})$, $\alpha = 0.5$, $r=10^6$, $s_i = i^{-\alpha}$ for $i=1,\dots,r$. Also plotted for 100 values of $\lambda^{\star} \in [10^{-2},1]$ is the upper bound from Theorem \ref{thm:mainGDvsRidge_Upper} with constant equal to unity instead of $1152$, and the lower bound from Theorem  \ref{thm:mainGDvsRidge_Lower} with constant equal to unity instead of $1/512$ and with $r=10^{10}$. 
    }
    \label{fig:GDVsRidge:upperlower}
\end{figure}

\textbf{Numerical Evaluation of Upper and Lower Bounds.}
Looking to Figure \ref{fig:GDVsRidge:upperlower}, 
we see the relative sub-optimality plotted against the optimal amount of regularization, alongside the theoretical upper and lower bounds from Theorems \ref{thm:mainGDvsRidge_Upper} and \ref{thm:mainGDvsRidge_Lower},  (these are more precise versions of Theorem \ref{thm:Informal}), and we are using constants of unity instead of $1152$ and $1/512$, respectively). 

We observe that the numerically computed relative sub-optimality can indeed be much smaller than unity when $\lambda_\star$ is small; in which case the grid of gradient descent estimators performs much better than ridge regression.
This supports our theoretical observations about the two estimators.
The upper bound captures the variation in the relative sub-optimality arising from the discretization. 

Currently, when the rank $r = 10^{6}$,  a gap remains between the upper and the lower bounds. The general qualitative behavior of the bounds appears correct, but the constants could be tightened with additional work. Tightening the dependence on the rank $r$ is left for future work.  

\begin{figure}[!h]
    \centering
    \includegraphics[width=0.5\textwidth]{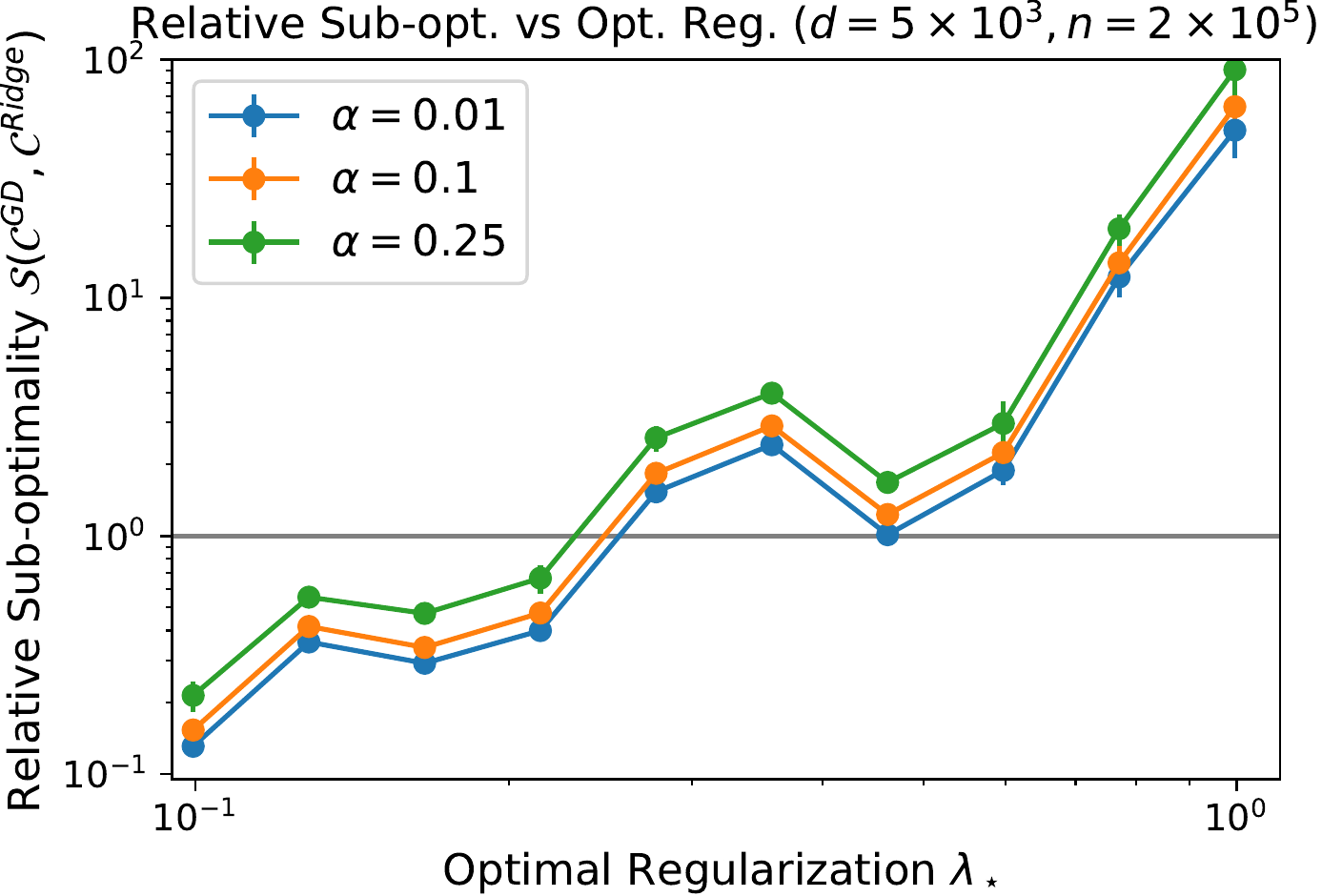}
    \caption{Estimate of relative sub-optimality $\mathcal{S}(\mathcal{C}^{\text{GD}}(\eta,k),\mathcal{C}^{\text{Ridge}}(\lambda_{\min},k))$ as a function of optimal regularization $\lambda_\star$ on a $\log-\log$ scale. Problem parameters: $\eta = 1$, $\sigma \in [10^{0.8},10^{1.3}]$, $\psi = 10$, $\alpha \in \{0.01,0.1,0.25\}$, $d = 5 \cdot 10^{3}$, $n=2\cdot 10^{5}$, $k = \lceil 1/(\eta \lambda_{\min})\rceil = 10$ where $\lambda_{\min} = 0.99\cdot 10^{1.6} d/(n \psi)$.  
    Fixed $\beta^{\star} = \mathbf{1} \sqrt{\psi/d}$ where $\mathbf{1}$ is a vector of ones. 
    For each estimator, the expected sub-optimality is estimated as an average over 20 replications, with the best estimator chosen thereafter. For gradient descent, we choose the best estimator from the class $\{\hat\beta_{\eta,s}\}_{s=1,\dots,k}$, for ridge regression we choose the best regularization out of $\{\lambda_\star + (1-\lambda_{\min})/(k-1),\ldots,\lambda_\star - (1-\lambda_{\min})/(k-1)\}$. 
    The error bars (short and barely visible) show two standard deviations over ten replications of this process.}
    \label{fig:GDVsRidge}
\end{figure}

\textbf{Random Design Setting.} Looking to Figure \ref{fig:GDVsRidge} an estimate of the relative sub-optimality is plotted against the optimal regularization $\lambda_\star = \sigma^2 d/(\psi n)$, for a particular choice of $n,d,k$.\footnote{In this experiment $\psi\neq 1$, but this is immaterial as it only amounts to a scaling of the results.} 
We observe a roughly linear relationship (on a $\log-\log$ scale), with the relative sub-optimality decreasing as the optimal regularization decreases. This is consistent with the theoretical results in Theorem \ref{thm:Informal}. The estimated relative sub-optimality does not decrease monotonically in $\lambda_\star$. A possible explanation is that the discretization grid is fixed, and thus, as the optimal regularization parameter varies, so does its distance to the closest grid point. Moreover, the relative sub-optimality is increasing in $\alpha$, which is consistent with our theoretical results for this range of $\alpha$.

\subsection{Eigenvalues with a Slow  Power Law Decay}
\label{gd-vs-ridge-slowdecay}
In this section we consider eigenvalues decaying according to a slow power law, i.e., $s_i = i^{-\alpha}$ for $\alpha \in (0,1)$. 
In this case we show, when the rank $r$ is sufficiently large, that the relative sub-optimality of gradient descent and ridge regression is on the order of $\Theta(\lambda_\star^4/\lambda_{\min}^2)$. This is summarized in the following two theorems, the first providing the upper bound and the second a lower bound. 
\begin{theorem}
\label{thm:mainGDvsRidge_Upper}
Let $2/k < \lambda_{\min} \leq \lambda_{\star} \leq 1$, 
$\delta = (1-\lambda_{\min})/(k-1)$, 
$\eta = 1/(k \lambda_{\min})$ and $s_i = i^{-\alpha}$ for $i=1,\dots,r$, with $\alpha \in (0,1)$. 
If $r \geq 2^{1/(1-\alpha)}(1+\lambda_\star^{-1/\alpha})$ and  
$\lambda_{\star} = \lambda_{\min} + (j+\ep)\delta $ for some $j \in \{0,\dots,k-2\}$ and $\ep \in [0,1]$,\footnote{When $\ep\in\{0,1\}$, the denominator equals zero; and in that case we use the convention that $1/0 = +\infty$.} 

\begin{align*}
    \mathcal{S}\left( \mathcal{C}^{\text{GD}}(\eta,k), 
    \mathcal{C}^{\text{Ridge}}(\lambda_{\min},k)
    \right)
    & \leq 
    \frac{1152}{\min\{1-\ep,\ep\}^2(1 - \lambda_{\min})^2}
    \cdot\left[ 
    A(k,r,\alpha,\lambda_\star,\lambda_{\min}) 
    + 
    \frac{\lambda_{\star}^4}{\lambda_{\min}^2}
    \right]
\end{align*}
where $A(k,r,\alpha,\lambda_\star,\lambda_{\min}) = \frac{k^2}{r^{1-\alpha}}
    \left( \frac{1}{\lambda_\star^{1/\alpha - 1}} 
    + \mathcal{J}_{\alpha,\lambda_{\star}}(k,\lambda_{\min})
    \right)$
and
\begin{align}\label{J}
    \mathcal{J}_{\alpha,\lambda_{\star}}(k,\lambda_{\min}) 
    = 
    1 + \begin{cases}
     \frac{\lambda_\star^{3-1/\alpha}}{3\alpha -1} 
    & \text{if  } \alpha > 1/3 \\
    3 \log(k \lambda_{\min} /\lambda_\star)
    & \text{ if } \alpha = 1/3 \\
     \frac{1}{1-3\alpha} \left( \frac{k \lambda_{\min}}{\lambda_{\star}^2} \right)^{1/\alpha - 3} 
    & \text{ if } \alpha < 1/3. 
    \end{cases}
\end{align}
\end{theorem}
The upper bound in Theorem \ref{thm:mainGDvsRidge_Upper} goes to infinity as $\ep \rightarrow 0$ or $\ep\to1$. 
We recall that this is due to $\lambda_{\star}$ converging to a grid point in the discretization of ridge regression. A similar phenomenon occurs as $\lambda_{\min} \rightarrow 1$, because $\lambda_{\star} \rightarrow 1$ converges to a grid point.   The bound also diverges as $\lambda_{\min} \rightarrow 0$, because gradient descent performs many wasteful iterations i.e., $k = 1/(\eta \lambda_{\min}) \rightarrow \infty$. Specifically, gradient descent gets compared to ridge regression with $\delta \rightarrow 0$, which converges to the Bayes optimal estimator, so the denominator in $\mathcal{S}\left( \mathcal{C}^{\text{GD}}(\eta,k),  \mathcal{C}^{\text{Ridge}}(\lambda_{\min},k) \right)$ goes to zero.

The bound in Theorem \ref{thm:mainGDvsRidge_Upper} consists of three terms.  The first two decrease with the number of eigenvalues $r$, with the third depending upon the factor $\lambda_{\star}^{4}/\lambda_{\min}^2$. Precisely, defining $r_{\lambda_\star,\alpha,k,\lambda_{\min}}^{\text{Upper}}$ as
\begin{align*}
    r_{\lambda_\star,\alpha,k,\lambda_{\min}}^{\text{Upper}} : = 
    \left(\frac{\lambda_{\min} k }{\lambda_\star^2}\right)^{2/(1-\alpha)} 
    \left( \frac{1}{\lambda_\star^{1/\alpha - 1}} 
+ \mathcal{J}_{\alpha,\lambda_{\star}}(k,\lambda_{\min})
\right)^{1/(1-\alpha)},
\end{align*}
then the bound is $O(\lambda_\star^4/\lambda_{\min}^2)$ provided the rank is sufficiently large that $r \geq r_{\lambda_\star,\alpha,k,\lambda_{\min}}^{\text{Upper}}$. 
Each term within the bound represents the ratio between the error of gradient descent and ridge regression in different parts of the spectrum (for a precise definition see Section \ref{b-ind}). The third term is associated to eigenvalues below the regularization level, i.e.,  $s_i < \lambda_{\star}$, and becomes the ``dominant'' source of error when the eigenvalues decay slowly and the rank is sufficiently large, i.e., $r \geq r_{\lambda_\star,\alpha,k,\lambda_{\min}}^{\text{Upper}}$. Gradient descent then has a smaller error on these eigenvalues than ridge regression by a factor of  $O(\lambda_{\star}^2/\lambda_{\min})$. For details on how this factor arises  see the discussion in Section \ref{b-ind-GD}.

We give a lower bound which shows that the dependence on $\lambda_{\star}$ and $\lambda_{\min}$ is sharp as the number $r$  of non-zero eigenvalues grows. 
\begin{theorem}
\label{thm:mainGDvsRidge_Lower}
Consider the setting of Theorem \ref{thm:mainGDvsRidge_Upper}.  
Moreover, suppose that $1/(\eta \lambda_{\star})= \ell + \kappa$ for some integer $\ell \in \{0,1,2,\dots,k-1\}$ and $\kappa \in (0,1)$. Then 
\begin{align*}
    \mathcal{S}\! \left( \mathcal{C}^{\text{GD}}(\eta,\!k), 
    \mathcal{C}^{\text{Ridge}}(\lambda_{\min},\! k)
    \right) \! \geq \! 
    \frac{\min\{ 1- \kappa,\kappa \}^2}{512  \min\{ 1 \! -\! \ep, \ep \}^2 (1 \! -  \! \lambda_{\min})^2}
    \!\!
    \left( 
    \!\! 1 \! - \! 
    \frac{2\left(4 k \lambda_{\min}/(\lambda_{\star}^2 \kappa)\right)^{\frac1\alpha - 1}}{(1-\alpha)r^{1-\alpha}} 
    \!
    \right)
    \! \!
    \frac{\lambda_{\star}^4}{\lambda_{\min}^2}.
\end{align*}
\end{theorem}

Theorem \ref{thm:mainGDvsRidge_Lower} implies, when the rank $r$ is large, the ratio of errors is lower bounded by a factor of $\lambda_{\star}^4/\lambda_{\min}^2$. 
Thus, Theorem \ref{thm:mainGDvsRidge_Upper} is sharp in the large $r$ regime. Precisely,  define $r^{\text{Lower}}_{\lambda_{\star},\alpha,k,\lambda_{\min}}$ as
\begin{align*}
    r^{\text{Lower}}_{\lambda_{\star},\alpha,k,\lambda_{\min}} := 
    \left(\frac{4}{1-\alpha}\right)^{1/(1-\alpha)} \left(\frac{4 k \lambda_{\min}}{\lambda_{\star}^2 \kappa }\right)^{1/\alpha}.
\end{align*}
If $r \geq r^{\text{Lower}}_{\lambda_{\star},\alpha,k,\lambda_{\min}}$, the lower bound is on the order of $\lambda_\star^{4}/\lambda_{\min}^2$. When considering Theorem \ref{thm:Informal} we see that it is then sufficient to pick $r_{\lambda_{\star},\alpha,k,\lambda_{\min}} = \max\{r^{\text{Lower}}_{\lambda_{\star},\alpha,k,\lambda_{\min}},r^{\text{Upper}}_{\lambda_{\star},\alpha,k,\lambda_{\min}}\} $, for the upper and lower bounds to be of the same order. 

We note that the
 lower bound goes to zero as $\kappa \rightarrow 1$ or $\kappa \rightarrow0$, due to the proof approach. In particular, we leverage that gradient descent must perform a discrete number of iterations, i.e., $\lfloor 1/(\eta \lambda_{\star})\rfloor$ or $\lceil 1/(\eta \lambda_{\star})\rceil$, each of which, respectively,  under- or overshrinks smaller eigenvalues ($s \leq \lambda_{\star})$. For more details see proof in Appendix \ref{sec:proof:GDBounds:Fine}.  
This is reasonable as we expect the upper bound to also decrease in this setting.

Following Theorem \ref{thm:mainGDvsRidge_Upper} and \ref{thm:mainGDvsRidge_Lower}, we now consider a low-dimensional setting where the ambient dimension $d =n^{q}$ for some $q \in (0,1)$. In this case some care is required when interpreting the bounds, as the optimal regularization $\lambda_{\star}  = O(n^{-(1-q)})$ and rank $r =O( n^{q})$ are coupled to the sample size $n$ and dimension $d=n^q$.     
\begin{corollary}
\label{cor:LowDimensional}
Consider the setting of Theorem \ref{thm:mainGDvsRidge_Upper}. Suppose $r = d = n^{q}$ for some $q \in (1/(\alpha+1),1]$, $\lambda_{\min} = \sigma_{\min}^2 d/n$ for $0 < \sigma_{\min} < \sigma/\sqrt{2} $ where $n \geq \max\{ 4^{1/(q(1-\alpha))}, (4^{1/(1-\alpha)} \sigma^{-2/\alpha})^{\alpha/(q(1+\alpha) - 1)}, \sigma_{\min}^{2/(1-q)} \}$ and $k > 2n^{1-q}/\sigma_{\min}^2$. Then, with $\tau_n=\frac{(\sigma^2/\sigma_{\min})^4}{n^{2(1-q)}} $ and
$B(\sigma,\alpha,k,\sigma_{\min}, n) =
    \sigma^{-2(\alpha^{-1} -1)} 
    +
    n^{2(q-1)} (k \sigma_{\min}^2/\sigma^4)^{1/\alpha - 3}$,
    \begin{align}
     \mathcal{S}\left( \mathcal{C}^{\text{GD}}(\eta,k), 
    \mathcal{C}^{\text{Ridge}}(\lambda_{\min},k)
    \right) 
     & \lesssim \!
    \frac{1}{\mathrm{Dist}_{\delta}(\lambda_{\star},\Gamma)^2}
    \Bigg[\!
    \frac{ k^2 \! \log\left(\! \frac{k \sigma_{\min}^2}{\sigma^2} \! \right)}{n^{q(\alpha^{-1}\! - \! \alpha) \! - \! \alpha^{-1} + 1} }
    B(\sigma,\alpha,k,\sigma_{\min}, n)
    \! + \!
    \tau_n
    \! \Bigg].\label{ubcor}
    \end{align}
Moreover under the assumptions Theorem \ref{thm:mainGDvsRidge_Lower} the following lower bound holds  
\begin{align}\label{lbcor}
    \mathcal{S}\left( \mathcal{C}^{\text{GD}}(\eta,k), 
    \mathcal{C}^{\text{Ridge}}(\lambda_{\min},k)
    \right) \gtrsim 
    \frac{1}{\mathrm{Dist}_{\delta}(\lambda_{\star},\Gamma)^2}
    \left( 
    1 - \frac{(4 k \sigma_{\min}^2/(\sigma^4 \kappa) )^{1/\alpha - 1}}{n^{q(1/\alpha - \alpha) - 1/\alpha +1}}
    \right)
   \tau_n.
\end{align}
\end{corollary}

As the optimal amount of regularization $\lambda_{\star}$ decreases with the sample size, the upper bound also decreases with the sample size. 
When $\sigma,\sigma_{\min},k$ are fixed, the upper bound decreases---up to logarithmic factors---at an $O(n^{- \min\{ q(1/\alpha-\alpha) - 1/\alpha + 1, 2(1-q) \} })$ rate, with $q > (\alpha \! +\!1)/(1 \! + \! \alpha(2 \! - \! \alpha)) $ yielding a rate matching the lower bound \eqref{lbcor} of $\Theta(n^{-2(1-q)})$. One interpretation of this is that the rates are optimal provided the rank or dimension, as controlled by $q$, is sufficiently large.

\subsection{Eigenvalues Decaying Exponentially or with a Fast Power Law }
\label{gd-vs-ridge-expdecay}
In this section we consider eigenvalues decaying at either an exponential rate i.e., for all $i=1,\ldots,r$, $s_i = \exp(-\rho(i-1))$, or a fast power law rate $s_i = i^{-\alpha}$ for $\alpha > 1$. In this case we show that ridge regression outperforms gradient descent as the number $k$ of models increases, regardless of the rank $r$. This is summarized within the following two theorems, the first considering exponentially decaying eigenvalues and the second considering eigenvalues that decay at a fast power law rate. 
\begin{theorem}
\label{thm:mainGDvsRidge_Upper:exp}
Let $\delta = (1-\lambda_{\min})/(k-1)$, $\eta = 1/( k \lambda_{\min})$ and $s_i = \exp(-\rho (i-1))$ for $i =1,2,\dots,r$ with $\rho > 0$, for $r \geq 1 + {\rho}^{-1} \log(1/\lambda_{\star})$.  Furthermore, suppose that $\lambda_{\star} = \lambda_{\min} + (j+\ep)\delta $ for some $j \in \{0,\dots,k-2\}$ and $\ep \in [0,1]$. If 
\begin{align}\label{expgdc}
\lambda_{\star} 
    < 
    \frac{1}{2} \left(2e + \frac{4 e^{1+\rho}}{1 - e^{-\frac{1}{64}}}\right)^{-1}
    \quad \text{ and }
    \quad
    k \geq \frac{ 1 }{\lambda_{\min}} \max\left\{ 64 \log\left(1 + \frac{1}{\lambda_{\star}}\right),\frac{8e}{\lambda_{\star}} \frac{1 + e^{-\rho}}{1-e^{-{64}^{-1}}}\right\},
\end{align}
then
\begin{align*}
    \mathcal{S}\left( \mathcal{C}^{\text{GD}}(\eta,k), 
    \mathcal{C}^{\text{Ridge}}(\lambda_{\min},k)
    \right)
    & \geq  
    \frac{C_{\rho}}{\min\{1-\ep,\ep\}^2}
    \left( \frac{ \lambda_{\star} k}{1-\lambda_{\min}}\right)^2,
\end{align*}
where $C_{\rho} = (1 - \exp(-{64}^{-1}))^3/[512e^2(1+e^{2\rho})(1 - e^{-{64}^{-1}} + 2 e^{\rho})^2] $.
\end{theorem}
\begin{theorem}
\label{thm:mainGDvsRidge_Upper:poly}
Suppose $\delta = (1-\lambda_{\min})/(k-1)$, $\eta = 1/( k \lambda_{\min})$ and $s_i = i^{-\alpha}$ for $i =1,2,\dots,r$ and $\alpha > 1$, for $r \geq \lambda_{\star}^{-1/\alpha}$.  Furthermore, suppose that $\lambda_{\star} = \lambda_{\min} + (j+\ep)\delta $ for some $j \in \{0,\dots,k-2\}$ and $\ep \in [0,1]$. If 
\begin{align}
    \lambda_{\star} 
    & < 
    \frac{1}{2^{1+2\alpha}} \left[ 1 + 2^{4+3\alpha}(1-e^{-\frac{1}{64}})^{-1}\right]^{\frac{-\alpha}{\alpha+1}} 
    \label{thm6cond1}
    \\
    k & \geq \frac{ 1 }{\lambda_{\min}} \max\left\{ 64 \log\left(1 + \frac{1}{\lambda_{\star}}\right),\frac{17\cdot 2^{\alpha}}{\lambda_{\star}^{1-1/\alpha}}
    \frac{\alpha+1}{\alpha - 1} (1-e^{-{64}^{-1}})^{-1}\right\},\label{thm6cond2}
\end{align}
then
\begin{align*}
    \mathcal{S}\left( \mathcal{C}^{\text{GD}}(\eta,k), 
    \mathcal{C}^{\text{Ridge}}(\lambda_{\min},k)
    \right)
    & \geq  
    \frac{C_{\alpha}}{\min\{1-\ep,\ep\}^2}
    \left( \frac{ \lambda_{\star} k}{1-\lambda_{\min}}\right)^2,
\end{align*}
where $C_{\alpha} = \frac{\alpha-1}{\alpha}(1 - \exp(-{32}^{-1}))^2/[32 (2^{(1+2\alpha)(2+\alpha(1+\alpha))}) (1+2\alpha)(1 + 2^{4 + 3\alpha} (1- e^{-\frac{1}{64}}))^{\frac{1+2\alpha}{\alpha+1}}] $.
\end{theorem}

The theorems state,  provided the optimal regularization is smaller than the specified constants, and the number of models $k$ is large with respect to $1/(\lambda_{\star}\lambda_{\min})$, 
that the relative sub-optimality is lower bounded by a factor independent of $k,\lambda_{\star}$ 
times $(\lambda_{\star}k)^2$. 
Thus, as $k$ grows, ridge regression outperforms gradient descent, since the ratio increases. 
This behavior aligns with what we have when we consider gradient descent estimators with all possible regularization parameters.

The condition on the number of models $k$ differs from the previous results, in that we now require $k \gtrsim (\lambda_{\min} \lambda_{\star})^{-1} \gtrsim \lambda_{\star}^{-2}$ whereas Theorem \ref{thm:mainGDvsRidge_Upper} and \ref{thm:mainGDvsRidge_Lower} required $k \gtrsim \lambda_{\min}^{-1}$. This is to ensure a sufficiently small step size $\eta$ which then ensures certain technical conditions needed in the proof.
We leave refining the lower bound to values of $k$ within the range $(\lambda_{\star}\lambda_{\min})^{-1} \gtrsim k \gtrsim \lambda_{\star}^{-1}$ to future work.  
We note that $C_{\rho},C_{\alpha}$ decay exponentially in $\rho$ and $\alpha$ respectively, and thus, the lower bounds remain non-vacuous as $k$ grows, provided $\rho \lesssim \log(k)$ and $\alpha \lesssim \log(k)$.

\subsection{Bounds for Individual Classes}
\label{b-ind}
In this section we present bounds for the class of estimators associated to gradient descent and ridge regression.
We denote the error for a class of estimators $\mathcal{C}$ in terms of the gap from optimally tuned ridge regression (contrast with \eqref{subopt})
\begin{align*}
    \mathcal{E}(\mathcal{C}) 
    := 
    \min_{\hbeta \in \mathcal{C}} \left\{ \E_{\beta_{\star},\epsilon}[L_{\beta_{\star}}(\hbeta)]
    - 
    \E_{\beta_{\star},\epsilon}[L_{\beta_{\star}}(\hbeta_{\lambda_{\star}})]\right\}.
\end{align*}
We define the spectral measure $\widehat{H}$ associated to the eigenvalues $\{s_i\}_{i=1}^{r}$, 
to have a cumulative distribution function, also denoted as $\widehat{H}$, for all $s\in\mathbb{R}$, 
given by
$\widehat{H}(s) := \frac{1}{d} \sum_{i=1}^{r} \mathbbm{1}(s_i \leq s)  + (1-\frac{r}{d})\mathbbm{1}(0 \leq s)$ where $\mathbbm{1}(\cdot)$ is the indicator function. 
The following results bound the sub-optimality for a given class  $\mathcal{E}(\mathcal{C})$ in terms of a function integrated against the the spectral measure $\widehat{H}$. 
This approach is rather natural, as the estimation error of optimally tuned ridge regression  can be written
\begin{align}
\label{equ:OptTunedRidge}
\E_{\beta_{\star},\epsilon}[L_{\beta_{\star}}(\hbeta_{\lambda_{\star}})] 
= 
\lambda_{\star} \int \frac{1}{s + \lambda_{\star}} d \widehat{H}(s).
\end{align}
In our case bounds on  $\mathcal{E}(\mathcal{C})$ take the form $\lambda_{\star} \int G(s) d \widehat{H}(s)$ for an ``error'' function $G : \mathbb{R}_{+} \rightarrow \mathbb{R}_{+}$ that describes ``where'' on the spectrum the estimator is sub-optimal. 
Precisely, the function $G$ will depend upon $\lambda_{\star}$ and encode properties of the class $\mathcal{C}$, e.g., the step size $\eta$ and number of iterations $t$ for gradient descent, or the discretization length $\delta$ for ridge regression. Since $\widehat{H}(s)$ is an  atomic measure, we adopt the convention, for any $b > a$, that $\int_b^{a} d \widehat{H}(s) = \int \mathbbm{1}(s)_{(a,b]}  d \widehat{H}(s) $ i.e., the integration range is closed on the right and open on the left.

The remainder of this section is then structured as follows. Section \ref{b-ind-RR} considers  classes arising from ridge regression with a uniformly discretized regularization parameter. Section \ref{b-ind-GD} considers classes arising from the iterations of gradient descent.

\subsubsection{Ridge Regression}
\label{b-ind-RR}
The following lemma provides upper and lower bounds on the class of estimators arising from ridge regression in the case of a \textit{fine discretization}, i.e., where $\delta = (1-\lambda_{\min})/(k-1) \leq \lambda_{\star}$.

\begin{proposition}[Ridge Regression Estimator, Fine Discretization]
\label{lem:RidgeBounds}
Let $\delta = (1-\lambda_{\min})/(k-1)$. 
Suppose $0 \leq \lambda_{\min} \leq \lambda_{\star} \leq 1$, $\lambda_\star \geq \delta > 0$ and $\lambda_{\star} = \lambda_{\min} + (j + \ep)\delta$ for some $j \in \{0,1,\dots,k-2\}$ and $\ep \in [0,1]$. Then 
\begin{align*}
    \frac{1}{2^5}\!
    \min\! \{1 \! - \! \ep,\ep\}^2
    \lambda_{\star}
    \!\!\! \int \! \! G^{\text{Ridge}}(s) d \widehat{H}(s) \! \leq \mathcal{E}(\mathcal{C}^{\text{Ridge}}(\lambda_{\min},k))
    \leq \! 4 \min\{ 1 \! - \! \ep, \ep \}^2 
    \lambda_{\star} \!\!\! \int \!\! G^{\text{Ridge}}(s) d \widehat{H}(s)
\end{align*}
where 
\begin{align*}
    G^{\text{Ridge}}(s) 
    = 
    \delta^2
    \begin{cases}
     \frac{1}{\lambda_{\star} s^2}   &  \text{ if } s \in (\lambda_{\star},\infty) \\
     \frac{s}{\lambda_{\star}^4} & \text{ if } s \in (0,\lambda_{\star}].
    \end{cases}
\end{align*}
\end{proposition}
Results for tuning ridge regression with either a coarse discretization $\delta > \lambda_{\star}$ or logarithmic grid are presented within Proposition \ref{lem:RidgeBounds:Coarse} and \ref{lem:log_grid} in Section \ref{sec:app:add_bounds_individual:RR}.

Now, let us rewrite the sub-optimality of ridge regression in terms of a multiplicative constant of optimally tuned ridge regression, so that for some constant $C > 0$, and $G^{\text{Ridge}}$ from Proposition \ref{lem:RidgeBounds},
\begin{align*}
    \min_{\hbeta \in \mathcal{C}^{\text{Ridge}}(\lambda_{\min},k) } 
    \E_{\beta_{\star},\epsilon}[L_{\beta_{\star}}(\hbeta)] 
    = \lambda_{\star} \int \Big( 1 + C G^{\text{Ridge}}(s)(s+\lambda_{\star}) \Big)  \frac{1}{s + \lambda_{\star}} d \widehat{H}(s).
\end{align*}
We can interpret $C G^{\text{Ridge}}(s)(s+\lambda_{\star})$ as a form of spectral dependent \emph{risk inflation}.  Recalling the discretization length $\delta = 1/k$, and using Propositions \ref{lem:RidgeBounds}, \ref{lem:RidgeBounds:Coarse} and \ref{lem:log_grid}, we can take a uniform bound to get the risk inflation of ridge regression with each discretization
\begin{align}
    \max_{s \in (0,s_1] } 
    G^{\text{Ridge}}(s)(s+\lambda_{\star}) 
    & = \frac{2}{k^2} \frac{1}{\lambda_{\star}^2}, 
    \qquad
    \max_{s \in (0,s_1] } 
    G^{\text{Coarse Ridge}}(s)(s+\lambda_{\star}) 
     = \frac{1}{k \lambda_{\star}}    \nonumber
    \\
    \max_{s \in (0,s_1] } G^{\text{log-ridge}}(s) (s+\lambda_{\star}) 
    & = \frac{\log^2(1/\lambda_{\min})}{4 (k-1)^2}.
    \label{equ:uniform_vs_log_grid}
\end{align}
We make two observations. 
First, for the uniform discretizations, the risk inflation initially decays at the ``slow'' rate of $(k \lambda_{\star})^{-1}$ for $k \leq \lambda_{\star}^{-1}$, and then moves to the fast   $(k \lambda_{\star})^{-2}$ rate for $k \geq \lambda_{\star}^{-1}$. This phase transition in the number of models $k$ also reflects what is observed in Section \ref{sec:SuboptCompGrid} for ridge regression with an orthogonal design. 
Second, as the log-grid relies on knowing a lower bound $\lambda_{\min} < \lambda_{\star}$ we see trade-off as to when a log-grid is beneficial over a uniform grid. Precisely, as the log-grid is favored approximately when $\log(1/\lambda_{\min}) \leq 1/\lambda_{\star}$ or $\lambda_{\min} \geq e^{-1/\lambda_{\star}}$.

\subsubsection{Gradient Descent}
\label{b-ind-GD}
The following presents the upper and lower bounds for gradient descent. Recall the definition of $\mathcal{C}^{\text{GD}}(\eta,t)$ provided at the beginning of Section \ref{gd-vs-ridge}. 

\begin{proposition}[Gradient Descent Estimators, Fine Discretization]
\label{lem:GDBounds:Fine}
Suppose $0 <  \lambda_{\star} \leq 1$ and $t \geq  \lceil1/(\eta \lambda_{\star})\rceil $ for a step size $0 < \eta \leq  1/\max\{s_1,\lambda_{\star}\}$. Then 
\begin{align*}
    \mathcal{E}(\mathcal{C}^{\text{GD}}(\eta,t)) 
    \leq 
    18 \lambda_{\star}\int G^{\text{GD}}(s) d \widehat{H}(s) 
\end{align*}
where
\begin{align}\label{ggd}
    G^{\text{GD}}(s)
    = 
    \begin{cases}
    \frac{\lambda_{\star}}{s^2} & \text{ for } s \in (\lambda_{\star},1/\eta] \\
    \frac{s^3}{\lambda_{\star}^4} & \text{ for } s \in (\eta \lambda_{\star}^2,\lambda_{\star}] \\
    \eta^2  s & \text{ for }  s \in (0,\eta \lambda_{\star}^2].
    \end{cases}
    \end{align}
Furthermore, if $\frac{1}{\eta \lambda_{\star}} = \ell  + \kappa$ for some integer $\ell$ and $\kappa \in (0,1)$, then for any $t\ge 1$
\begin{align*}
    \mathcal{E}(\mathcal{C}^{\text{GD}}(\eta,t)) 
    \geq 
    \frac{1}{16} 
    \min\{1-\kappa,\kappa\}^2
    \lambda_{\star}
    \int_{0}^{\eta \lambda_{\star}^2 \kappa/4} G^{\text{GD}}(s) d \widehat{H}(s).
\end{align*}
\end{proposition}

For ridge regression 
as studied in Proposition \ref{lem:RidgeBounds}, the discretization length tuning parameter $\delta$ simply scales the error function.
In contrast, for gradient descent, 
the stepsize $\eta$ does not simply scale the error function $G^{\text{GD}}$. 
One interpretation is that ridge regression converges \emph{uniformly} as the number of models increases. However, for gradient descent a non-zero error remains regardless of how small the step size is (see Theorem \ref{thm:GDLower2} in Appendix \ref{sec:proof:Subopt:FastDecay}). 

Following the discussion after Theorem \ref{thm:mainGDvsRidge_Upper} in Section \ref{gd-vs-ridge-slowdecay}, we can see how the factor of $\lambda_\star^4/\lambda^2_{\min}$ arises within the upper bound of the relative sub-optimality. Specifically, consider the step size $\eta = 1/(k \lambda_{\min})$ for gradient descent and discretization length $\delta = 1/k$ for ridge regression. Taking the ratio of error functions, we then have for eigenvalues $s < \eta \lambda_\star^2$ that 
$G^{\text{GD}}(s) /G^{\text{Ridge}}(s) = \lambda_\star^{4}/\lambda_{\min}^2 $. When considering relative sub-optimality, the corresponding ratio of expectations $\int G^{\text{GD}}(s) d\widehat{H}(s) / \int G^{\text{Ridge}}(s) d\widehat{H}(s)$ arises, and therefore, we need to take into account where the spectrum $\widehat{H}(s)$ ``places mass''.

The lower bound in Proposition \ref{lem:GDBounds:Fine} holds over ``part of the spectrum'', in that the integral is over the range up to $\eta \lambda_{\star}^2 \kappa/4$. 
This lower bound is sharp when there are many eigenvalues within the interval $[0,\eta \lambda_{\star}^2 \kappa/4]$, e.g., when the eigenvalues decay at a slow power law rate. For lower bounds when there is a faster decay in the eigenvalues, it is natural to focus on the range $(0,\lambda_{\star}]$. In Theorem \ref{thm:GDLower2} in Appendix \ref{sec:proof:Subopt:FastDecay}, a similar lower bound is given for an integral over the interval $(u_{\min}(\eta),s_1]$, for a certain function $u_{\min}$ of the stepsize, depending on finer notions of the spectrum. 

For completeness a bound for gradient descent in the coarse regime $t \lesssim 1/(\eta \lambda_{\star})$ is presented within Proposition \ref{lem:GDBounds:Coarse} in Section \ref{sec:app:add_bounds_individual:GD} of the supplementary material.

\section{Minimax Estimation under Orthogonal Designs}
\label{mmx-ortho-sec}

For orthogonal designs, we can obtain even more precise results.
We find the \emph{exact minimax optimal} classes of estimators over a range of signal-to-noise ratios. 
Further, we can find the \emph{exact} sub-optimality of ridge regression and gradient descent (Section \ref{sec:SuboptCompGrid}).

To study the minimax optimal class of $k$ estimators, we consider a parameter set $\Theta$ of design matrices $X$, signal strengths $\psi$, noise levels $\sigma$, and a collection $E$ of classes  $T_k  = \{\Phi_1,\ldots,\Phi_k\}$ of $k$ spectral shrinkers  $\Phi_j:[0,\infty)\to [0,\infty)$, $j\in[k]$,
each determining a corresponding class of estimators 
$\widehat{T_k} =\{\hbeta_{\Phi_1},\ldots,\widehat{\beta}_{\Phi_k}\}$. 
\footnote{We will sometimes also call $T_k$ a collection, a set, a $k$-tuple, or---in the scalar case---a grid of estimators.}
Evaluating the \emph{minimax optimal excess risk} over the class  $T_k \in E$ amounts to finding, for $[k] = \{1,\ldots,k\}$, the following \emph{min-max-min} problem:  
\beq\label{mmxrisk}
\inf_{T_k  = \{\Phi_1,\ldots,\Phi_k\}\in E} \,\sup_{\theta\in \Theta} \,\min_{j\in[k]}\,
\mathcal{E}(\widehat{\beta}_{\Phi_j}).
\eeq
Here we are minimizing over classes $T_k\in E$ the worst-case risk over the unknown parameters $\theta\in\Theta$ of the \emph{best estimator $\widehat{\beta}_{\Phi_j}$ in the class $T_k$}. 
This is a generalization of the usual notion of minimax excess risk in statistical decision theory and learning, which corresponds to classes with size $k=1$, where the inner minimization is vacuous
\cite{lehmann1998theory,lehmann2005testing,berger2013statistical}.

Since solving problem \eqref{mmxrisk} could be very challenging in general, we consider parameter spaces $\Theta$ such that all $s_i=s > 0$ are fixed to some known value, and $r=n$. This corresponds to an orthogonal
design where $X^\top X/n = s I_d$. 
We assume that the noise level is fixed, and without loss of generality we consider  $\sigma=1$.
Thus, for some lower and upper bounds $0 \le \psi_- \le \psi_+$ on the signal strength $\psi$, we study the parameter space 
\beq\label{Th}
\Theta = \{\psi \in [\psi_-,\psi_+], 
X^\top X/n = s I_d, n=d\}.
\eeq

We have the following result, proved in Section \ref{sec:mmxgridpf}.

\begin{theorem}[Minimax Optimal Class for Orthogonal Designs]\label{mmxgrid} Let the parameter set $\Theta$ be defined in \eqref{Th} and let $\sigma=1$. 
Consider the collection $E$ of all classes of $k$ estimators 
$T_k = \{\Phi_1,\ldots,\Phi_k\}$ determined by spectral shrinkers.
There is a minimax optimal class, solving the min-max-min problem \eqref{mmxrisk}, of the form
$$\widehat{T_k} = \left\{[1- \Phi_j(s)]/s \cdot X^\top Y/n, j=1,\ldots, k\right\},$$ 
where
$\Phi_j$ are constant functions defined, for all $s\ge 0$, by
\begin{align*}\Phi_j(s) = \left[\frac1{x_+}+\left(j-\frac12\right)c\right]^2-\frac{c^2}{4},
\end{align*}
and $c = \frac1k(\frac{1}{{x_-}}-\frac{1}{{x_+}})$, with $x_{\pm} = \sqrt{1+s\psi_{\pm}}$.
The minimax optimal excess risk, i.e., the optimal value of the objective in \eqref{mmxrisk}, is
\begin{align*}
   M(\psi_-,\psi_+,s,k)
   &= \frac{1}{k^2s}\left(\frac{1}{\sqrt{1+s\psi_-}}-\frac{1}{\sqrt{1+s\psi_+}}\right)^2.
\end{align*}

\end{theorem}
The proof is based on a combination of classical mathematical analysis and explicit calculations, including an ``infinite descent'' argument to find the minimizer.
Specifically, after some calculations, we can reduce the problem to evaluating the quantity given in equation \eqref{kp1} in the appendix. This is a minimum over $k+1$ terms.
We argue by contradiction that the solution is achieved when the $k+1$ quantities in \eqref{kp1} are equal. 
To show this, we suppose that at the minimum, the $k+1$ terms are not all equal, and choose the maximal term. Then, by careful analysis, we verify that by changing one of the optimization variables in that term, we can decrease the maximal term to equal its larger neighbor, which is a  contradiction to our assumption that we were at the minimum. Surprisingly, the mini-max optimal grid aligns with gradient descent with a decreasing learning rate. For details on this see the discussion within Section \ref{sec:app:compgrid}.

\subsection{Sub-optimality compared to the Minimax Grid}
\label{sec:SuboptCompGrid}
It is of interest to characterize the sub-optimality of ridge regression and gradient descent to the minimax optimal grid.  The following result is proved in Section \ref{pf:err}. Denote by $\{x\}$ the fractional part of the scalar $x$, and by $\lfloor x\rfloor$ its floor, the largest integer less than or equal to $x$. Denote by $C_{\text{Ridge}}(\delta) = \{\widehat{\beta}_{\lambda}\}_{\lambda \in \{\delta,2\delta,\dots,1\}}$ the class of estimators equal to ridge regression with regularization over a grid $\{\delta,2\delta,\dots,1\}$.

\begin{proposition}[Max-Min  Excess Risk of Ridge Regression]\label{err}
Under the conditions of Theorem \ref{mmxgrid}, let $\delta=1/k$, $\phi=k/\psi_--(\lfloor k/\psi_+ \rfloor +1 )$, and  $\psi_->k/(k+1)$. Then, the worst-case excess risk of the best estimator in the class $C_{\text{Ridge}}(\delta)$ of ridge regression estimators with  regularization over the grid $\{\delta, 2\delta,\dots,1\}$ $=$ $\{1/k,\ldots,1\}$ equals
$$
\sup_{\psi\in [\psi_-,\psi_+]}\min\limits_{\widehat{\beta} \in C_{\text{Ridge}}(\delta)} \mathcal{E}(\widehat{\beta})
= \frac{1}{s}\max\left\{Q_\star(j,\tau^{\star}_1)^2,Q_\star(j+1,\tau^{\star}_2)^2\right\}
$$
where $j = \lfloor k/\psi_+\rfloor$, $a= ks$, and
$$
Q_\star(j,\tau) =
\frac{a/[j+\tau]}{(1+a/[j+\tau])^{1/2}}
 \min\left(\frac{\tau}{j+a},\frac{1-\tau}{j+1+a}\right),
 $$
\begin{align*}
    \tau^\star_1 = 
          \max\left(\left\{ \frac{k}{\psi_+}\right \},\frac{\lfloor k/\psi_{+} \rfloor+a}{2(\lfloor k/\psi_{+} \rfloor+a)+1}\right)
          \!\! \quad\quad \!\!
           \text{and}
           \!\! \quad\quad \!\! 
          \tau^{\star}_2 
          = \min\left(\phi,\frac{\lfloor k/\psi_+\rfloor+1+a}{2(\lfloor k/\psi_+\rfloor+1+a)+1}\right).
\end{align*}
\end{proposition}

The condition $\psi_->k/(k+1)$ requires that $\psi_-$ is sufficiently large, while $\lfloor k/\psi_+ \rfloor +1 \le k/\psi_-$ requires that $\psi_+$ is sufficiently large compared to $\psi_-$ and $k$.

If $s = 1$, we can verify $\tau_1^*\in[1/3,1]$, and $a/[j+\tau_1^*]=k/[\lfloor k/\psi_+\rfloor+\tau_1^*] =O( \psi_+)$, where $A=O(B)$ means that $A\le CB$ for a constant $C$ not depending on $k$. Also
$$
\frac{\tau_1^*}{j+a}
=
O\left(\frac{1}{\lfloor k/\psi_+\rfloor+k}\right)
=
O\left(\frac{1}{k}\right).
$$
This shows that $Q_\star(j,\tau^{\star}_1) = O(1/k)$.
The same reasoning also shows that $Q_\star(j+1,\tau^{\star}_2) = O(1/k)$.
Thus, the upper bound on the excess risk of ridge is of order at most $O(1/k)$. As we will see numerically, the rate decays linearly in $1/k$ in certain ranges of $k$, and linearly in $1/k^2$ in other ranges. This is a subtle and perhaps unexpected behavior that is uncovered by our analysis.

A similar approach allows us to characterize the sub-optimality of gradient descent. We have the following result, proved in Section \ref{sec:ergdpf}. 
Let $C_{\text{GD}}(\eta,k) = \{\hat\beta_{\eta,i}\}_{i =1,\dots,k}$ denote the class of estimators associated to the first $k$ iterations of gradient descent initialized at $\hat\beta_{\eta,0}=0$ zero, with step size $\eta > 0$. We denote $a \vee b = \max\{a,b\}$ as well as $a \wedge b = \min\{a,b\}$. 

\begin{proposition}[Max-Min Excess Risk of Gradient Descent]\label{ergd}
Under the conditions of Theorem \ref{mmxgrid} and Proposition \ref{err}, consider  gradient descent with step size $0<\eta<1/s$. Let 
$$b=1- \eta s \in (0,1), \!\quad 
j = \lfloor \log_{1/b}(1+s\psi_-) \rfloor,\!\quad 
\tau^\star = \log_{1/b}(2/(b+1)),\!\quad 
\tau_- = \{\log_{1/b}(1+s\psi_-)\}$$ 
and suppose that 
$(1-\eta s)^{k}<(1+s\psi_+)^{-1}.$

Then, the worst-case excess risk of the best iterate of gradient descent out of the first $k$ iterates equals
$$
\sup_{\psi\in [\psi_-,\psi_+]}\min\limits_{\widehat{\beta} \in C_{\text{GD}}(\eta,t)} \mathcal{E}(\widehat{\beta})
=      b^{j  - \tau^\star} \max\big\{   Q_\star(\tau_- \vee \tau^{\star})^2, 
    b\cdot Q_\star( \tau_-\wedge \tau^{\star})^2
    \big\},
$$
where
\begin{align*}
    Q_\star(\tau) = \min\{ 1 - b^{\tau}, b (b^{\tau - 1} - 1) \}.
\end{align*}
\end{proposition}

The condition $(1-\eta s)^{k}<(1+s\psi_+)^{-1}$ requires $\eta$ and $k$ to be large enough. This assumption ensures that the gradient descent iterates reach far enough so that the optimal stopping point is before the last iteration. Relaxing this condition requires a more careful study of the behavior of the last iteration, which we do not expect to lead to major new insights. If $s=1$ and $\eta = C/k$ for some constant $C> 0$, then $(1-\eta s)^{k} \leq \exp(-C)$, and thus the aforementioned assumption is satisfied for $C$ sufficiently large with respect to $\psi_+$. From Proposition \ref{ergd} we then find that the excess risk is $O(\eta^2) =O(1/k^2)$, since the maximum of $Q_\star(\tau)$ is achieved at $\tau^{\star}$, and $b^{\tau^{\star}} = (1+b)/2 = 1 - \eta/2$.

\subsection{Experiments}
\begin{figure}%[H]
    \centering
    \includegraphics[width=0.32\textwidth]{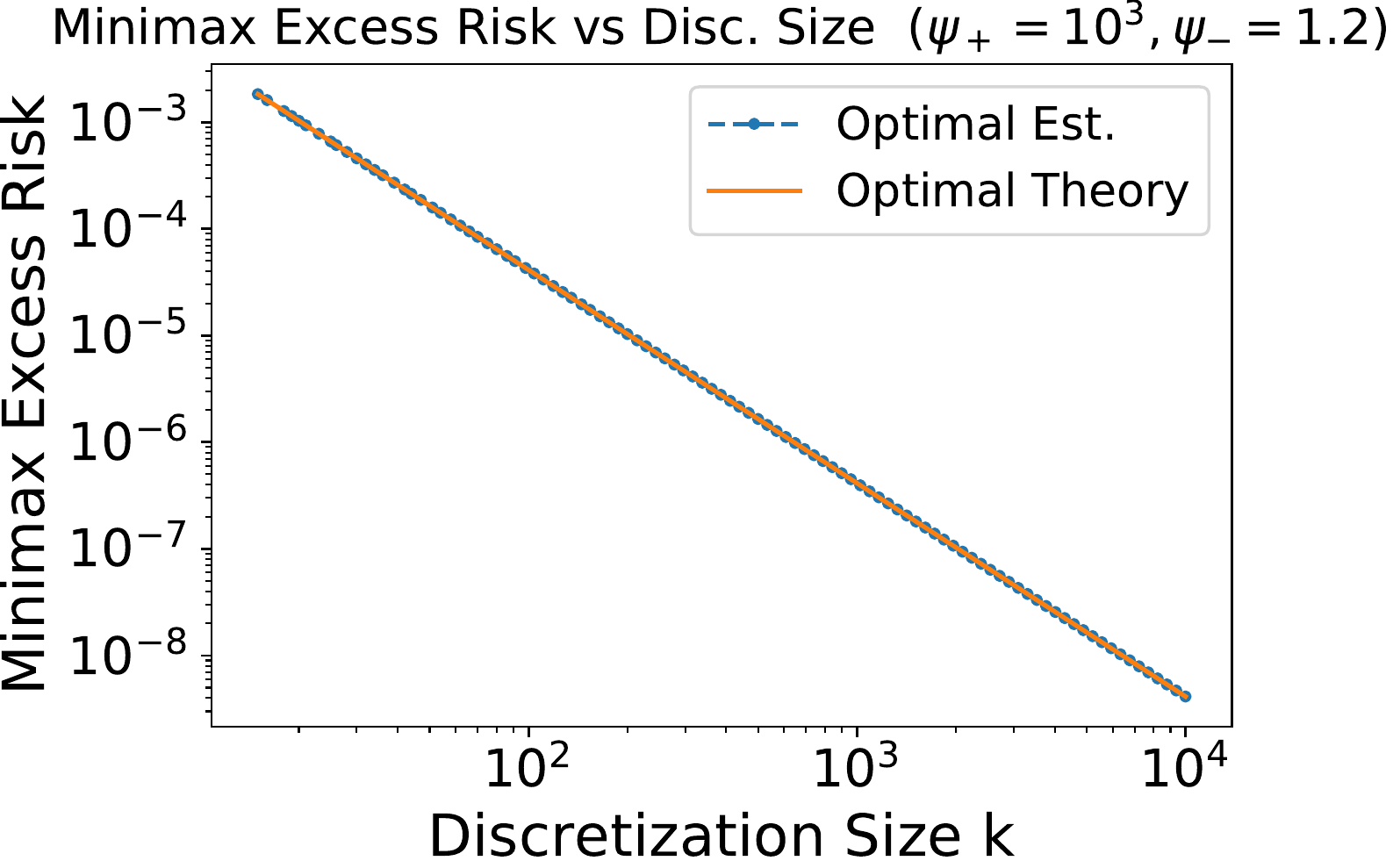}
    \includegraphics[width=0.32\textwidth]{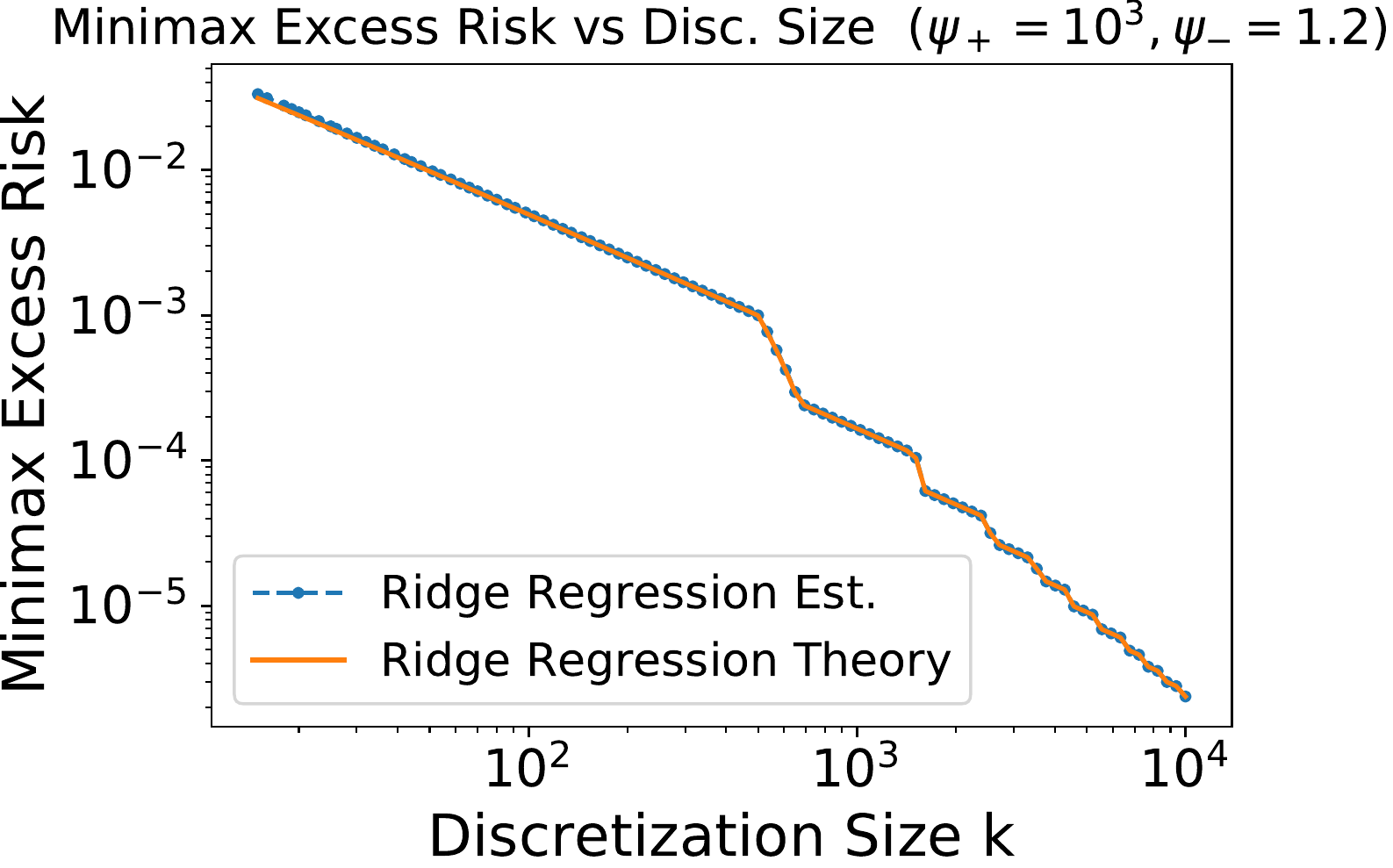}
    \includegraphics[width=0.32\textwidth]{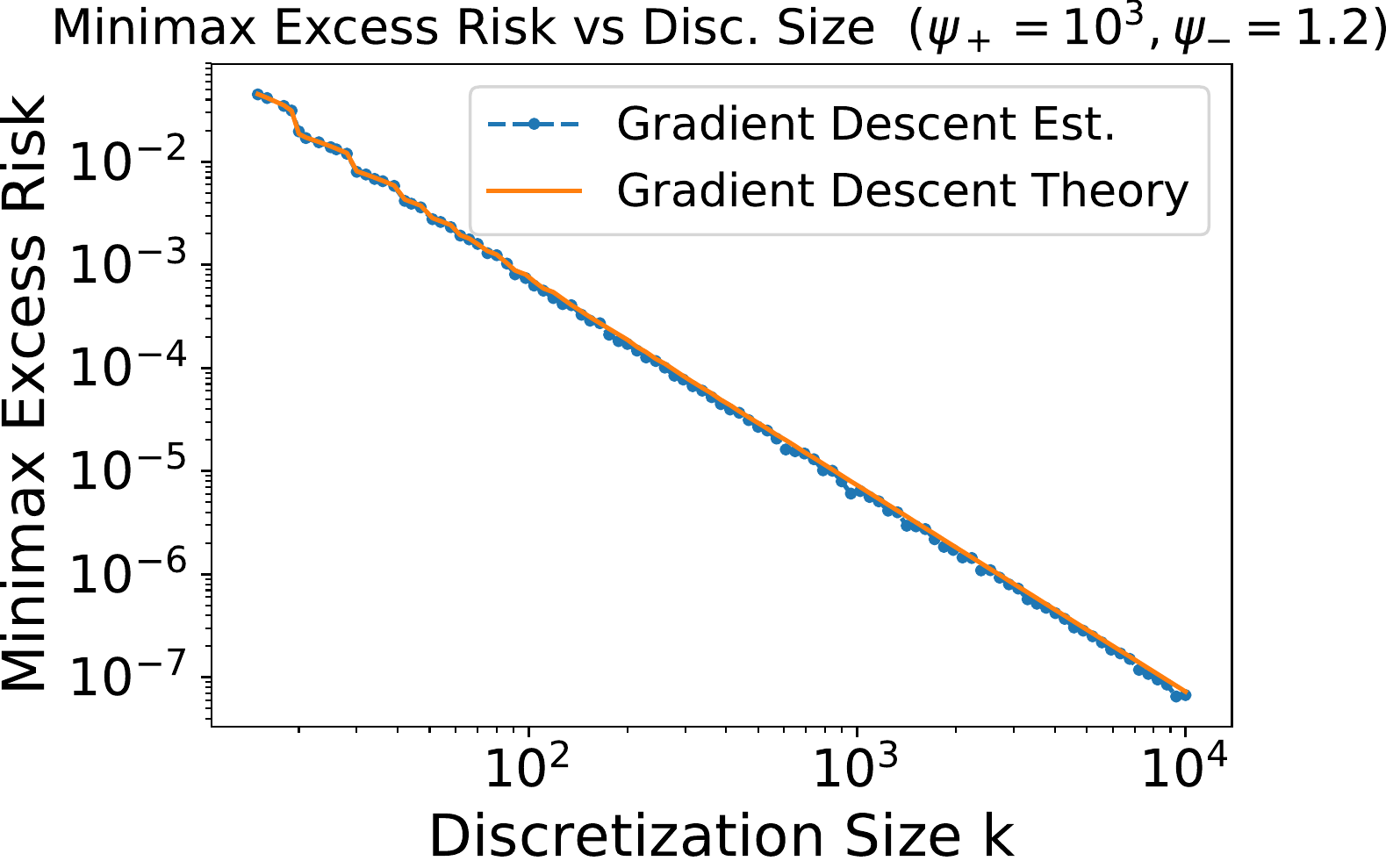}
    \caption{Plots of minimax excess risk versus discretization size $k \in [10^{1.2},10^{4}]$ for different methods: \textit{Left}: Optimal grids, \textit{Middle}: Ridge regression, \textit{Right}: Gradient descent. \textit{Orange Line}: Indicates estimate of minimax excess risk by discretizing $[\psi_{-},\psi_{+}]$ into $10^5$ points (to evaluate supremum in \eqref{mmxrisk}), \textit{Blue Line}: Theoretical quantity.   Gradient descent step size $\eta = 8/k$.}
    \label{fig:MinimaxVsTheory}
\end{figure}

\begin{figure}
\centering
\includegraphics[width=0.48\textwidth]{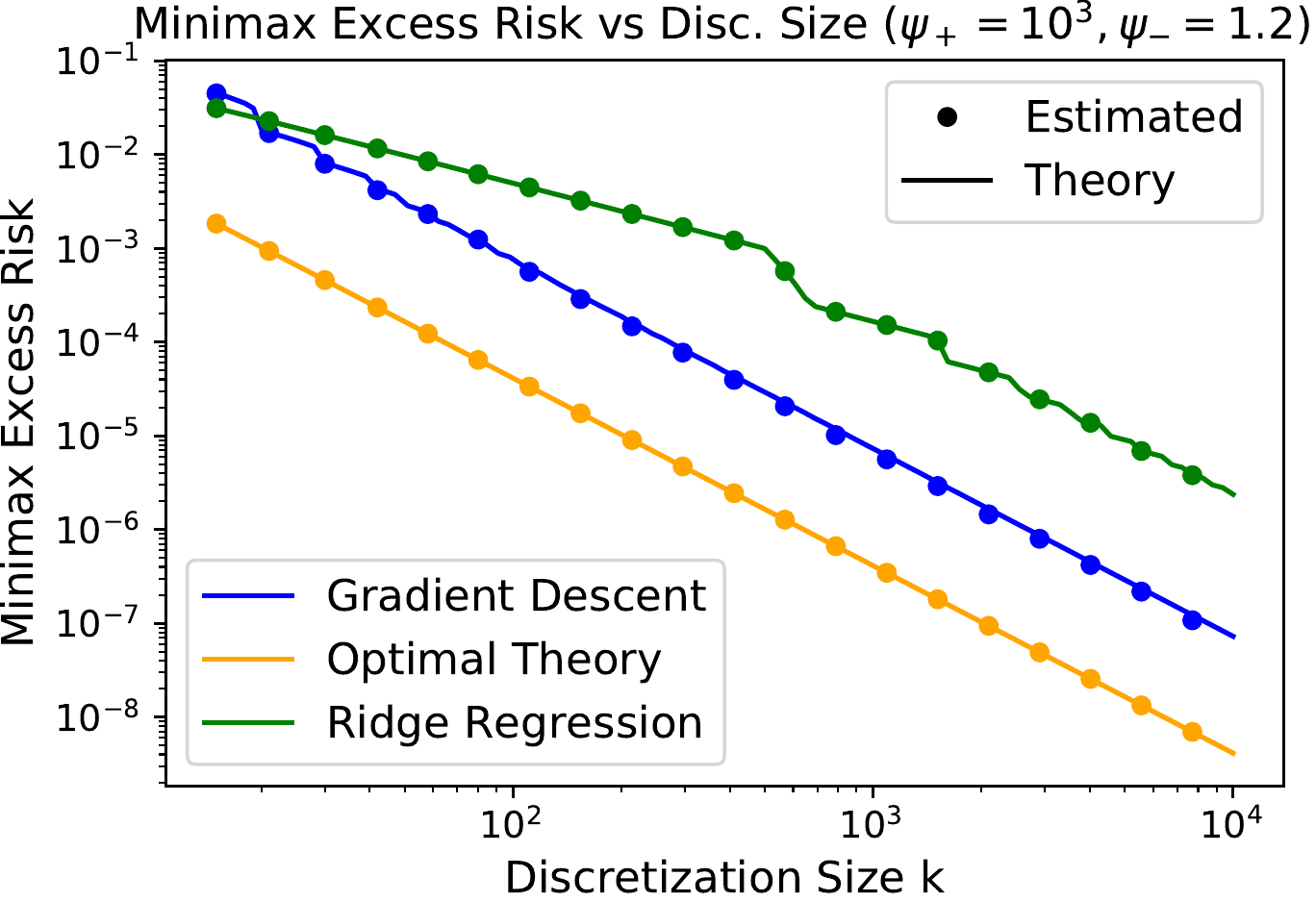}
\includegraphics[width=0.49\textwidth]{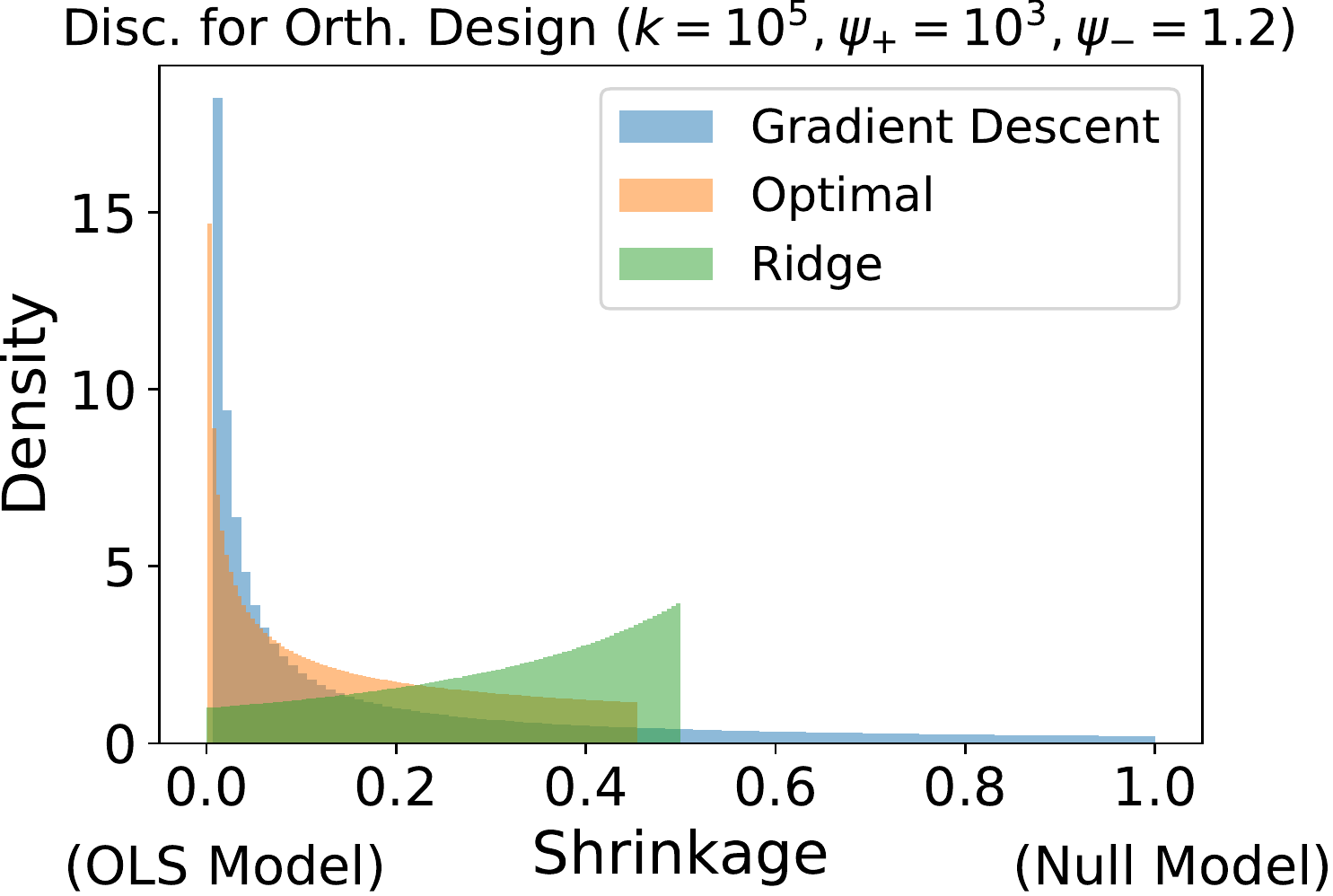}
\caption{\textit{Left}: Lines indicate the theoretical minimax risk for each class from Figure \ref{fig:MinimaxVsTheory} plotted side-by-side, with points showing numerical estimates. \textit{Right}: Histogram of spectral shrinkers for the orthogonal design case with various methods: \textit{gradient descent}, \textit{Optimal}: optimal design, and \textit{Ridge}: ridge regression with uniform discretization. \textit{OLS Model}: Ordinary Least Squares Model.}
\label{fig:orthogonalhist}
\end{figure}

Numerical results are provided to support the results for orthogonal designs.  
Figure \ref{fig:MinimaxVsTheory} shows an agreement between the theoretical results and numerically approximated minimax excess risks. The left plot of Figure \ref{fig:orthogonalhist} shows that gradient descent decays at the same $O(k^{-2})$ rate as the optimal grid in Theorem \ref{mmxgrid}. 
In contrast, ridge regression's excess risk decays at the rate $O(k^{-1})$ for small $k$, and then, seemingly, switches between rates of decay of $O(k^{-2})$ and $O(k^{-1})$. The right plot of Figure \ref{fig:orthogonalhist}, shows how different methods distribute their shrinkers (i.e., the values $\phi_j$) in the interval $(0,1)$. Gradient descent and the optimal grid place more mass towards the ordinary least squares (OLS) model than ridge regression tuned over a uniform grid.

\section{Conclusion}
We have introduced a 
framework to compare the statistical performance of \emph{classes of estimators}. 
In linear regression with a random ground truth,
we have compared gradient descent and ridge regression. 
We have found that, 
if the rank of the empirical covariance matrix is sufficiently large, the \emph{rate} at which the eigenvalues decays influences whether gradient descent or ridge regression is better, with a \emph{slower decay} leading to more favorable behavior for gradient descent over. 
This is surprising, as for the infinitely large set of all regularization parameters, the situation is reversed.
Our work lays foundations to investigate how classes of estimators---e.g., those obtained by implicit or explicit regularization---should be constructed and adapted to a specific problem instance, instead of considering optimally tuned estimators.

\section*{Acknowledgements}

ED was supported in part by the NSF grants BIGDATA (1837992),  CAREER (2046874), and the NSF--Simons Foundation Award on the Mathematical and Scientific Foundations of Deep Learning (2031895). PR was supported in part by the Alan Turing Institute under the EPSRC grant EP/N510129/1.

\appendix 

\section{Additional Results and Discussion}
\label{sec:app:additional_results_discussion}
In this section we present additional results and discussion to support the main body of the manuscript. This section is structured as follows. 
\begin{itemize}
    \item Section \ref{sec:app:lit} discusses additional related literature related to shrinkage estimators for least squared regression (Section \ref{sec:lit})
    \item Section \ref{sec:app:log_grid} compares the performance of ridge regression tuned with a logarithmic grid to gradient descent when the eigenvalues decay slowly (Theorem \ref{thm:log-ridgeVsGD}). 
    \item Section \ref{sec:app:add_discussion} includes discussion of the framework investigated within this work as well as future research direction (Section \ref{gd-vs-ridge}).  
    \item Section \ref{sec:ProofSketch} presents a proof sketch of Theorem \ref{thm:Informal} which compares ridge regression to gradient descent with eigenvalues that decay at different rates (Section \ref{gd-vs-ridge}). 
    \item Section \ref{sec:app:add_bounds_individual} considers the performance of
    individual classes of ridge regression and gradient descent models in the coarse regime, as well as ridge regression tuned with a logarithmic grid (Section \ref{b-ind}). 
    \item Section \ref{sec:app:compgrid} investigates the structure of spectral shrinkers in the orthogonal design case, and shows that gradient descent with a polynomial decreasing step size equals the minimax grid (Section \ref{mmx-ortho-sec}). 
\end{itemize}

\subsection{Related Literature}
\label{sec:app:lit}

Our work is related to the study of the asymptotic performance of ridge regression in a high-dimensional setting where both the ambient dimension and number of data points go to infinity in proportion to one another, e.g., \cite{tulino2004random,karoui2011geometric,karoui2013asymptotic,dicker2013optimal,dicker2016ridge,dobriban2018high,hastie2019surprises,ali2019continuous,ali2020implicit,richards2021asymptotics,wu2020optimal,liu2019ridge} etc. Assuming the data follows a linear model or certain nonlinear models, the works  \cite{dicker2016ridge,dobriban2018high,hastie2019surprises,ali2019continuous,ali2020implicit,richards2021asymptotics,wu2020optimal,liu2019ridge} study the estimation and prediction error of ridge regression in terms of integrals against the eigenvalue distribution of the empirical covariance matrix. By making assumptions on the covariate structure, they leverage the convergence of the empirical eigenvalue distribution to variants of the Marchenko–Pastur law \cite{marchenko1967distribution,bai2010spectral} to obtain results in terms of population quantities. While our results are in terms of the eigenvalue spectrum, we alternatively focus on the fixed design setting with spectral distributions that are distinct from those classically studied in random matrix theory, such as eigenvalues with a power law or exponential decay. 
Such distributions have been observed empirically in many different applications, including in finance, biology, and signal processing. For examples, see Section 8 of the arXiv version of the paper \cite{dobriban2017sharp}.

\sloppy Our focus differs from the literature on  tuning regularization for least squares via cross-validation \cite{arlot2010survey,kale2011cross,miolane2018distribution,xu2019consistent,hastie2019surprises,patil2021uniform}, which considers empirical model selection criteria and their consistency. Similarly, the works \cite{birge2001gaussian,baraud2009gaussian,arlot2009data,JMLR:v10:arlot09a,bellec2020cost} investigate selection criteria aiming to recover the performance of the best-in-class model (often referred to as oracle within). In contrast, we focus on choosing the class of models which arises from, for instance, a 
\emph{finite grid of regularization parameters}.

For an orthogonal design matrix, the least squares problem reduces to mean estimation. However, unlike the usual setting where the mean parameter is fixed, here we consider a random-effects model where the mean parameter is isotropically distributed, and are interested in the average-case behavior over this prior. In this sense, our work has a Bayesian interpretation \cite{lehmann1998theory,berger2013statistical}.

There is a great deal of other work studying statistical properties of optimization methods. \cite{skouras1994estimation} show that in linear models, early stopped gradient flow on the ridge regression objective with a sufficiently small positive regularization parameter can outperform ridge regression in estimation error for \emph{a fixed design and regression parameters}. In contrast, for the random-effects models we consider, ridge regression is average-case (i.e., Bayes) optimal under certain assumptions.  Other connections between algorithmic regularization and penalization have also been studied \cite{qian2019connections}.

The related work \cite{Dhillon2013Risk} shows that the estimation risk of ridge regression and principal component regression  are within a multiplicative factor of four. Our approach is distinctly different, beyond just considering a different pair of estimators (early stopped gradient descent in place of principal component regression). First, \cite{Dhillon2013Risk} compares estimators with the \emph{same} hyperparameter choice (denoted $\lambda > 0 $), while we compare each estimator with their \emph{best} hyperparameter choice from two given sets. This is motivated from the best hyperparameter choice for principal component regression (or any general estimator) potentially not aligning with that for ridge regression. Second, our analysis also accounts for finding the best regularization parameter (by searching over a grid, say). This depends upon the estimator's sensitivity with respect to its hyperparameter, specifically, the change due to the grid resolution around the optimal hyperparameter choice, and can result in estimators outperforming ridge regression (Theorem \ref{thm:Informal}). As principal component regression can be represented via a spectral shrinker, one future research direction is to extend the results of \cite{Dhillon2013Risk} into the framework developed within our work.

\subsection{Ridge Regression with Logarithmic Grid versus Gradient Descent} 
\label{sec:app:log_grid}

If $\Gamma = \{\lambda_{j}\}_{j=1}^{k}$ is a grid of strictly positive values, then we denote $\log(\Gamma) = \{\log(\lambda_{j})\}_{j=1}^{k}$. 
Additionally,  define for $k \geq 2$ and $\lambda_{\min} > 0$ the collection of ridge regression estimators tuned over a logarithmic grid as 
$$\mathcal{C}^{\text{log-Ridge}}(\lambda_{\min},k) := \{\widehat{\beta}_{\lambda_j} : \lambda_{j} = \exp\left(-(j-1)\log(1/\lambda_{\min})/(k-1) \right), \, j=1,\dots,k\}.$$
With this we have the following result.
\begin{theorem}
\label{thm:log-ridgeVsGD}
Suppose $0 < \sigma_{\min} < \sigma/\sqrt{2}$ with $\lambda_{\min} = d \sigma^2_{\min} /n$. For $k > 1$ let $\eta = 1/(k \lambda_{\min})$ and assume $k \geq 3/\lambda_{\min}$ and $(\eta \lambda_\star)^{-1} \not\in \{0,1,2,3,\dots,k-1\}$. Let $\delta = \log(1/\lambda_{\min})/(k-1)$ with $\lambda_j = \exp\big( -(j-1) \delta\big)$ for $j \geq 1$ and $\Gamma = \{\lambda_j \}_{j=1}^{k}$. Suppose $s_i = i^{-\alpha}$ for $i=1,\dots,r$,  $\alpha \in (0,1)$. Then there exists a threshold $r_{\alpha,\lambda_{\star},k,\lambda_{\min}} > 1$, such that if $r \geq r_{\alpha,\lambda_{\star},k,\lambda_{\min}}$, then 
\begin{align*}
    \mathcal{S}\left( \mathcal{C}^{\text{GD}}(\eta,k), 
    \mathcal{C}^{\text{$\log$-Ridge}}(\lambda_{\min},k)
    \right)  
    & \simeq 
    \frac{1}{\mathrm{Dist}_{\delta}(\log(\lambda_{\star}),\log(\Gamma))^2} \frac{ 1}{\log^2(n/(d \sigma_{\min}^2))} \frac{\sigma^4}{\sigma_{\min}^4}.
\end{align*}
\end{theorem}
The above theorem shows that the relative sub-optimality scales in this case as $\Theta( \log^{-2}(1/\lambda_{\star}))$
when other parameters are fixed,
when $\sigma_{\min} = \Theta(\sigma)$ and the eigenvalues decay according to a slow polynomial power $s_i= i^{-\alpha}$ for $1 \geq \alpha > 0$. 
This reflects Theorem \ref{thm:Informal} above, in that the relative sub-optimality goes to zero when the optimal regularization $\lambda_{\star} = d \sigma^2 /n$ goes to zero, i.e., $d/n \rightarrow 0$ or $\sigma \rightarrow 0$. We note the change from a polynomial decay in $\lambda_{\star}$ in Theorem \ref{thm:Informal} to a logarithmic decay in Theorem \ref{thm:log-ridgeVsGD}. 
Thus we expect that the performance of ridge regression over a logarithmic grid is closer to that of gradient descent, as $\lambda_{\star} \rightarrow 0$ in this case. 

\subsection{Discussion of our framework: Gradient Descent versus Ridge Regression }
\label{sec:app:add_discussion}
Our framework aims to address the issue that  hyperparameter tuning can strongly affect the behavior of algorithms. Already in a very simple setting, namely in random effects linear models, we show that this leads to a sophisticated behavior of the relative performance of two truly fundamental classes of regularization methods: ridge regression and gradient descent.
This discovery is only possible due to our framework of comparing best in class algorithms; which in our view demonstrates the power of our perspective.

There are a number of important further considerations to discuss.
\begin{enumerate}
    \item In future work, certain key algorithmic considerations need to be taken into account for a more thorough comparison between methods. For instance, as we have shown, our framework can be extended to incorporate that some algorithms are computationally more expensive than others, by simply changing the class of algorithms to match the computational cost. In addition, it will be important to apply it to study more sophisticated classes of algorithms, and for more general data generating models. 
    \item At the moment, our framework compares methods according to their true expected performance on the test distribution. In practice, the methods would be compared
    after having selected the best-in-class on a finite training dataset, and
    according to their performance on a finite a test set. 
    This will introduce some amount of noise, and thus has to be taken into account in the analysis. However we expect that with sufficiently large training and test sets, our conclusions would largely be unchanged. 
    The analysis of the noise due to the finiteness of the training and test set remains an important future problem.
\end{enumerate}

\subsection{Proof Sketch for Theorem \ref{thm:Informal}}
\label{sec:ProofSketch}
We begin by following the proof of Lemma \ref{lem:ErrorDecomp} which, if we denote $M_{\Phi}(s) = 1 - \Phi(s) s$, expresses the excess estimation error of $\widehat{\beta}_{\Phi}$  as
\begin{align*}
    \mathcal{E}(\widehat{\beta}_{\Phi}) 
    = 
   \frac{ \lambda_{\star} }{d} \sum_{i=1}^{r} \left( \frac{1}{\lambda_{\star}} + \frac{1}{s_i} \right)\left( M_{\Phi}(s_i) - \frac{\lambda_{\star}}{\lambda_{\star} + s_i} \right)^2.
\end{align*}
We then consider the appropriate
ridge regression classes  $\mathcal{C}^{\text{Ridge}}(\lambda_{\min},k)$ =$\{\widehat\beta_{\Phi^{\text{Ridge}}_i}\}_{i=1,\dots,k}$ where $M_{\Phi_i^{\text{Ridge}}}(s) = (\lambda_{\min} + (i-1) \delta)/(s + (\lambda_{\min} + (i-1) \delta))$; and gradient descent classes $\mathcal{C}^{\text{GD}}(\eta,k) = \{\widehat\beta_{\Phi^{\text{GD}}_i}\}_{i=1,\dots,k}$ where $M_{\Phi_i^{\text{GD}}}(s)=(1-\eta s)^{i}$. We begin by bounding the best-in-class error for each method, i.e., the performance of the best model in each class $\min_{\hbeta \in \mathcal{C}} \mathcal{E}(\mathcal{C})$, after which we take the ratio. Bounding the performance of the best-in-class models is described in the following two paragraphs, with their ratio considered in a third paragraph.

\textbf{Best-in-Class for Ridge Regression.} If we denote $\lambda_i = \lambda_{\min} + (i-1)\delta$ for $i=1,\dots,k$, and recall that $\lambda_{\star} = \lambda_{\min} + (j+\ep)\delta$, the best-in-class error for ridge regression reduces to evaluating the following minimum
\begin{align*}
    \mathcal{E}(\mathcal{C}^{\text{Ridge}}(\lambda_{\min},k)) 
    & = 
    \lambda_{\star} \min_{t=1,\dots,k} \frac{1}{d} \sum_{i=1}^{r} \left( \frac{1}{\lambda_{\star}} + \frac{1}{s_i} \right)\left( \frac{\lambda_{t}}{\lambda_{t} + s} - \frac{\lambda_{\star}}{\lambda_{\star} + s_i} \right)^2\\
    & = 
    \lambda_{\star} \min_{\kappa =\{-\ep, 1-\ep\} }
    \frac{1}{d} \sum_{i=1}^{r} \left( \frac{1}{\lambda_{\star}} + \frac{1}{s_i} \right) \frac{s_i^2 \kappa^2 \delta^2}{(\lambda_{\star} + \kappa \delta + s_i)^2(\lambda_{\star} + s_i)^2},
\end{align*}
where the second equality arises from monotonicity, allowing  the minimum to be reduced to the two regularization choices in $\{\lambda_i\}_{i=1,\dots,k}$ closest to $\lambda_{\star}$. When $\ep = 1/2$ and the discretization length is smaller than the optimal regularization $\delta \leq \lambda_{\star}$, each term in the series is then on the order $\delta^2 \frac{s}{\lambda_{\star}} \frac{1}{(s+\lambda_{\star})^3}$, with the factor in-front of $\delta^2$ describing the sensitivity to the discretization length. When $s > \lambda_{\star}$ this factor is roughly $\frac{1}{\lambda_{\star} s^2}$, and for $s < \lambda_{\star}$ is roughly $\frac{s_i}{\lambda_\star^4}$, allowing us to then write
\begin{align}
\label{equ:Ridge:Sketch}
    \mathcal{E}(\mathcal{C}^{\text{Ridge}}(\lambda_{\min},k)) 
    &
    \simeq 
    \lambda_{\star} \delta^2 \frac{1}{d}\left( \sum_{ s_i > \lambda_{\star}} \frac{1}{\lambda_{\star} s^2} + 
    \sum_{ s_i \leq \lambda_{\star}} \frac{s_i}{\lambda_{\star}^4}
    \right).
\end{align}

\textbf{Best-in-Class for Gradient Descent.} The best-in-class error in this case reduces to evaluating the following minimum 
\begin{align*}
    \mathcal{E}(\mathcal{C}^{\text{GD}}(\eta,k)) 
    = 
    \lambda_{\star} \min_{t=1,\dots,k} \frac{1}{d} \sum_{i=1}^{r} \left( \frac{1}{\lambda_{\star}} + \frac{1}{s_i} \right)\left( (1-\eta s_i)^{t} - \frac{\lambda_{\star}}{\lambda_{\star} + s_i} \right)^2.
\end{align*}
The first step is to switch from an iteration viewpoint, i.e., $t=1,\dots,k$, to an eigenvalue perspective. 
We consider the optimal per-eigenvalue-iteration $t^{\star}$ as the possibly non-integer-valued iteration number which matches the optimally tuned ridge regression estimator for eigenvalue $s$, i.e., for which $(1-\eta s)^{t^{\star}(s)} = \lambda_{\star} /(\lambda_{\star} + s)$.
If we choose $t = \lceil (\eta \lambda_{\star})^{-1} \rceil \geq (\eta \lambda_{\star})^{-1} = \max_{s \in [0,1/\eta] }t^{\star}(s)$ we can then bound 
\begin{align*}
    \mathcal{E}(\mathcal{C}^{\text{GD}}(\eta,k)) 
    \lesssim 
    \lambda_{\star} \frac{1}{d}
    \left( \sum_{i : s_i > \lambda_{\star}} \frac{\lambda_{\star}}{s_i^2} 
    + 
    \sum_{i : s_i \leq \lambda_{\star}} \frac{1}{s_i} 
    \left(1 - (1-\eta s)^{\lceil (\eta \lambda_{\star})^{-1} \rceil - t^{\star}(s_i)}\right)^2
    \right),
\end{align*}
where for $s> \lambda_{\star}$ we have $0 \leq \frac{\lambda_{\star}}{\lambda_{\star} + s} - (1-\eta s)^{t} \leq \frac{\lambda_{\star}}{s}$ and for $s \leq \lambda_{\star}$ we bound $\frac{\lambda_{\star}}{\lambda_{\star} + s} - (1-\eta s)^{t} = \frac{\lambda_{\star}}{\lambda_{\star} + s_i}( 1- (1- \eta s)^{ t - t^{\star}(s)}) \leq ( 1- (1- \eta s)^{ t - t^{\star}(s)})$. The second series depends upon the smoothness of $t^{\star}(s)$ through the difference $\lceil (\eta \lambda_{\star})^{-1} \rceil - t^{\star}(s)$. A technical result (Proposition \ref{prop:IterationsLimit}) then shows, when $\eta \lambda_{\star} \leq 1$, that for $s \in [0,\lambda_{\star}]$ 
\begin{align*}
    \frac{1}{\eta \lambda_{\star}} - t^{\star}(s) \leq \frac{3}{2} \frac{s}{\eta \lambda_\star^2}.
\end{align*}
Using this, with $\lceil x \rceil \leq x+1$ alongside Bernouli's inequality $(1+x)^{r} \geq 1 + xr$ for $r \geq 1$ and $x \geq -1$, then yields the upper bound $1 - (1-\eta s)^{\lceil (\eta \lambda_{\star})^{-1} \rceil - t^{\star}(s_i)} \lesssim \eta s(1 + \frac{s}{\eta \lambda_{\star}^2})$. Plugging this in we arrive at 
\begin{align}
\label{equ:GDBound:Sketch}
    \mathcal{E}(\mathcal{C}^{\text{GD}}(\eta,k)) 
    \lesssim
    \frac{\lambda_{\star} }{d}
    \left( \sum_{i : s_i > \lambda_{\star}} \frac{\lambda_{\star}}{s_i^2} 
    + 
    \sum_{i : s_i \leq \lambda_{\star}} \max\big\{ \eta^2 s_i,\frac{s_i^3}{\lambda_{\star}^4} \big\}
    \right).
\end{align}

\textbf{Ratio of Best-in-Class Performances.} Taking the ratio of \eqref{equ:GDBound:Sketch} and \eqref{equ:Ridge:Sketch} we find   
\begin{align}
\label{equ:ratio:sketch}
    \mathcal{S}(\mathcal{C}^{\text{GD}}(\eta,k),\mathcal{C}^{\text{Ridge}}(\lambda_{\min},k)) 
    \lesssim \frac{1}{\delta^2} \frac{ \sum_{i : s_i > \lambda_{\star}} \frac{\lambda_{\star}}{s_i^2} }{ 
    \sum_{s_i \leq \lambda_{\star}} \frac{s_i}{\lambda_{\star}^4}
    }
    + 
    \frac{1}{\delta^2} \frac{\sum_{i : s_i \leq \lambda_{\star}} 
    \max\big\{ \eta^2 s_i,\frac{s_i^3}{\lambda_{\star}^4}\big\} 
    }{
    \sum_{i : s_i \leq \lambda_{\star}} \frac{s_i}{\lambda_{\star}^4}
    }.
\end{align}
When $\sum_{s_i \leq \lambda_{\star}} s_i$ is large, e.g., $s_i = i^{-\alpha}$ for $\alpha \in (0,1)$, and $r$ large, the second term will then dominate. Since the maximum in the numerator is $\eta^2 s_i$ for $s_i \leq \eta \lambda_{\star}^2$, the second term is on the order of $\frac{\eta^2}{\delta^2} \lambda_{\star}^4 = \frac{\lambda_{\star}^4}{\lambda_{\min}^2} \frac{1}{(1-\lambda_{\min})^2}$, where we recall that $\eta = (k \lambda_{\min})^{-1}$ and $\delta = (1-\lambda_{\min})/(k-1)$. More generally, the bound for this expression arises from the error of ridge regression with respect to eigenvalues below the regularization level---$s_i \leq \lambda_{\star}$---scaling as $\frac{\delta^2}{\lambda_{\star}^4}$, whereas the error of gradient descent for $s_i \leq \eta \lambda_{\star}^2$ scales as $\eta^2$.

\subsection{Bounds for Individual Classes}
\label{sec:app:add_bounds_individual}
In this section present additional bounds for classes of estimators associated to gradient descent and ridge regression, complimenting results contained within Section \ref{b-ind}. 
Section \ref{sec:app:add_bounds_individual:RR} has results associated to ridge regression. Section \ref{sec:app:add_bounds_individual:GD} presents results associated to gradient descent. 

\subsubsection{Ridge Regression}
\label{sec:app:add_bounds_individual:RR}
We provide the bound for the case of a \textit{coarse discretization} $\delta =(1-\lambda_{\min})/(k-1) > \lambda_{\star}$. 
\begin{proposition}[Ridge Regression Estimator, Coarse Discretization]
\label{lem:RidgeBounds:Coarse}
Let $\lambda_{\min} = 0$, $\delta = (1-\lambda_{\min})/(k-1)$ and suppose $\lambda_{\star} = \delta - \delta'$ for some $\delta/2 \geq \delta' > 0$. Then 
\begin{align*}
    \frac{1}{4} \left(\frac{\delta'}{\delta}\right)^2  \lambda_{\star}
    \int G^{\text{Coarse Ridge}}(s) d \widehat{H}(s) \leq \mathcal{E}(\mathcal{C}^{\text{Ridge}}(0,k))
    \leq  \left(\frac{\delta'}{\delta}\right)^2  \lambda_{\star} \int G^{\text{Coarse Ridge}}(s) d \widehat{H}(s)
\end{align*}
where 
\begin{align}\label{cri}
    G^{\text{Coarse Ridge}}(s)
    = 
    \frac{s}{\lambda_{\star}(s+\lambda_{\star})}
    \min\left\{\frac{\delta^2}{s^2},1\right\}.
\end{align}
\end{proposition}
We now provide the bound for a logarithmic grid. 
\begin{proposition}[Ridge Regression Estimator, Log-Grid Discretization]
\label{lem:log_grid}
Assume for some $k > j_{\star} > 2$ and $\epsilon \in (0,1)$ that $\lambda_{\star} = \exp\big( (j_{\star}-2 + \epsilon)\log(\lambda_{\min}) /(k-1)  \big)$
as well as that $k \geq 1 + 3 \log(1/\lambda_{\min})$. Then 
\begin{align*}
     \frac{\min\{\epsilon, \! 1\!- \! \epsilon\}}{12}  \lambda_{\star} \!\! \int \! \! G^{\text{$\log$-Ridge}}\!(s) d \widehat{H}(  s  ) \! \leq \! \mathcal{E}(\mathcal{C}^{\text{$\log$-Ridge}}\! (\lambda_{\min},\! k))
    \! \leq \!  \min\{\epsilon,\! 1 \! - \!\epsilon\} \lambda_{\star} \!\! \int \!\! G^{\text{$\log$-Ridge}}\!(s) d \widehat{H}(s)
\end{align*}
where 
\begin{align*}
    G^{\text{$\log$-Ridge}}(s) 
    = 
    \frac{\log^2(1/\lambda_{\min})}{(k-1)^2} 
    \frac{s\lambda_{\star}}{(\lambda_{\star} + s)^3}. 
\end{align*}
\end{proposition}

\subsubsection{Gradient Descent}
\label{sec:app:add_bounds_individual:GD}
We now provide an upper bound on the error of gradient descent in a coarse discretization regime when $t \lesssim 1/(\eta \lambda_{\star})$. 
\begin{proposition}[Gradient Descent Estimators, Coarse Discretization]
\label{lem:GDBounds:Coarse}
Let $\lambda_{\min} = 0$, $0 < \eta < 1/\max\{s_1,\lambda_{\star}\}$ and $t < \frac{1-\eta \lambda_{\star}}{2 \eta \lambda_{\star}}$. Then 
\begin{align*}
    \mathcal{E}(\mathcal{C}^{\text{GD}}(\eta,t)) 
    \leq 
    2 \lambda_\star \int G^{\text{Coarse GD}}(s) d \widehat{H}(s) 
\end{align*}
where 
\begin{align}\label{coarsegd}
    G^{\text{Coarse GD}}(s) 
    = 
    \begin{cases}
    \frac{\lambda_{\star}}{s^2} & \text{ for } s \in \big(\frac{1}{\lambda_{\star}(\eta t)^2},\frac{1}{\eta}\big] \\
    \frac{ 1}{\lambda_{\star} (\eta t s)^4} & \text{ for } 
    s \in \big( \frac{1}{\eta t},\frac{1}{\lambda_{\star}(\eta t)^2}\big]\\
    \frac{1}{\lambda_{\star}} & \text{ for } s \in (\lambda_{\star},\frac{1}{\eta t}] \\
    \frac{s}{\lambda_{\star}^2} 
    & \text{ for }  s \in (0,\lambda_{\star}].
    \end{cases}
\end{align}
\end{proposition}

\subsection{Comparing Grids for Orthogonal Designs}
\label{sec:app:compgrid}
We can compare other grids with the minimax optimal grid via their spectral filters, defined as $\Phi_j$, for $j=1,\dots,k$. 
Since the design is isotropic, any spectral shrinker is a constant over the spectrum.
Thus, by fixing $s$ from \eqref{Th}, we can view the filters as constants in $s$, and only a function of $j$, denoting them as $\phi_j$:
$$
\phi_j = \Phi_j(s).
$$
We can express the spectral filters of several methods as follows.
    \benum
    \item {\bf Ridge Regression}: For a general design, the spectral filter is $\Phi_j(s) = \lambda_j/(\lambda_j+s)$ for some $\lambda_j>0$. As a function of $j$, for a uniform discretization we have $\lambda_j = j/k$. Thus, the filter values are $\phi_j = a'j/(a'j+b')$, $j=1,\ldots, k$, for certain $a',b'>0$.
    \item {\bf Gradient Descent}: For a general design, the spectral filter is $\Phi_j(s) = (1-\eta s)^j$ as a function of the iteration $j$, for some $0<\eta<1/s$. 
    It is convenient to sort these values in increasing order in $j$, consistent with the other methods, writing $\phi_j  =b_0/b^j$ for some  $b_0>0$ and $b=1/(1-\eta s)>1$.
    \item {\bf Minimax Optimal}: We have shown in Theorem \ref{mmxgrid} that the minimax optimal grid is $\phi_j = (A+B j)^2-C^2$, for certain $A,B,C>0$.
    \eenum
    Thus as a function of $j$, the grids have the form:
         \begin{align*}
         a'j/(a'j+b'),\qquad b_0/b^j,\qquad(A+B j)^2-C^2.
         \end{align*}
    They are normalized to take values in the range $[0,1]$.
    Most of the values of the GD grid are concentrated around its minimum. This and the minimax optimal grid ares \emph{convex}\footnote{We recall that a sequence $(a_k)_{k\ge 1}$ is convex if $a_{k+1}+a_{k-1}\ge 2 a_k$ for all $k\ge 2$, see e.g., \cite{wayne1973convex}.} in $j$. The ridge grid is \emph{concave} in $j$. 
    
 {\bf GD with Decreasing Learning Rate can be Minimax Optimal.} It turns out that taking an appropriately decreasing learning rate in GD leads to the minimax optimal grid. One can verify that the minimax optimal grid values $\phi_j$ can be written as $\phi_j = a_{j-1} a_j$, where $a_j={x_+}^{-1}+jc$, $j=0,\ldots, k$ is an arithmetic progression going from $x_{+}^{-1}$ to $x_{-}^{-1}$. Since we sort the grids in increasing order, GD with learning rate $\eta_j\in(0,1/s)$ at the $j$-th iteration obeys the recursion $\phi_{j+1}  = \phi_j/(1-\eta_js)$. To match this with the minimax grid, we need
    \begin{align*}
         1-\eta_js&=\phi_j/\phi_{j+1}
         =[a_{j-1} a_j]/[a_{j} a_{j+1}]
         =a_{j-1}/a_{j+1}\\
         &=[{x_+}^{-1}+(j-1)c]/[{x_+}^{-1}+(j+1)c]
         =1-2c/[{x_+}^{-1}+(j+1)c].
    \end{align*}
    This shows that the matching learning rate schedule for GD has the form
    \begin{align*}
         \eta_j
         &=\frac{2c}{s\cdot [{x_+}^{-1}+(j+1)c]}.
    \end{align*}
    Such learning rate schedules are well known in optimization \cite{nesterov2013introductory}. For instance, for such a learning rate one can get 
	an $O(1/j)$ convergence rate, where $j$ is the number of iterations, for stochastic gradient descent on strongly convex functions, see e.g.,   Theorems 5, 6 of \cite{vaswani2018fast} for results in this direction.

In contrast, representing the minimax optimal grid as ridge regression estimators over a certain regularization grid does not appear to lead to a cleanly expressed grid.

\section{Proofs for Gradient Descent and Ridge Regression with General Designs}
In this section we present the proofs associated to the estimation error of gradient descent and ridge regression for general design matrices. 
\begin{itemize}
    \item Section  \ref{sec:Proof:FirstLemma} presents the proof of Lemma \ref{lem:ErrorDecomp} which decomposes the estimation error of any estimator $\widehat{\beta}_{\Phi}$ into the estimation error of optimally tuned ridge regression plus a non-zero term. 
    \item Section \ref{sec:ProofInformal}  summarizes the relative sub-optimality between ridge regression and gradient descent with eigenvalues that decay at different rates (Theorem \ref{thm:Informal}). 
    \item Section \ref{sec:proof:LowDimensional}  considers upper and lower bounds on the relative sub-optimality between ridge regression and gradient descent with eigenvalues that decay at a slow power law, in the low dimensional case (Corollary \ref{cor:LowDimensional}). 
    \item Section \ref{sec:Proof:Subopt:Slow} compares the performance of ridge regression and gradient descent when the eigenvalues decay at a slow power law rate, in the general case (Theorem  \ref{thm:mainGDvsRidge_Upper} and  \ref{thm:mainGDvsRidge_Lower}). 
    \item  Section \ref{sec:sub-opt-limit} considers formal results associated to the limits presents within the contribution section (Corollary \ref{RiskLimit:Formal} and  \ref{RiskLimit:Lower:Formal}). 
    \item Section  \ref{sec:proof:RidgeBounds} has estimation error bounds for classes of estimators constructed from ridge regression with a discretized regularization parameter (Proposition  \ref{lem:RidgeBounds} and \ref{lem:RidgeBounds:Coarse}).
    \item Section \ref{sec:Proof:GradientDescentEstimators} has estimation error bounds for classes of estimators constructed from gradient descent iterates (Proposition \ref{lem:GDBounds:Fine} and \ref{lem:GDBounds:Coarse}).
    \item Section \ref{sec:proof:Subopt:FastDecay} has results that compare the performance of ridge regression and gradient descent when the eigenvalues decay either exponentially or with a fast power law rate (Theorem \ref{thm:mainGDvsRidge_Upper:exp} and \ref{thm:mainGDvsRidge_Upper:poly}).
    \item Section \ref{sec:LogGridRidge} summarizes the relative sub-optimality between ridge regression tuned with a logarithmic grid and gradient descent when the eigenvalues decay at a slow power law rate (Theorem \ref{thm:log-ridgeVsGD}). 
\end{itemize}

\subsection{Proof of Lemma \ref{lem:ErrorDecomp}}
\label{sec:Proof:FirstLemma}
For simplicity of notation, we assume without loss of generality that $\psi=1$. This can always be achieved by rescaling the problem.
Performing a standard bias and variance decomposition and using the independence of $\epsilon$ and $\beta_{\star}$, recalling \eqref{hbphi}, gives
\begin{align*}
    & \E_{\beta_{\star},\epsilon}[L_{\beta_{\star}}(\widehat{\beta}_{\Phi})]
    = 
   \E_{\beta_{\star},\epsilon}[\|\widehat{\beta}_{\Phi} - \E_{\epsilon}[\widehat{\beta}_{\Phi}] \|_2^2]
    + 
    \E_{\beta_{\star}}[ \|\E_{\epsilon}[\widehat{\beta}_{\Phi}] - \beta_{\star} \|_2^2]\\
    & = \frac{\sigma^2}{n} \trace\left( \Phi\left(\frac{X^{\top} X}{n}\right) \frac{X^{\top} X}{n} \Phi\left(\frac{X^{\top} X}{n}\right)\right) \\
    &\quad\quad 
    + 
    \frac{1}{d} \trace\left( 
    \left(I - \Phi\left(\frac{X^{\top}X}{n}\right) \frac{X^{\top}X}{n} \right)
    \left(I - \Phi\left(\frac{X^{\top}X}{n}\right) \frac{X^{\top}X}{n}\right)^{\top} \right)\\
    & = 
      \frac{  d - r}{d}
    + 
     \frac{1}{d} \sum_{i=1}^{r} ( 1 - \Phi(s_i) s_i)^2
    + 
    \frac{\sigma^2}{n} 
    \sum_{i=1}^{r} \Phi(s_i)^2 s_i.
\end{align*}
Denoting $M_\Phi(s) = 1 - \Phi(s)s$ and $\zeta= 1/\lambda_{\star} = n/ (d \sigma^2)$, we then find 
\begin{align*}
     \E_{\beta_{\star},\epsilon}[L_{\beta_{\star}}(\widehat{\beta}_{\Phi})] 
    =
    \frac{ d - r}{d}
    + 
    \frac{\sigma^2}{n} \left[ 
    \sum_{i=1}^{r}
    \zeta M_\Phi(s_i)^2 + 
    \frac{(1 - M_\Phi(s_i))^2}{s_i}
    \right].
\end{align*}
Since for $s > 0 $
\begin{align*}
    \zeta M_\Phi(s)^2 + \frac{(1 - M_\Phi(s))^2}{s}
    & = 
    \frac{(1 + \zeta s)}{s} \left( M_\Phi(s) - \frac{1}{1 + \zeta s} \right)^2
    + 
    \frac{1}{1/\zeta + s},
\end{align*}
we find
\begin{align*}
     \E_{\epsilon,\beta_{\star}}[L_{\beta_{\star}}(\widehat{\beta}_{\Phi})]
    & =
    \frac{d - r}{d}
    + 
     \sigma^2  \frac{1}{n} \sum_{i=1}^{r} \frac{1}{1/\zeta + s_i}
    + 
     \frac{\sigma^2}{n} \sum_{i=1}^{r} 
    \frac{(1 + \zeta s_i )}{s_i} \left( M_\Phi(s_i) - \frac{1}{1 + \zeta s_i} \right)^2\\
    & = 
    \frac{d}{n} \sigma^2 \int \frac{1}{s+\lambda_{\star}} d\widehat{H}(s)
    + 
     \frac{\sigma^2}{n} \sum_{i=1}^{r} 
    \left(\frac{1}{s_i} + \frac{1}{\lambda_{\star}}\right) 
    \left( M_\Phi(s_i) - \frac{\lambda_{\star} }{\lambda_{\star} + s_i} \right)^2,\nonumber
\end{align*}
where we recall $\lambda_{\star}=1/\zeta$ and the probability measure $\widehat{H}(s)$ associated to the eigenvalues $\{s_i\}_{i=1}^{r}$ with a point mass of $(1-r/d)$ at zero. Thus Lemma \ref{lem:ErrorDecomp} holds because the rightmost term is non-negative, and zero when the estimator $\widehat{\beta}_{\Phi}$ equals ridge regression with penalization $\lambda= \lambda_{\star}$, i.e. $\Phi(s_i) = \frac{s_i}{s_i + \lambda_{\star}}$ for $i=1,\dots,r$. 

For a general $\psi$ possibly different from unity, we obtain with $\zeta = (n\psi)/(d\sigma^2)$ that
\begin{align}
     \E_{\epsilon,\beta_{\star}}[L_{\beta_{\star}}(\widehat{\beta}_{\Phi})]
    & -
       \int \frac{\lambda_{\star}}{s+\lambda_{\star}} d\widehat{H}(s)
    =
     \frac{\sigma^2}{n} \sum_{i=1}^{r} 
    \frac{(1 + \zeta s_i )}{s_i} \left( M_\Phi(s_i) - \frac{1}{1 + \zeta s_i} \right)^2.\label{errdec}
\end{align}

\subsection{Proof of Theorem \ref{thm:Informal}}
\label{sec:ProofInformal}

The proof of Theorem  \ref{thm:Informal} is split into two cases: \textbf{Slow Power Law Decay}, and \textbf{Fast Power Law or Exponential decay}, which we present in turn. 

\textbf{Slow Power Law Decay}.
Since $s_1=1$, the condition $k > \max\{s_1,2\}/\lambda_{\min}$ from Theorem \ref{thm:mainGDvsRidge_Upper} holds, thus from that result, if
\begin{align*}
    r \geq
    &
    M:=\max  \left\{ 
    \left(\frac{\lambda_{\min} k }{\lambda_\star^2}\right)^{2/(1-\alpha)} 
        \left( \frac{1}{\lambda_\star^{1/\alpha - 1}} 
    + \mathcal{J}_{\alpha,\lambda_{\star}}(k,\lambda_{\min})
    \right)^{1/(1-\alpha)},
    2^{1/(1-\alpha)}(1+\lambda_\star^{-1/\alpha})
    \right\}
\end{align*}
then we have the upper bound 
\begin{align*}
    \mathcal{S}\left( \mathcal{C}^{\text{GD}}(\eta,k), 
    \mathcal{C}^{\text{Ridge}}(\lambda_{\min},k)
    \right)
    & \leq 
    \frac{2304}{\min\{1-\ep,\ep\}^2(1 - \lambda_{\min})^2}
    \frac{\lambda_{\star}^4}{\lambda_{\min}^2}.
\end{align*}
Meanwhile from Theorem \ref{thm:mainGDvsRidge_Lower} we see that once 
\begin{align*}
    r \geq 
    \left(\frac{4}{1-\alpha}\right)^{1/(1-\alpha)} \left(\frac{4 k \lambda_{\min}}{\lambda_{\star}^2 \kappa }\right)^{1/\alpha}
\end{align*}
then we have the lower bound 
\begin{align*}
    \mathcal{S}\left( \mathcal{C}^{\text{GD}}(\eta,k), 
    \mathcal{C}^{\text{Ridge}}(\lambda_{\min},k)
    \right)
    & \geq
    \frac{1}{2\cdot512} 
    \frac{\min\{ 1- \kappa,\kappa \}^2}{\max\{ 1-\ep, \ep \}^2 (1-\lambda_{\min})^2}
    \frac{\lambda_{\star}^4}{\lambda_{\min}^2}.
\end{align*}
It therefore suffices to take
\begin{align*}
    & L_{\alpha,\lambda_{\star},k,\lambda_{\min}} = 
    \max\left\{ M,
    \left(\frac{4}{1-\alpha}\right)^{1/(1-\alpha)} \left(\frac{4 k \lambda_{\min}}{\lambda_{\star}^2 \kappa }\right)^{1/\alpha}
    \right\}
\end{align*}
to have matching upper and lower bounds up to constant factors with 
\begin{align*}
\mathcal{S}\left( \mathcal{C}^{\text{GD}}(\eta,k), \mathcal{C}^{\text{Ridge}}(\lambda_{\min},k)\right) \simeq 
\frac{1}{ \min\{1-\ep,\ep\}^2 (1-\lambda_{\min})^2}\frac{\lambda_{\star}^4}{\lambda_{\min}^2}.    
\end{align*}
From $\sigma_{\min} \leq \sigma/\sqrt{2}$ we have that $\lambda_{\min} = d \sigma_{\min}^2/n \leq d \sigma^2/(2n) = \lambda_{\star}/2 \leq 1/2$ and therefore $1\geq (1-\lambda_{\min})^2 \geq 1/4$. 
Meanwhile, since $\lambda_{\star} = \lambda_{\min}+ (j+\ep)\delta$ for an integer $k-2 \geq j \geq 0$ and $\ep \in (0,1)$, we have $\text{Dist}(\lambda_{\star},\Gamma) = \delta \min\{\ep,1-\ep\}$ and as such $\text{Dist}_{\delta}(\lambda_{\star},\Gamma) = \text{Dist}(\lambda_{\star},\Gamma)/\delta = \min\{\ep,1-\ep\}$.
This then gives us $\mathcal{S}\left(\mathcal{C}^{\text{GD}}(\eta,k), \mathcal{C}^{\text{Ridge}}(\lambda_{\min},k)\right) \simeq \text{Dist}_{\delta}(\lambda_{\star},\Gamma)^{-2} \lambda_{\star}^4/\lambda_{\min}^2 $ as required.

\textbf{Fast Power Law or Exponential Decay}: 
We now consider the case when the eigenvalues decay more quickly, beginning with the case of eigenvalues that decay exponentially, so $s_i = \exp(-\rho(i-1))$ for some $\rho > 0$ for $i=1,\dots,r$. The condition on $k$ from \eqref{expgdc} in Theorem \ref{thm:mainGDvsRidge_Upper:exp} is satisfied by assumption, whereas the upper bound on $\lambda_{\star}$ from \eqref{expgdc} is on the order of a constant. 
The condition $r \geq 1  + \rho^{-1}\log(1/\lambda_{\star})$  in Theorem \ref{thm:mainGDvsRidge_Upper:exp} holds due to the assumption $s_r \leq \lambda_{\star}$. Given this we have from Theorem \ref{thm:mainGDvsRidge_Upper:exp} the lower bound 
\begin{align}
\label{equ:Lowerbound:Tmp}
    \mathcal{S}\left( \mathcal{C}^{\text{GD}}(\eta,k), 
    \mathcal{C}^{\text{Ridge}}(\lambda_{\min},k)
    \right)
    & \geq  
    \frac{C_{\rho}}{\min\{1-\ep,\ep\}^2}
    \left( \frac{ \lambda_{\star} k}{1-\lambda_{\min}}\right)^2
\end{align}
for some constant $C_{\rho} >0$ depending exponentially $\rho$, i.e., $C_{\rho} =\Theta( \exp(-\rho))$, and independent of $n,p,k,\sigma$. Following the previous case, going from $\min\{1-\ep,\ep\}$ to $\text{Dist}_{\delta}(\lambda_{\star},\Gamma)$, we arrive at the lower bound 
$\mathcal{S}\left( \mathcal{C}^{\text{GD}}(\eta,k), \mathcal{C}^{\text{Ridge}}(\lambda_{\min},k)\right)  \gtrsim \text{Dist}_{\delta}(\lambda_{\star},\Gamma)^{-2}(k\lambda_\star)^2$ as $k$ grows, provided $\rho \lesssim \log(k)$.

Let us now consider eigenvalues that decay at a faster power law rate so $s_i = i^{-\alpha}$ for $\alpha > 1$ and $i=1,\dots,r$. From Theorem \ref{thm:mainGDvsRidge_Upper:poly} we see once again that the condition \eqref{thm6cond1}
on $k$ is implied by the assumption of the theorem,
and the upper bound \eqref{thm6cond2} on $\lambda_{\star}$ is on the order of a constant. Similarly, $r \geq \lambda_{\star}^{-1/\alpha}$  holds due to $s_r \leq \lambda_{\star}$. Given this, the lower bound matches \eqref{equ:Lowerbound:Tmp} for a (different) constant $C_{\alpha}> 0$ depending expoentially on $\alpha$ so $C_{\alpha} = \Theta(\exp(-\alpha))$, and thus, we have the desired lower bound provided $\alpha \lesssim \log(k)$.

\subsection{Proofs for Relative Sub-optimality with Slow Power Law Decaying Eigenvalues}
\label{sec:Proof:Subopt:Slow}
In this section we present proofs for the results associated to bounding the relative sub-optimality of gradient descent and ridge regression when the eigenvalues decay at a slow power law rate. Section  \ref{sec:proof:mainGDvsRidge_Upper} presents the proof of Theorem \ref{thm:mainGDvsRidge_Upper} which gives an upper bound, meanwhile Section \ref{sec:prooof:mainGdvsRidge_Lower} presents the proof of Theorem \ref{thm:mainGDvsRidge_Lower} which gives the lower bound. 

\subsubsection{Proof of Theorem \ref{thm:mainGDvsRidge_Upper}}
\label{sec:proof:mainGDvsRidge_Upper}
Plugging the expression for $\delta$ in Proposition \ref{lem:RidgeBounds} from Section \ref{b-ind}, with $\delta = (1-\lambda_{\min})/(k-1)$ (where the condition $\lambda_\star \ge \delta$ holds because $\lambda_\star \ge \lambda_{\min} \ge (1-\lambda_{\min})/(k-1) = \delta$, as $\lambda_{\min} > 1/k$) immediately yields the lower bound for ridge regression 
\begin{align*}
    \mathcal{E}(\mathcal{C}^{\text{Ridge}}(\lambda_{\min},k)) 
    \geq \frac{1}{2^{5}} \min\{ 1-\ep,\ep\}^2\lambda_{\star} 
    \int G^{\text{Ridge}}(s) d \widehat{H}(s),
\end{align*}
where 
\begin{align}
    \int G^{\text{Ridge}}(s) d \widehat{H}(s) 
    & = 
    \frac{(1 - \lambda_{\min})^2}{(k-1)^2}
    \frac{1}{d} 
    \left( 
    \sum_{i : s_i > \lambda_{\star}}
    \frac{1}{\lambda_{\star} s_i^2} 
    + 
    \sum_{i : s_i \leq \lambda_{\star}}
    \frac{s_i}{\lambda_{\star}^4}
    \right) \nonumber 
    \\
    & \geq 
    \frac{(1 - \lambda_{\min})^2}{k^2}
    \frac{1}{d} 
    \frac{1}{\lambda_{\star}^4}
    \left( 
    \lambda_{\star}^3 
    \sum_{i : s_i > \lambda_{\star}}
    \frac{1}{s_i^2}
    + 
    \sum_{i : s_i \leq \lambda_{\star}} s_i 
    \right).
    \label{equ:RidgeIntegral:calc}
\end{align}
We now turn to gradient descent. Since $\eta = \frac{1}{k \lambda_{\min}}$, using Proposition \ref{lem:GDBounds:Fine} with $t=k$, where the condition $t \ge \lceil 1/(\eta \lambda_\star)\rceil$ holds because $t = k = 1/(\eta \lambda_{\min}) \ge \lceil1/(\eta \lambda_\star)\rceil$, we find 
\begin{align*}
    \mathcal{E}(\mathcal{C}^{\text{GD}}(\eta,t)) 
    \leq 
    18 \lambda_{\star}\int G^{\text{GD}}(s) d \widehat{H}(s),
\end{align*}
where 
\begin{align}
\label{equ:GDSpectralInt}
    \int G^{\text{GD}}(s) d \widehat{H}(s)
    & = 
    \frac{1}{d}
    \left( \sum_{i : s_i > \lambda_{\star}} \frac{\lambda_{\star}}{s_i^2}
    + 
    \sum_{i : \lambda_{\star}^2/(k \lambda_{\min}) < s_i \leq \lambda_{\star} }
    \frac{s_i^3}{\lambda_{\star}^4}
    + 
    \frac{1}{k^2 \lambda_{\min}^2} 
    \sum_{i: s_i \leq \lambda_{\star}^2/(k \lambda_{\min})}
    s_i
    \right)\\
    \nonumber
    & = 
    \frac{1}{d} 
    \frac{1}{\lambda_{\star}^4}
    \left( \lambda_{\star}^3 \sum_{i : s_i > \lambda_{\star}} 
    \frac{1}{s_i^2}
    + 
    \sum_{i : \lambda_{\star}^2/(k \lambda_{\min}) < s_i \leq \lambda_{\star} }
    s_i^3
    + 
    \frac{\lambda_{\star}^4}{k^2 \lambda_{\min}^2} 
    \sum_{i: s_i \leq \lambda_{\star}^2/(k \lambda_{\min})}
    s_i
    \right).
\end{align}
Dividing this by the bound in \eqref{equ:RidgeIntegral:calc} gives, using $\lambda_{\star}^2/(k \lambda_{\min}) \le \lambda_{\star}$ 
\begin{align*}
    & \frac{\int G^{\text{GD}}(s) d \widehat{H}(s) }{
    \int G^{\text{Ridge}}(s) d \widehat{H}(s) 
    }  \\
    & \le 
    \frac{1}{(1-\lambda_{\min})^2}
    \left[
    k^2 
    \frac{ \lambda_{\star}^3 
    \sum_{i : s_i > \lambda_{\star}}
    \frac{1}{s_i^2}}{
    \left( 
    \lambda_{\star}^3 
    \sum_{i : s_i > \lambda_{\star}}
    \frac{1}{s_i^2}
    + 
    \sum_{i : s_i \leq \lambda_{\star}} s_i 
    \right)
    }\right.\\
    & \quad\quad\quad\quad\quad\quad \left.+  
    k^2 \frac{\sum_{i : \lambda_{\star}^2/(k \lambda_{\min}) < s_i \leq \lambda_{\star} } 
    s_i^3
    }{ 
    \left( 
    \lambda_{\star}^3 
    \sum_{i : s_i > \lambda_{\star}}
    \frac{1}{s_i^2}
    + 
    \sum_{i : s_i \leq \lambda_{\star}} s_i 
    \right)
    } 
    + 
    \frac{\lambda_{\star}^4}{ \lambda_{\min}^2} 
    \frac{\sum_{i: s_i \leq \lambda_{\star}^2/(k \lambda_{\min})}
    s_i}{
    \left( 
    \lambda_{\star}^3 
    \sum_{i : s_i > \lambda_{\star}}
    \frac{1}{s_i^2}
    + 
    \sum_{i : s_i \leq \lambda_{\star}} s_i 
    \right)
    }
    \right] \\
    & \leq 
    \frac{1}{(1-\lambda_{\min})^2}
    \left[
    k^2 
    \frac{ \lambda_{\star}^3 
    \sum_{i : s_i > \lambda_{\star}}
    \frac{1}{s_i^2}}{
    \sum_{i : s_i \leq \lambda_{\star}} s_i 
    }
    +  
    k^2 \frac{\sum_{i : \lambda_{\star}^2/(k \lambda_{\min}) < s_i \leq \lambda_{\star} } 
    s_i^3
    }{ 
    \sum_{i : s_i \leq \lambda_{\star}} s_i 
    }
    +  \frac{\lambda_{\star}^4}{\lambda_{\min}^2}
    \right].
\end{align*}
Since $s_i = i^{-\alpha}$ for $i=1,\dots,r$, we have the following bounds:
\begin{align}
    \sum_{i : s_i > \lambda_{\star}} \frac{1}{s_i^2} 
    & \leq \lambda_{\star}^{-(2+1/\alpha)} \label{equ:PartialSumBounds1} \\
    \sum_{i : s_i \leq \lambda_{\star}} s_i & \geq 
    \frac{r^{1-\alpha}}{2}
    \label{equ:PartialSumBounds2}
    \\
    \sum_{i : \lambda_{\star}^2/(k \lambda_{\min}) < s_i \leq \lambda_{\star} } 
    s_i^3 & \leq 
    \mathcal{J}_{\alpha,\lambda_\star}(k,\lambda_{\min}),
    \label{equ:PartialSumBounds3}
\end{align}
where $\mathcal{J}_{\alpha,\lambda_\star}$ is defined in \eqref{J}.
Plugging in we then find
\begin{align*}
    & \frac{\int G^{\text{GD}}(s) d \widehat{H}(s) }{
    \int G^{\text{Ridge}}(s) d \widehat{H}(s) 
    }
     \leq 
     \frac{2}{(1-\lambda_{\min})^2}
     \left[ 
     \frac{k^2}{r^{1-\alpha}} 
     \left( \frac{1 }{\lambda_\star^{1/\alpha - 1}} 
     + \mathcal{J}_{\alpha,\lambda_\star}(k,\lambda_{\min}) 
     \right)
     + 
     \frac{\lambda_{\star}^{4}}{\lambda_{\min}^2}
     \right]
\end{align*}
as required.

We now show the bounds \eqref{equ:PartialSumBounds1}, \eqref{equ:PartialSumBounds2}, \eqref{equ:PartialSumBounds3}. For \eqref{equ:PartialSumBounds1} we have 
\begin{align*}
    \sum_{i : s_i > \lambda_{\star}} \frac{1}{s_i^2}
    = 
    \sum_{i<\lambda_{\star}^{-1/\alpha}} i^{2\alpha}
    \leq  
    \sum_{i=1}^{\lambda_{\star}^{-1/\alpha}} i^{2\alpha}
    \leq 
    \lambda_{\star}^{-1/\alpha} \lambda_\star^{-2}
    = 
    \lambda_{\star}^{-(2+1/\alpha)}.
\end{align*}
For \eqref{equ:PartialSumBounds2}  we have from $i^{-\alpha} \geq \int_{i}^{i+1} x^{-\alpha} dx$ that
\begin{align*}
    \sum_{i:s_i \leq \lambda_\star} s_i 
    & = 
    \sum_{i=\lceil \lambda_\star^{-1/\alpha} \rceil }^{r} i^{-\alpha}
     \geq \int_{\lceil \lambda_\star^{-1/\alpha} \rceil}^{r+1} x^{-\alpha}dx\\
    & = 
    \frac{1}{1-\alpha}\left( (r+1)^{1-\alpha} - (\lceil \lambda_\star^{-1/\alpha}\rceil )^{1-\alpha} \right)\\
    & \geq 
    \frac{1}{1-\alpha} \left( r^{1-\alpha} - (\lceil \lambda_\star^{-1/\alpha}\rceil )^{1-\alpha} \right)
     \geq 
    \frac{r^{1-\alpha}}{2(1-\alpha)},
\end{align*}
for the last inequality we used the assumption that $r \geq 2^{1/(1-\alpha)} (1+\lambda_\star^{-1/\alpha})$. 

\sloppy Finally for \eqref{equ:PartialSumBounds3}  we have with $a=\lambda_{\star}^2/(k\lambda_{\min})$, $b= \lambda_{\star}$ as well as  with $i^{-3\alpha} \leq \int_{i-1}^{i} x^{-3\alpha} dx$ that
\begin{align*}
    \sum_{i : a < s_i \leq b } 
    s_i^3
    & = 
    \sum_{i= \lceil b^{-1/\alpha} \rceil }^{ \lceil a^{-1/\alpha} \rceil -1 }
    i^{-3\alpha} 
    = 
    \lceil b^{-1/\alpha} \rceil^{-3\alpha} 
    + 
    \sum_{i= \lceil b^{-1/\alpha} \rceil +1 }^{ \lceil a^{-1/\alpha} \rceil -1 }
    i^{-3\alpha}\\
    & \leq \lceil b^{-1/\alpha} \rceil^{-3\alpha} +
 \int_{\lceil b^{-1/\alpha} \rceil}^{\lceil a^{-1/\alpha} \rceil -1 }
    x^{-3\alpha}dx
     \leq \lceil b^{-1/\alpha} \rceil^{-3\alpha} +
  \int_{b^{-1/\alpha}}^{a^{-1/\alpha}}
    x^{-3\alpha}dx\\
    & \leq 
    1
    + 
    \begin{cases}
    3 \log(k \lambda_{\min} /\lambda_\star)
    & \text{ if } \alpha = 1/3 \\
     \frac{\lambda_\star^{3-1/\alpha}}{1-3\alpha}
     \left[ \left( \frac{\lambda_\star}{k \lambda_{\min}} \right)^{3-1/\alpha} - 1\right] & \text{ if } \alpha \not= 1/3
    \end{cases}\\
    & \leq 
    \mathcal{J}_{\alpha,\lambda_\star}(k,\lambda_{\min}),
\end{align*}
where we used that $\lceil \lambda_\star^{-1/\alpha} \rceil^{-3\alpha}   \leq \lambda_\star^3 \le 1$. This yields the desired inequality.

\subsubsection{Proof of Theorem \ref{thm:mainGDvsRidge_Lower}}
\label{sec:prooof:mainGdvsRidge_Lower}
Following the steps in Theorem \ref{thm:mainGDvsRidge_Upper} but applying the upper bound for ridge regression from Proposition \ref{lem:RidgeBounds} with $\delta= (1-\lambda_{\min})/(k-1)$ gives 
\begin{align}\label{r-ub}
    \mathcal{E}(\mathcal{C}^{\text{Ridge}}(\lambda_{\min},k))
    \leq 
    4 \min\left\{ 1-\ep, \ep \right\}^2 
    \lambda_\star
    \int G^{\text{Ridge}}(s) d \widehat{H}(s).
\end{align}
Moreover \eqref{equ:PartialSumBounds1}, \eqref{equ:PartialSumBounds2} and the condition $r \geq 2^{1/(1-\alpha)}(1+\lambda_\star^{-1/\alpha})$ 
yield
\begin{align}
\label{equ:eigenval_tail_lower_bound}
    \sum_{i: s_i \leq \lambda_{\star}} s_i 
    \geq  \frac{r^{1-\alpha}}{2} 
    \geq
    \lambda_{\star}^{1-1/\alpha}
    = 
    \lambda_{\star}^{3}
    \lambda_{\star}^{-(2+1/\alpha)}
     \geq 
    \lambda_{\star}^{3}\sum_{i : s_i > \lambda_{\star}} \frac{1}{s_i^2}.
\end{align}
Therefore 
\begin{align*}
    \int G^{\text{Ridge}}(s) d \widehat{H}(s)
    \leq 
    \frac{(1- \lambda_{\min})^2}{(k-1)^2}
    \frac{1}{d} \frac{1}{\lambda_{\star}^4}
    2 
    \sum_{i : s_i \leq \lambda_{\star}} s_i.
\end{align*}

For gradient descent, Proposition \ref{lem:GDBounds:Fine} yields the lower bound
\begin{align}
\label{equ:GD:lower_calc_1}
    \mathcal{E}(\mathcal{C}^{\text{GD}}(\eta,k))
    & \geq 
        \frac{1}{16} 
    \min\{1-\kappa,\kappa\}^2
    \lambda_{\star}
    \int_{0}^{\eta \lambda_{\star}^2 \kappa/4} G^{\text{GD}}(s) d \widehat{H}(s).
\end{align}
Plugging in that $\eta = \frac{1}{k \lambda_{\min}}$ as well as that $G^{\text{GD}}(s) = \eta^2 s$  for $s \leq \eta \lambda_\star^2 $ by definition \eqref{ggd}  gives 
\begin{align}
\label{equ:GD:lower_calc_2}
    \int_{0}^{\eta \lambda_{\star}^2 \kappa/4} G^{\text{GD}}(s) d \widehat{H}(s)
    = 
    \frac{1}{d} 
    \sum_{i : s_i \leq \eta \lambda_\star^2 \kappa/4 } 
    \eta^2 s_i
    = \left( \frac{1}{\lambda_{\min} k} \right)^2
    \frac{1}{d}
    \sum_{i : s_i \leq \frac{\lambda_{\star}^2 \kappa }{4 k \lambda_{\min}} }
    s_i.
\end{align}
Taking the ratio with the ridge regression bound and plugging in the above bounds yields 
\begin{align*}
    &\mathcal{S}(\mathcal{C}^{\text{GD}}(\eta,k),\mathcal{C}^{\text{Ridge}}(\lambda_{\min},k))
    \geq 
    \frac{1}{64} \frac{\min\{1-\kappa,\kappa\}^2}{\min\{1-\ep,\ep\}^2} \frac{\int_{0}^{\eta \lambda_{\star}^2 \kappa/4} G^{\text{GD}}(s) d \widehat{H}(s)}{\int G^{\text{Ridge}}(s) d \widehat{H}(s)}\\
    & \geq 
    \frac{1}{128} \frac{\min\{1-\kappa,\kappa\}^2}{\min\{1-\ep,\ep\}^2} 
    \left(1 - \frac{1}{k}\right)^2 
    \frac{1}{(1-\lambda_{\min})^2}
    \frac{\lambda_{\star}^4}{\lambda_{\min}^2} 
    \frac{ \sum_{i : s_i \leq \frac{\lambda_{\star}^2 \kappa }{4 k \lambda_{\min}} }
    s_i}{ 
    \sum_{i : s_i \leq \lambda_{\star}} s_i
    }\\
    & \geq 
    \frac{1}{512} 
    \frac{\min\{ 1- \kappa,\kappa \}^2}{\min\{ 1-\ep, \ep \}^2}
    \frac{1}{(1-\lambda_{\min})^2}
    \frac{\lambda_{\star}^4}{\lambda_{\min}^2}
    \frac{ \sum_{i : s_i \leq \frac{\lambda_{\star}^2 \kappa }{4 k \lambda_{\min}} }
    s_i}{ 
    \sum_{i : s_i \leq \lambda_{\star}} s_i
    }.
\end{align*}
We bound
\begin{align*}
    \sum_{i : s_i \leq \frac{\lambda_{\star}^2 \kappa }{4 k \lambda_{\min}} }
    s_i
    & = 
    \sum_{i : s_i \leq \lambda_{\star}} s_i 
    - 
    \sum_{i : \lambda_{\star} \geq s_i >
    \frac{\lambda_{\star}^2 \kappa }{4 k \lambda_{\min}}} 
    s_i 
     \geq
    \sum_{i : s_i \leq \lambda_{\star}} s_i 
    - 
    \sum_{i :  s_i >
    \frac{\lambda_{\star}^2 \kappa }{4 k \lambda_{\min}}} 
    s_i \\
    & \geq 
    \sum_{i : s_i \leq \lambda_{\star}} s_i  
    - 
    \frac{1}{1-\alpha} \left( \frac{4 k \lambda_{\min}}{\lambda_{\star}^2 \kappa }\right)^{1/\alpha - 1},
\end{align*}
where we have used that  $\sum_{i : s_i > b} s_i = \sum_{i =1}^{\lceil b^{-1/\alpha}\rceil-1 } i^{-\alpha} \leq 1+ \int_1^{b^{-1/\alpha}} x^{-\alpha} dx \leq \frac{1}{1-\alpha} b^{1-1/\alpha}$.

Therefore, using \eqref{equ:PartialSumBounds2}  we find 
\begin{align}
\label{equ:series_ratio_lower_bound}
    \frac{ \sum_{i : s_i \leq \frac{\lambda_{\star}^2 \kappa }{4 k \lambda_{\min}} }
    s_i}{ 
    \sum_{i : s_i \leq \lambda_{\star}} s_i
    }
    \geq 
    1 - 
    \frac{2}{1-\alpha} \left( \frac{4 k \lambda_{\min}}{\lambda_{\star}^2 \kappa }\right)^{1/\alpha - 1}
    \frac{1}{r^{1-\alpha}}.
\end{align}
This then leads to the required lower bound.

\subsection{Limits of the Relative Sub-optimality}
\label{sec:sub-opt-limit}
In this section we present and prove Corollary \ref{RiskLimit:Formal} and \ref{RiskLimit:Lower:Formal} which are summarized within the main contributions. We begin by presenting both corollaries, with their proofs given after. We begin with Corollary \ref{RiskLimit:Formal}, for eigenvalues with a slow power law decay.

\begin{corollary}
\label{RiskLimit:Formal}
Consider the setting of Theorem \ref{thm:mainGDvsRidge_Upper}. Let $r  = \min(d,n)$ and $d,n \rightarrow \infty$ so that $\gamma = d/n$ is fixed, and $\lambda_{\min,k} = \lambda_{\star,k}/2$. There is a constant $C_{\lambda_{\star,k}}$ depending on the discretization we have 
\begin{align*}
    \limsup_{\substack{ d,n \rightarrow \infty \\ d/n = \gamma}}\mathcal{S}\big( \mathcal{C}^{\text{GD}}(\eta,k), 
    \mathcal{C}^{\text{Ridge}}(\lambda_{\min},k)
    \big)
\lesssim C_{\lambda_{\star,k}} \lambda_{\star,k}^2.
\end{align*}
When there is $\delta \in (0,1]$ and $\epsilon \in (0,1)$ such that for all $k$, there is $j_k \in \{0,1\dots,\lfloor (k-1)^{1-\delta} \rfloor \}$ such that
$\lambda_{\star, k} = 2 (\ep + j_k)/(k-1 + \ep + j_k) $,  then $\lim_{k \rightarrow \infty} C_{\lambda_{\star, k}} \lambda_{\star, k}^2 = 0$. 
\end{corollary}

We now present Corollary \ref{RiskLimit:Lower:Formal}, for eigenvalues that decay quickly, i.e., exponentially, or with a fast power law decay. 
\begin{corollary}
\label{RiskLimit:Lower:Formal}
Consider the setting of Theorem \ref{thm:mainGDvsRidge_Upper:exp} (exponentially decaying eigenvalues) or Theorem \ref{thm:mainGDvsRidge_Upper:poly} (eigenvalues with a fast power law decay). Let $r  = \min(d,n)$ and $d,n \rightarrow \infty$ so that $\gamma := d/n$ is fixed, and $\lambda_{\min,k} = \lambda_{\star,k}/2$. There is a constant $C_{\lambda_{\star,k}}$ depending on the discretization, such that we have 
\begin{align*}
    \liminf_{\substack{ d,n \rightarrow \infty \\ d/n = \gamma}}\mathcal{S}\big( \mathcal{C}^{\text{GD}}(\eta,k), 
    \mathcal{C}^{\text{Ridge}}(\lambda_{\min},k)
    \big)
\gtrsim C_{\lambda_{\star,k}} \lambda_{\star,k}^2 k^2.
\end{align*}
When there is $\delta^{\prime} \in (0,1/2)$, $\delta^{\prime\prime} \in (0, 1/2-\delta^{\prime})$ and $\epsilon \in (0,1)$ such that for all $k$, there is $j_k \in \{\lceil (k-1)^{1/2 + \delta^{\prime}}\rceil \dots, \lfloor (k-1)^{1-\delta^{\prime\prime}}\rfloor  \}$ such that
$\lambda_{\star, k} = 2 (\ep + j_k)/(k-1 + \ep + j_k) $,  then $\lim_{k \rightarrow \infty} C_{\lambda_{\star, k}} \lambda_{\star, k}^2 k^2 = \infty$.

\end{corollary}

\subsubsection{Proof of Corollary \ref{RiskLimit:Formal}}
In the case $r = \min(d,n)$ and $d,n \rightarrow \infty$, we immediately have $r \rightarrow \infty$. The upper bound from Theorem \ref{thm:mainGDvsRidge_Upper} then simplifies to
\begin{align*}
    \limsup_{\substack{ d,n \rightarrow \infty \\ d/n = \gamma}} \mathcal{S}\left( \mathcal{C}^{\text{GD}}(\eta,k), 
    \mathcal{C}^{\text{Ridge}}(\lambda_{\min},k)
    \right)
    & \leq 
    \frac{1152}{\min\{1-\ep,\ep\}^2(1 - \lambda_{\min})^2}
    \frac{\lambda_{\star}^4}{\lambda_{\min}^2},
\end{align*}
where we recall that $\lambda_{\star} = \lambda_{\min} + (j+\ep)\delta$, with $\epsilon \in (0,1)$, $j \in \{0,1\dots,k-2\}$ and $\delta = (1-\lambda_{\min})/(k-1)$.  
For $\epsilon \in (0,1)$ and  $k \geq 2$,  consider the sequence for $\lambda_{\star,k} = \lambda_{\min} + (j_{k} + \ep)(1-\lambda_{\min})/(k-1)$ where $j_{k} \in \{0,1,2,\dots,\lceil (k-1)^{1-\delta}\rceil \}$ with $\delta \in (0,1]$. 
When  $\lambda_{\min,k} = \lambda_{\star,k}/2 $, we have
\begin{align*}
    \lambda_{\star, k} & = \lambda_{\min} + (\ep + j_k) \frac{(1-\lambda_{\min})}{k-1}
    =  \frac{\lambda_{\star, k}}{2} \Big(1 - \frac{\ep +j_k}{k-1} \Big) + \frac{\ep + j_k}{(k-1)} 
\end{align*}
and thus
\begin{align*}
    \lambda_{\star, k}  = 2 \frac{\frac{\ep + j_k }{k-1}}{ 1 + \frac{\ep + j_k}{k-1}}
    = 2 \frac{ \ep + j_k}{k-1 + \ep + j_k}.
\end{align*}
Plugging in this sequence we find since $1-\lambda_{\min} \geq 1/2$
\begin{align*}
    \limsup_{\substack{ d,n \rightarrow \infty \\ d/n = \gamma}} \mathcal{S}\left( \mathcal{C}^{\text{GD}}(\eta,k), 
    \mathcal{C}^{\text{Ridge}}(\lambda_{\min},k)
    \right)
    \leq \frac{73728}{\min\{1-\ep,\ep\}^2 } \Big( 
     \frac{ \ep + j_k}{k-1 + \ep + j_k}
    \Big)^2.
\end{align*}
Taking the limit as $k \rightarrow \infty$ we then see that $j_k/(k-1) \leq (k-1)^{1-\delta}/(k-1) = (k-1)^{-\delta} \rightarrow 0$ as required.

\subsubsection{Proof of Corollary \ref{RiskLimit:Lower:Formal}}
Let us begin by considering Theorem \ref{thm:mainGDvsRidge_Upper:exp}, for eigenvalues that decay exponentially. With $r = \min(d,n)$ and $n,d \rightarrow \infty$ we have 
\begin{align}
\label{equ:fast:limit}
    \liminf_{\substack{ d,n \rightarrow \infty \\ d/n = \gamma}}\mathcal{S}\big( \mathcal{C}^{\text{GD}}(\eta,k), 
    \mathcal{C}^{\text{Ridge}}(\lambda_{\min},k)
    \big)
    \geq \frac{C_{\rho}}{\min\{1-\ep,\ep\}^2} 
    \left( \frac{\lambda_{\star}k}{1-\lambda_{\min}}\right)^2,
\end{align}
where we recall that $\lambda_{\star} = \lambda_{\min} + (j+\ep)\delta$ with $\epsilon \in (0,1)$, $j \in \{0,1\dots,k-2\}$ and $\delta = (1-\lambda_{\min})/(k-1)$.  Recalling that $\lambda_{\min}  = \lambda_{\star}/2$, for condition \eqref{expgdc} in Theorem \ref{thm:mainGDvsRidge_Upper:exp} to be satisfied it is sufficient that
\begin{align}
\label{equ:fast:limit:conditions}
    \lambda_{\star} \leq C_{\rho}^{\prime}
    \quad\quad \text{and} \quad\quad 
    k \geq \frac{C^{\prime\prime}_{\rho}}{\lambda_{\star}^2}, 
\end{align}
for some constants $C_{\rho}^{\prime},C_{\rho}^{\prime\prime} > 0 $. 

Following the proof of Corollary \ref{RiskLimit:Formal}, let us now consider the sequence 
$\lambda_{\star, k} =\lambda_{\min} + (j_{k} + \ep)\delta =  2 (\ep + j_k)/(k-1 + \ep + j_k) $ for some $\epsilon \in (0,1)$ and   $j_k \in \{\lceil (k-1)^{1/2 + \delta^{\prime}}\rceil \dots, \lfloor (k-1)^{1-\delta^{\prime\prime}}\rfloor  \}$ with $\delta^{\prime} \in (0,1/2)$ and $\delta^{\prime\prime} \in (0, 1/2-\delta^{\prime})$. We need to show that the conditions in \eqref{equ:fast:limit:conditions} hold for this sequence once $k$ is sufficiently large, that is,  
\begin{align*}
\lambda_{\star,k} = \frac{\ep + j_{k}}{\ep + k + j_{k}} \leq C_{\rho}^{\prime} 
\quad\quad \text{and} \quad\quad 
    k \geq \frac{C_{\rho}^{\prime\prime}}{\lambda_{\star,k}^2} = \frac{C^{\prime\prime}_\rho}{4} \Big( \frac{k-1}{\ep + j_{k}} + 1\Big)^2.
\end{align*}
From $j_{k} \leq (k-1)^{1-\delta^{\prime\prime}}$ we see that $\lambda_{\star,k} = (\ep + j_{k})/(\ep + j_{k} + k-1) \leq (\ep + (k-1)^{1-\delta^{\prime\prime}})/(k-1) \leq 2(k-1)^{-\delta^{\prime\prime}} $, and thus the first condition is satisfied for $k \geq 1 + (2/C_{\rho}^{\prime})^{1/\delta^{\prime\prime}}$. 
Since $j_{k} \geq (k-1)^{1/2+\delta^{\prime}}$ for $\delta^{\prime} \in (0,1/2)$, we have for sufficiently large $k$  that
\begin{align*}
    \frac{k-1}{\ep + j_{k}} + 1
    \leq 2(k-1)^{1/2-\delta^{\prime}}.
\end{align*}
Thus, for $k \geq (C_{\rho}^{\prime\prime})^{1/(2\delta^{\prime})}$ we have
\begin{align*}
    \frac{C^{\prime\prime}_\rho}{4} \Big( \frac{k-1}{\ep + j_{k}} + 1\Big)^2
    \leq 
    C_{\rho}^{\prime\prime}(k-1)^{1 - 2 \delta^{\prime}} \leq C_{\rho}^{\prime\prime}k^{1 - 2 \delta^{\prime}}
    \leq k
\end{align*}
as required for the second condition. 

Using that $1-\lambda_{\min} \geq 1/2$ and taking the limit as $k\to\infty$ we find
\begin{align*}
    \liminf_{k \rightarrow \infty} \lim_{\substack{ d,n \rightarrow \infty \\ d/n = \gamma}}\mathcal{S}\big( \mathcal{C}^{\text{GD}}(\eta,k), 
    \mathcal{C}^{\text{Ridge}}(\lambda_{\min},k)
    \big)
    \geq 4 \frac{C_{\rho}}{\min\{1-\ep,\ep\}^2} 
    \lim_{k \rightarrow \infty}\big( \lambda_{\star, k}, k\big)^2 
    = \infty,
\end{align*}
where in the last step we used that, since $k-1 \geq j_{k} \geq (k-1)^{1/2+ \delta^{\prime}}$ and $\ep \in (0,1)$, we have
\begin{align*}
    \lambda_{\star,k} k 
    = 
   2 k \frac{ \ep + j_k}{k-1 + \ep + j_k} 
    \geq 
    k \frac{j_{k}}{k + j_{k}}
    \geq k \frac{ (k-1)^{1/2 + \delta^{\prime}} }{k + k-1   }
    \geq \frac{1}{2} (k-1)^{1/2+\delta^{\prime}}.
\end{align*}

We now consider Theorem \ref{thm:mainGDvsRidge_Upper:poly}, for eigenvalues with a fast power law decay. The proof now follows the same set of steps as above, but with $C_{\rho},C_{\rho}^{\prime},C_{\rho}^{\prime\prime}$ swapped for certain $C_{\alpha},C_{\alpha}^{\prime},C_{\alpha}^{\prime\prime}$. Specifically, in the limit as $d,n \rightarrow \infty$ with $d/n = \gamma $ we arrive at the same limit \eqref{equ:fast:limit} but with the constant $C_{\rho}$ swapped for $C_{\alpha}$. Similarly, conditions \eqref{thm6cond1} and \eqref{thm6cond2} in Theorem  \ref{thm:mainGDvsRidge_Upper:poly} are satisfied when \eqref{equ:fast:limit:conditions} holds with appropriate constants $C_{\alpha}^{\prime}, C_{\alpha}^{\prime\prime} > 0$.

\subsection{Proof of Corollary \ref{cor:LowDimensional}}
\label{sec:proof:LowDimensional}
Since 
$r = d = n^{q}$ for $q \in (0,1)$, $\lambda_{\star} = \sigma^2 / n^{1-q}$. Moreover, as $q > 1/(1+\alpha)$, the condition  $n^{q} = r \geq 2^{1/(1-\alpha}(1 + \lambda_{\star}^{-1/\alpha}) = 2^{1/(1-\alpha)}(1 + \sigma^{-2/\alpha} n^{(1-q)/\alpha}) $ in Theorem \ref{thm:mainGDvsRidge_Upper} is then satisfied since
$n \geq \max\{ 4^{1/(q(1-\alpha)}, (4^{1/(1-\alpha)} \sigma^{-2/\alpha})^{\alpha/(q(1+\alpha) - 1)}\}$. Applying Theorem \ref{thm:mainGDvsRidge_Upper} with $\alpha > 1/3$ and upper bounding
\begin{align*}
    \mathcal{J}_{\alpha,\lambda_{\star}}(k,\lambda_{\min}) 
    = 
    1 + \frac{\lambda_\star^{3-1/\alpha}}{3\alpha - 1}
    \leq 
    1+ \frac{1}{3 \alpha - 1} 
    =
    \frac{3\alpha}{3\alpha - 1}
    \leq 
    \frac{3}{3\alpha - 1}
\end{align*}  
yields, using $1/\lambda_\star^{1/\alpha - 1} \ge 1$,
\begin{align*}
    \mathcal{S}\left( \mathcal{C}^{\text{GD}}(\eta,k), 
    \mathcal{C}^{\text{Ridge}}(\lambda_{\min},k)
    \right)
    & \leq
    \max\left\{1, \frac{3}{3\alpha -1} \right\}\\
    & \quad\quad\quad \times 
    \frac{2304}{\min\{1-\ep,\ep\}^2(1 - \lambda_{\min})^2}
    \left[ 
    \frac{1}{r^{1-\alpha}}
    \frac{k^2}{\lambda_\star^{1/\alpha - 1}} 
    + 
    \frac{\lambda_{\star}^4}{\lambda_{\min}^2}
    \right] .
\end{align*}

The quantity in the square brackets can be written in terms of $\sigma,\sigma_{\min},k,n$ as follows 
\begin{align*}
    \frac{1}{r^{1-\alpha}}
    \frac{k^2}{\lambda_\star^{1/\alpha - 1}} 
    + 
    \frac{\lambda_{\star}^4}{\lambda_{\min}^2} 
    & = 
    \frac{k^2}{\sigma^{2(1/\alpha - 1)}}
    \frac{1}{n^{q(1/\alpha -\alpha) - 1/\alpha + 1} }
    + 
    \frac{1}{n^{2(1-q)}} \left( \frac{\sigma^4}{\sigma^2_{\min}}\right)^2.
\end{align*}
In the case $\alpha = 1/3$ we alternatively have $\mathcal{J}_{\alpha,\lambda_{\star}}(k,\lambda_{\min}) \leq 5\log(k \lambda_{\min}/\lambda_{\star})$. Following an identical set of steps to the above yields the upper bound 
\begin{align*}
    \mathcal{S}\left( \mathcal{C}^{\text{GD}}(\eta,k), 
    \mathcal{C}^{\text{Ridge}}(\lambda_{\min},k)
    \right)
    & \leq
    5\log(k \lambda_{\min}/\lambda_{\star})
    \frac{2304}{\min\{1-\ep,\ep\}^2(1 - \lambda_{\min})^2}\\
    & \quad\quad \cdot \left[ \frac{k^2}{\sigma^{2(1/\alpha - 1)}}
    \frac{1}{n^{q(1/\alpha -\alpha) - 1/\alpha + 1} }
    + 
    \frac{1}{n^{2(1-q)}} \left( \frac{\sigma^4}{\sigma^2_{\min}}\right)^2
    \right].
\end{align*}

Combining these bounds yields the result for $\alpha \geq 1/3$. Let us now consider the case $\alpha < 1/3$, in which case we have $\mathcal{J}_{\alpha,\lambda_{\star}}(k,\lambda_{\min}) \leq 2 \left( \frac{k \lambda_{\min}}{\lambda_\star^2} \right)^{1/\alpha - 3}/(1-3\alpha)$. Following the same steps as above, by applying Theorem \ref{thm:mainGDvsRidge_Upper}, we find 
\begin{align*}
    & \mathcal{S}\left( \mathcal{C}^{\text{GD}}(\eta,k), 
    \mathcal{C}^{\text{Ridge}}(\lambda_{\min},k)
    \right)
     \leq
    \frac{1152}{\min\{1-\ep,\ep\}^2(1 - \lambda_{\min})^2}\\ 
    & \quad\quad 
    \cdot 
    \left[
    \frac{k^2}{\sigma^{2(1/\alpha - 1)}}
    \frac{1}{n^{q(1/\alpha -\alpha) - 1/\alpha + 1} }
    + \frac{2}{1-3\alpha}
    \frac{k^2 (k \sigma_{\min}^2/\sigma^4 )^{1/\alpha - 3}}{n^{q(1-\alpha) - (1-q)(1/\alpha - 3)}}
    +
    \frac{1}{n^{2(1-q)}} \left( \frac{\sigma^4}{\sigma^2_{\min}}\right)^2
    \right]\\
    & = 
    \frac{1152}{\min\{1-\ep,\ep\}^2 (1 - \lambda_{\min})^2}
    \Bigg[
    \frac{k^2}{n^{q(1/\alpha -\alpha) - 1/\alpha + 1} }
    \left(
    \frac{1}{\sigma^{2(1/\alpha - 1)}}
    + \frac{2}{1-3\alpha}
    \frac{(k \sigma_{\min}^2/\sigma^4 )^{1/\alpha - 3}}{n^{2(1-q)}}
    \right)\\
    & \qquad+
    \frac{1}{n^{2(1-q)}} \left( \frac{\sigma^4}{\sigma^2_{\min}}\right)^2
    \Bigg].
\end{align*}
Following the proof of Theorem \ref{thm:Informal}, we can use that $\sigma_{\min} \leq \sigma/\sqrt{2}$ to say $1 \geq (1-\lambda_{\min})^2 \geq 1/4$ as well as rewrite $\min\{1-\ep,\ep\} = \text{Dist}_{\delta}(\lambda_{\star},\Gamma)$.  Combining the two upper bounds thereafter leads to \eqref{ubcor}. 

For the lower bound, from Theorem \ref{thm:mainGDvsRidge_Lower}  with $r = n^{q}$ and $\lambda_\star = \sigma^2/n^{1-q}$ we find 
\begin{align*}
    & \mathcal{S}\left( \mathcal{C}^{\text{GD}}(\eta,k), 
    \mathcal{C}^{\text{Ridge}}(\lambda_{\min},k)
    \right) \\
    & \geq 
    \frac{1}{512} 
    \frac{\min\{ 1- \kappa,\kappa \}^2}{\min\{ 1-\ep, \ep \}^2 (1-\lambda_{\min})^2}
    \left( 
    1 - \frac{2}{1-\alpha} \left( \frac{4 k \lambda_{\min}}{\lambda_{\star}^2 \kappa }\right)^{1/\alpha - 1}
    \frac{1}{r^{1-\alpha}}
    \right)
    \frac{\lambda_{\star}^4}{\lambda_{\min}^2}\\
    & = 
    \frac{1}{512} 
    \frac{\min\{ 1- \kappa,\kappa \}^2}{\min\{ 1-\ep, \ep \}^2 (1-\lambda_{\min})^2}
    \left( 
    1 - \frac{2}{1-\alpha} \frac{(4 k \sigma_{\min}^2/(\sigma^4 \kappa) )^{1/\alpha - 1}}{n^{q(1/\alpha - \alpha) - 1/\alpha +1}}
    \right)
    \frac{(\sigma^4/\sigma_{\min}^2)^2}{n^{2(1-q)}}.
\end{align*}
To obtain the required bound, we can simply follow the previous steps to bound $1 \geq (1-\lambda_{\min})^2 \geq 1/4$ due to  $\sigma_{\min} \leq \sigma/\sqrt{2}$, and write $\min\{1-\ep,\ep\} = \text{Dist}_{\delta}(\lambda_{\star},\Gamma)$.

\subsection{Proofs for Individual Classes of Ridge Regression Estimators}
\label{sec:proof:RidgeBounds}
In this section we present upper and lower bounds on the estimation error of model classes produced by ridge regression.  Section \ref{sec:proof:RidgeBounds:Fine} presents the proof of Proposition \ref{lem:RidgeBounds}, which provides upper and lower bounds when the discretization length is \emph{smaller} than the optimal amount of regularization.  Section \ref{sec:proof:Ridgebounds:Coarse} presents the proof of Proposition \ref{lem:RidgeBounds:Coarse}, which gives upper and lower bounds when the discretization length is \emph{larger} than the optimal amount of regularization. 

\subsubsection{Proof of Proposition \ref{lem:RidgeBounds}}
\label{sec:proof:RidgeBounds:Fine}
Writing $\lambda = \lambda_{\star} + \Delta$ for some $\Delta  > - \lambda_{\star} $ we find  
\begin{align*}
    \E_{\beta_{\star},\epsilon}[L_{\beta_{\star}}(\widehat{\beta}_{\lambda})] 
    -
    \lambda_{\star} \int \frac{1}{s + \lambda_{\star}}
    d \widehat{H}(s)
    & = 
    \frac{\sigma^2}{n} \sum_{i=1}^{r} \left( \frac{1}{s_i} + \frac{1}{\lambda_{\star}} \right) 
    \left( \frac{\lambda}{\lambda + s_i} - \frac{\lambda_{\star}}{\lambda_{\star} + s_i} \right)^2 \\
    & = 
    \frac{\sigma^2}{n\lambda_{\star}} \sum_{i=1}^{r}
    \frac{s_i \Delta^2}{(\lambda_{\star} + \Delta + s_i)^2(\lambda_{\star} + s_i)}.
\end{align*}
The above function is increasing in $\Delta$ for $\Delta > 0$ and decreasing in $\Delta$ for $-\lambda_\star< \Delta < 0$. 
Recall that $\lambda_{\star} = \lambda_{\min} + j \delta + \delta'$ for an integer $j \in \{0,1,\dots,k-2\}$, where $\frac{1-\lambda_{\min}}{\delta}=k-1$ and $\delta \geq \delta' \geq 0$. Therefore,   
\begin{align*}
    & 
    \mathcal{E}(\mathcal{C}^{\text{Ridge}}(\lambda_{\min},k)) 
     = 
    \min_{\lambda \in \{\lambda_{\min}, \lambda_{\min} + \delta, \lambda_{\min} + 2\delta,\dots,1\} }
    \left(
    \E_{\beta_{\star},\epsilon}[L_{\beta_{\star}}(\widehat{\beta}_{\lambda})] 
    -
    \lambda_{\star} \int \frac{1}{s + \lambda_{\star}}
    d \widehat{H}(s)
    \right)\\
    & = 
    \min_{\Delta \in \{- j \delta - \delta', -(j-1)\delta - \delta',\dots,- \delta',\delta - \delta', 2\delta- \delta',\dots, (k-1-j)\delta - \delta'\}}
    \frac{\sigma^2}{n\lambda_{\star}} \sum_{i=1}^{r}
    \frac{s_i \Delta^2}{(\lambda_{\star} + \Delta + s_i)^2(\lambda_{\star} + s_i)}\\
    & = 
    \min_{\Delta \in \{- \delta', \delta- \delta'\}}
    \frac{\sigma^2}{n\lambda_{\star}} \sum_{i=1}^{r}
    \frac{s_i \Delta^2}{(\lambda_{\star} + \Delta + s_i)^2(\lambda_{\star} + s_i)},
\end{align*}
where for the second equality we rewrote the minimum in terms of the difference $\lambda - \lambda_{\star}$, 
and for the third we used the monotonicity in $\Delta$ noted previously. Since $\delta' =\ep \delta$ for $\ep \in [0,1]$, we find
\begin{align*}
    & \mathcal{E}(\mathcal{C}^{\text{Ridge}}(\lambda_{\min},k)) 
    = 
    \min_{\kappa \in \{- \ep,(1- \ep) \} }
    \frac{\sigma^2}{n\lambda_{\star}} \sum_{i=1}^{r}
    \frac{s_i \kappa^2 \delta^2}{(\lambda_{\star} + \kappa \delta + s_i)^2(\lambda_{\star} + s_i)}.
\end{align*}
Since by assumption $0< \delta \leq \lambda_{\star}$ we then find for any $s\ge 0$
\begin{align*}
    (1 - \ep) \frac{\delta}{\lambda_{\star} + s} 
    \geq \frac{(1- \ep)\delta}{\lambda_{\star} + (1- \ep) \delta + s} \geq 
    \frac{1- \ep}{2} \frac{\delta }{\lambda_{\star} + s}\\
    \frac{\ep}{1- \ep}
    \frac{\delta}{\lambda_{\star} +  \delta}
    \geq 
    \Big| \frac{- \ep \delta}{\lambda_{\star}  - \ep \delta + s}
    \Big|
    = \frac{\ep \delta}{\lambda_{\star} - \ep \delta + s}
    \geq 
    \ep \frac{\delta}{\lambda_{\star} + s}.
\end{align*}
This provides the upper and lower bound 
\begin{align*}
    & \mathcal{E}(\mathcal{C}^{\text{Ridge}}(\lambda_{\min},k))  
    \geq \frac{\min\{1-\ep,\ep\}^2 }{4}
    \frac{\sigma^2}{n\lambda_{\star}}
    \sum_{i=1}^{r} 
    \frac{s_i \delta^2}{(\lambda_{\star} + s_i)^3}\\
    & \mathcal{E}(\mathcal{C}^{\text{Ridge}}(\lambda_{\min},k)) 
    \leq 
    \min\left\{1- \ep,\frac{\ep}{1- \ep} \right\}^2
    \frac{\sigma^2}{n\lambda_{\star}}
    \sum_{i=1}^{r} 
    \frac{s_i \delta^2}{(\lambda_{\star} + s_i)^3}.
\end{align*}
The final bound is arrived at by noting that for $s > 0$ we have
\begin{align*}
    \frac{1}{8} G^{\text{Ridge}}(s) 
    \leq 
    \frac{s\delta^2}{\lambda_{\star}(\lambda_{\star} + s)^3}
    \leq 
     G^{\text{Ridge}}(s)
\end{align*}
and $\min\{1- \ep,\frac{\ep}{1- \ep}\} \leq \min\{ 1- \ep , 2\ep\} \leq 2 \min\{1- \ep,\ep\}$.

\subsubsection{Proof of Proposition \ref{lem:RidgeBounds:Coarse}}
\label{sec:proof:Ridgebounds:Coarse}
Following the proof of Proposition \ref{lem:RidgeBounds}, write $\lambda = \lambda_{\star} + \Delta$ for some $1 - \lambda_{\star} \geq \Delta \geq 0 $. We then have with $\lambda_{\min} = 0$, $\lambda_{\star} = \delta - \delta'$, and recalling that $\delta' < \delta$
\begin{align*}
    \mathcal{E}(\mathcal{C}^{\text{Ridge}}(0,k)) 
     & = 
    \min_{\lambda \in \{0,  \delta, 2\delta,\dots,1\} }
    \left(
    \E_{\beta_{\star},\epsilon}[L_{\beta_{\star}}(\widehat{\beta}_{\lambda})] 
    -
    \lambda_\star \int \frac{1}{s + \lambda_{\star}}
    d \widehat{H}(s)
    \right)\\
    & =
    \min_{\Delta \in \{ - \delta + \delta', \delta',\delta + \delta',, \dots,1 - \delta + \delta'\}}
    \frac{\sigma^2}{n\lambda_{\star}} \sum_{i=1}^{r}  
    \frac{s_i \Delta^2}{(\lambda_{\star} + \Delta + s_i)^2(\lambda_{\star} + s_i)}\\
    & = 
    \min_{\Delta \in \{ - \delta + \delta', \delta'\}}
    \frac{\sigma^2}{n\lambda_{\star}} \sum_{i=1}^{r}  
    \frac{s_i \Delta^2}{(\lambda_{\star} + \Delta + s_i)^2(\lambda_{\star} + s_i)},
\end{align*}
where we recall that the function is monotonic for $\Delta > 0$ and $\Delta < 0$. Picking $\Delta = \delta'$ we find the upper bound
\begin{align*}
    & \min_{\Delta \in \{ - \delta + \delta', \delta'\}}
    \frac{\sigma^2}{n\lambda_{\star}} \sum_{i=1}^{r}  
    \frac{s_i \Delta^2}{(\lambda_{\star} + \Delta + s_i)^2(\lambda_{\star} + s_i)}
    \leq  \delta'^2
    \frac{\sigma^2}{n\lambda_{\star}} \sum_{i=1}^{r}  
    \frac{s_i^2 }{(\delta + s_i)^2(\lambda_{\star} + s_i)^2 }.
\end{align*}
Meanwhile for the lower bound, taking the minimum into the sum
\begin{align*}
    & \min_{\Delta \in \{ - \delta + \delta', \delta'\}}
    \frac{\sigma^2}{n\lambda_{\star}} \sum_{i=1}^{r}  
    \frac{s_i \Delta^2}{(\lambda_{\star} + \Delta + s_i)^2(\lambda_{\star} + s_i)}\\
    &\geq 
    \frac{\sigma^2}{n\lambda_{\star}} \sum_{i=1}^{r}  
    \min_{\Delta \in \{ - \delta + \delta', \delta'\}}
    \frac{s_i \Delta^2}{(\lambda_{\star} + \Delta + s_i)^2(\lambda_{\star} + s_i)}\\
    & =
    \delta'^2 \frac{\sigma^2}{n\lambda_{\star}} \sum_{i=1}^{r}  
    \frac{s_i  }{(\delta + s_i)^2(\lambda_{\star} + s_i) }.
\end{align*}
We used that for $s > 0$
\begin{align*}
    \frac{\delta'}{\delta + s}
    \leq 
    \frac{\delta - \delta'}{s} 
    & \iff 
    \delta' s \leq \delta^2 - \delta' \delta - \delta' s + s \delta \\
    & \iff 
    s(2 \delta' - \delta) \leq \delta(\delta- \delta'),
\end{align*}
which holds since $\delta' \leq \delta/2$. From \eqref{cri} we have
\begin{align*}
    \frac{1}{4} G^{\text{Coarse Ridge}}(s) \leq
    \frac{s \delta^2 }{\lambda_{\star}(\delta + s)^2(\lambda_{\star} + s)}
    \leq 
    G^{\text{Coarse Ridge}}(s),
\end{align*}
which yields the result.

\subsection{Proofs for Individual Classes of Gradient Descent Estimators}
\label{sec:Proof:GradientDescentEstimators}
In this section we present the proofs for bounds on the estimation error of gradient descent. Before providing the proof of these results, we need to define the following function, for $s\in (0,1/\eta)$: 
\begin{align}\label{tstar}
    t^{\star}(s) := \frac{\log(1 + s/\lambda_{\star})}{-\log(1-\eta s)}.
\end{align}
By continuity, we may also define $t^{\star}(0) := \lim_{s\searrow 0} t^{\star}(s)  = 1/(\eta\lambda_\star)$.
The follow proposition presents a bound on $t^{\star}(s)$ when eigenvalues are below the level of regularization, i.e.,  $ s \in [0,\lambda_{\star}]$. We note that $t^{\star}(s)$ is well defined on this interval since we assume $\eta \le 1/\max\{s_1,\lambda_{\star}\}$. If $s_1=\lambda_\star$ and $\eta = 1/s_1$, by continuity we define $t^\star(\lambda_\star)=0$.
\begin{proposition}[Discretization Limit]
\label{prop:IterationsLimit}
For $u \in [0,\lambda_{\star}]$, we have
\begin{align*}
    \frac{1}{\eta \lambda_{\star}} - t^{\star}(u) 
    \leq 
    \frac{u}{\lambda_{\star}}\left(1+\frac{1}{2 \eta \lambda_{\star}} \right).
\end{align*}
\end{proposition}

Given this result, we can now present the proofs for bounding the estimation error of models produced by gradient descent. Section \ref{sec:proof:GDBounds:Fine} presents the proof of Proposition \ref{lem:GDBounds:Fine} which considers the fine discretization case.  Section \ref{sec:proof:GDBounds:Coarse} presents the proof Proposition \ref{lem:GDBounds:Coarse}, which considers the coarse discretization case. Section \ref{sec:proof:technicalProp} then presents the proof of Proposition \ref{prop:IterationsLimit}. 

\subsubsection{Proof of Proposition \ref{lem:GDBounds:Fine}}
\label{sec:proof:GDBounds:Fine}
    Let us begin by proving the upper bound. Since the sub-optimality gap evaluates the minimum over all models up to iteration $t$, it is enough to consider $t = \lceil \frac{1}{\eta \lambda_{\star}} \rceil$. Following the proof of Lemma \ref{lem:ErrorDecomp} we find
    \begin{align}\label{gdid}
        \E_{\beta_{\star},\epsilon}[L_{\beta_{\star}}(\widehat{\beta}_{\eta,t})]
    -
    \lambda_\star \int \frac{1}{s + \lambda_{\star}}
    d \widehat{H}(s)
    = 
    \frac{\sigma^2}{n}
    \sum_{i=1}^r 
    \left( \frac{1}{s_i} + \frac{1}{\lambda_{\star}} \right)
    \left( (1-\eta s_i)^{t} - \frac{\lambda_{\star}}{\lambda_{\star} + s_i}
    \right)^2.
    \end{align}
    We now set to bound $( (1-\eta s)^{t} - \frac{\lambda_{\star}}{\lambda_{\star} + s})^2$ splitting into two cases: $s \in (\lambda_{\star},1/\eta]$ and $ s \in (0,\lambda_{\star}]$. For the upper bound we will use both cases, whereas for the lower bound we will only consider $s \in (0, \lambda_\star]$. 
    
    Let us begin with the upper bound when $s \in (\lambda_{\star},1/\eta)$. One can check that $s \rightarrow t^{\star}(s)$ from \eqref{tstar} is a decreasing  function of $s$ over the entire interval $s\in[0,1/\eta)$ with the maximum achieved at the limit $\lim_{s \rightarrow 0} t^{\star}(s) = \frac{1}{\eta \lambda_{\star}}$. Noting that $\frac{\lambda_{\star}}{\lambda_{\star} + s} = (1-\eta s)^{t^{\star}(s)}$ as well as that $t \geq \frac{1}{\eta \lambda_{\star}} \geq t^{\star}(s)$ we immediately find for $s \in (\lambda_{\star},1/\eta)$
    \begin{align*}
        0 & \leq \frac{\lambda_{\star}}{\lambda_{\star} + s} - (1-\eta s)^{t}  \leq \frac{\lambda_{\star}}{\lambda_{\star} + s}. 
    \end{align*}
    When $s =1/\eta$ we get $\frac{\lambda_{\star}}{\lambda_{\star} + s} - (1-\eta s)^{t} = \frac{\lambda_{\star}}{\lambda_{\star} + s}$. 
    Using this upper bound the terms within the sum in \eqref{gdid} can be bounded when $s \in (\lambda_{\star},1/\eta]$ as
    \begin{align}
    \label{b1}
    \left( \frac{1}{s} + \frac{1}{\lambda_{\star}} \right)
    \left( (1-\eta s)^{t} - \frac{\lambda_{\star}}{\lambda_{\star} + s}
    \right)^2
    \leq 
    \left( \frac{1}{s} + \frac{1}{\lambda_{\star}} \right)
    \left(\frac{\lambda_\star}{\lambda_\star + s }\right)^2
    = 
    \frac{1}{s}\frac{\lambda_\star}{\lambda_\star + s }
    \leq \frac{\lambda_\star}{s^2}.
    \end{align}
     Let us now consider $ s \in (0,\lambda_{\star}]$. From   Proposition \ref{prop:IterationsLimit} we have 
    $t- t^{\star}(s) = \lceil \frac{1}{\eta \lambda_{\star}} \rceil - t^{\star}(s) \leq 1 + \frac{1}{\eta \lambda_{\star}} - t^{\star}(s) \leq 1 + \frac{s}{2\eta \lambda_{\star}^2} + \frac{s}{\lambda_{\star}} \leq 1 + \frac{3}{2} \frac{s}{\eta \lambda_{\star}^2}$  
    where we have $\eta \lambda_\star^2 \leq 1$ from the assumption on the step size $\eta \le 1/\max\{s_1,\lambda_{\star}\}$ and $\lambda_{\star} \le 1$.
     
    Therefore, using that $(1-\eta s)^{t - t^{\star}(s)}\geq (1- \eta s)^{1 + 2 \frac{s}{\eta \lambda_{\star}^2} }$ since $0 \leq t - t^{\star}(s) \leq 1 + 3/2\cdot \frac{s}{\eta \lambda_{\star}^2}$ and $1-\eta s \leq  1$, and from Bernoulli's inequality $1 - (1-x)^{\alpha} \leq \alpha x$ for $\alpha \ge1 $ and $x\in [0,1]$, we find
    \begin{align*}
        0 & \leq \frac{\lambda_{\star}}{\lambda_{\star} + s} - (1-\eta s)^{t}\nonumber
         = 
        (1-\eta s)^{t^{\star}(s)} - (1-\eta s)^{t} \nonumber
         = (1-\eta s)^{t^{\star}(s)}
        \left( 1 - (1-\eta s)^{t - t^{\star}(s)} 
        \right) \nonumber\\
        & \leq 
        1 - (1-\eta s)^{t - t^{\star}(s)}
        \leq 
        1 - (1- \eta s)^{1 + \frac{3}{2} \frac{s}{\eta \lambda_{\star}^2} } \nonumber\\
        & \leq \eta s 
        \left(1 + \frac{3}{2} \frac{s}{\eta \lambda_{\star}^2}
        \right)
         \leq 
        3
        \max\left\{ \eta s, \left(\frac{s}{\lambda_{\star}}\right)^2 \right\}.
    \end{align*}
  Thus, for $ s \in (0,\lambda_{\star}]$ 
    \begin{align}
    \label{b2}
        \left( \frac{1}{s} + \frac{1}{\lambda_{\star}} \right)
    \left( (1-\eta s)^{t} - \frac{\lambda_{\star}}{\lambda_{\star} + s}
    \right)^2
        \leq \frac{18}{s} 
        \max\left\{ \eta s, \left(\frac{s}{\lambda_{\star}}\right)^2 \right\}^2.
    \end{align}
    
     Combining \eqref{b1}, \eqref{b2} we then find for $s \in (0,1/\eta]$
     \begin{align*}
         \left( \frac{1}{s} + \frac{1}{\lambda_{\star}} \right)
    \left( (1-\eta s)^{t} - \frac{\lambda_{\star}}{\lambda_{\star} + s}
    \right)^2
    \leq  
    18 G^{\text{GD}}(s),
     \end{align*}
     as required. The same bound holds for any $t \geq \lceil \frac{1}{\eta \lambda_{\star}} \rceil$ since $\frac{\lambda_{\star}}{\lambda_{\star} + s} - (1-\eta s)^{t} \geq \frac{\lambda_{\star}}{\lambda_{\star} + s} - (1-\eta s)^{\lceil \frac{1}{\eta \lambda_{\star}} \rceil} \ge 0$. 
     
    We now prove the lower bound. Recall  $\frac{1}{\eta \lambda_{\star}} = k + \kappa$ for some $\kappa \in (0,1)$. We first consider $t^{\text{Upper}} = \lceil \frac{1}{\eta \lambda_{\star}} \rceil = \frac{1}{\eta \lambda_{\star}} + 1-\kappa$ and $t^{\text{Lower}} = \lfloor \frac{1}{\eta \lambda_{\star}} \rfloor  = \frac{1}{\eta \lambda_{\star}} - \kappa $ in turn. The bound for other values of  $t$ follows as a direct consequence, due to the monotonicity of the objective we consider. 
    
    \textbf{Case 1} Suppose $t = t^{\text{Upper}}$. Then $t \geq t^{\star}(s)$ for any $s \in (0,\lambda_{\star}]$, so  we have
    \begin{align*} 
        \frac{\lambda_{\star}}{\lambda_{\star} + s} -  (1-\eta s)^{t^{\text{Upper}}} & = 
        (1-\eta s)^{t^{\star}(s)} - (1-\eta s)^{t^{\text{Upper}}}\\
        & = 
        (1-\eta s)^{t^{\star}(s)}
        \left( 1 - (1-\eta s)^{t^{\text{Upper}} - t^{\star}(s)}\right)\\
        & = 
        \left( \frac{\lambda_{\star}}{\lambda_{\star} + s} \right) 
        \left( 1 - (1- \eta s)^{ 1- \kappa + \frac{1}{\eta \lambda_{\star}} - t^{\star}(s)}
        \right)\\
        & = 
        \left( \frac{\lambda_{\star}}{\lambda_{\star} + s} \right) 
        \left( 1 - (1 - \eta s)^{1-\kappa} (1- \eta s)^{\frac{1}{\eta \lambda_{\star}} - t^{\star}(s)}
        \right)\\
        & \geq 
        \frac{1}{2} 
        \left( 1 - (1-\eta s)^{1-\kappa} \right)
         \geq 
        \frac{1}{2} (1-\kappa) \eta s, 
    \end{align*}
where we have used that $\frac{1}{\eta \lambda_{\star}} \geq t^{\star}(s)$ as well as that $s \leq \lambda_{\star}$ and therefore $\frac{\lambda_{\star}}{\lambda_{\star} + s} \geq \frac{1}{2}$. In the final inequality we used that $1 - (1-x)^{\alpha} \geq \alpha x$ for $\alpha \in [0,1]$. Plugging in this bound we find for $ s \in (0,\lambda_{\star}]$
\begin{align*}
    \left( \frac{1}{s} + \frac{1}{\lambda_{\star}} \right) 
    \left( (1-\eta s)^{t^{\text{Upper}}} - \frac{\lambda_{\star}}{\lambda_{\star} + s}\right)^2 
    & \geq 
    \frac{1}{4} (1-\kappa)^2 \eta^2 s, 
\end{align*}
as required. 
The same bound holds for any $t \geq t^{\text{Upper}}$ since $\frac{\lambda_{\star}}{\lambda_{\star} + s} - (1-\eta s)^{t} \geq \frac{\lambda_{\star}}{\lambda_{\star} + s} - (1-\eta s)^{t^{\text{Upper}}}$.  The lower bound on the expression in \eqref{gdid} for $t \geq t^{\text{Upper}}$ is then arrived at by considering the eigenvalues less than or equal to $\lambda_{\star}$, specifically,  
\begin{align*}
    \eqref{gdid}
    & \geq 
    \frac{(1-\kappa)^2 \sigma^2}{4n} \sum_{i : 0 < s_i \leq \lambda_\star} 
    \eta^2 s_i \\
    & = 
    \frac{ (1-\kappa)^2 \lambda_{\star}}{4d}\sum_{i : 0 < s_i \leq \lambda_\star} 
    \eta^2 s_i 
     \geq \frac{1}{4} (1-\kappa)^2 
    \lambda_{\star}
    \int_0^{\kappa \eta \lambda_\star^2/4} G^{\text{GD}}(s) d\widehat{H}(s).
\end{align*}

\textbf{Case 2} Let us suppose that $t = t^{\text{Lower}} = \frac{1}{\eta \lambda_{\star}} - \kappa$. From  Proposition \ref{prop:IterationsLimit}, we find $t^{\star}(s)  \geq \frac{1}{\eta \lambda_{\star}} - \frac{s}{\eta\lambda_{\star}^2} - \frac{s}{\lambda_{\star}} \geq \frac{1}{\eta \lambda_{\star}} - 2 \frac{s}{\eta \lambda_{\star}^2}$,  where we used once again used that $\eta \lambda_\star^2 \leq 1$.
Therefore, if $s \in (0,\frac{\kappa}{4} \eta\lambda_{\star}^2]$, we have $t^{\star}(s) \geq \frac{1}{\eta \lambda_{\star}} - \frac{\kappa}{2} \geq \frac{1}{\eta \lambda_{\star}} - \kappa = t^{\text{Lower}}$. As such for $s \in (0,\frac{\kappa}{4} \eta\lambda_{\star}^2]$ we have
\begin{align*}
    (1- \eta s)^{t^{\text{Lower}}} - (1-\eta s)^{t^{\star}(s)}
    & = 
    (1- \eta s )^{t^{\text{Lower}}} (1 - (1-\eta s)^{t^{\star}(s) - t^{\text{Lower}}})\\
    & = 
    (1- \eta s )^{t^{\text{Lower}}} 
    (1 - (1-\eta s)^{t^{\star}(s) - \frac{1}{\eta \lambda_{\star}} + \kappa })\\
    & \geq 
    (1- \eta s )^{t^{\text{Lower}}} 
    (1 - (1-\eta s)^{ \kappa  - \frac{\kappa}{2} })
 \geq \frac{1}{4} \eta s \kappa,
\end{align*}
where in the final inequality we used that $(1-\eta s)^{t^{\text{Lower}}} \geq (1-\eta s)^{t^{\star}(s)} = \frac{\lambda_{\star}}{\lambda_{\star} + s} \geq \frac{1}{2}$ (as $s \leq \lambda_{\star}$) as well as that $1 - (1-x)^{\alpha} \geq \alpha x$ for $\alpha \in [0,1]$. Thus we find for $s \in (0,\frac{\kappa}{4} \eta\lambda_{\star}^2]$
\begin{align*}
    \left( \frac{1}{s} + \frac{1}{\lambda_{\star}} \right) 
    \left( (1-\eta s)^{t^{\text{Lower}}} - \frac{\lambda_{\star}}{\lambda_{\star} + s}\right)^2
    & \geq 
    \frac{1}{16} \eta^2 s \kappa^2. 
\end{align*}
All other cases  $t \leq t^{\text{Lower}}$ then arise from the above inequality since  $(1-\eta s )^{t} \geq(1-\eta s)^{t^{\text{Lower}}}$ in this case. Similar to the case of $t^{\text{Upper}}$, the lower bound on \eqref{gdid} for $t \leq t^{\text{Lower}}$ is arrived at by considering the eigenvalues in the interval $(0,\frac{\kappa}{4} \eta \lambda_{\star}^2]$, precisely, 
\begin{align*}
    \eqref{gdid}
    \geq \frac{\kappa^2 \lambda_\star }{16} 
    \frac{1}{d} \sum_{i : 0 < s_i \leq \kappa \eta \lambda_\star^2/4 } \eta^2 s_i = 
    \frac{\kappa^2 \lambda_\star }{16}  \int_0^{\kappa \eta \lambda_\star^2/4 } 
    G^{\text{GD}}(s) d \widehat{H}(s)
\end{align*}
Taking a minimum over the lower bounds for $t=t^{\text{Upper}},t^{\text{Lower}}$ in \eqref{gdid} then yields the result.

\subsubsection{Proof of Proposition \ref{lem:GDBounds:Coarse}}
\label{sec:proof:GDBounds:Coarse}
For $u \in (\lambda^\star, 1/\eta)$, we parametrize the number of iterations as $t = \frac{1}{\eta \sqrt{\lambda_{\star}u}}$. We can recover the stated result by setting $u = 1/((t \eta)^{2} \lambda_{\star})$. We present a bound for the eigenvalues in each of the sets $(u,1/\eta]$, $(\sqrt{\lambda_{\star} u},u]$, $(\lambda_{\star},\sqrt{\lambda_{\star}u}]$ and $(0,\lambda_{\star}]$. 

\textbf{Case 1} Consider $s \in (u,1/\eta]$. 
Using the inequality $\log(1+x) \leq \frac{x}{\sqrt{1+x}}$ as well as that $\frac{1}{-\log(1-\eta u)} \leq  \frac{1}{\eta u}$ we have $t^{\star}(u) \leq \frac{ u/\lambda_{\star}}{\eta u \sqrt{1 + u/\lambda_{\star}}} = \frac{1}{\eta \sqrt{\lambda_{\star}}\sqrt{\lambda_{\star} + u}} \leq \frac{1}{\eta \sqrt{\lambda_{\star} u}} =t$. Using that $t^{\star}(\cdot)$ is decreasing, we have $t =  \frac{1}{\eta \sqrt{\lambda_{\star} u }}  \geq  t^{\star}(u) \geq t^{\star}(s)  $ for $s \in (u,1/\eta)$. Therefore for $s \in (u,1/\eta)$ we have 
\begin{align*}
    0 \leq (1 - \eta s)^{t^{\star}(s)} - (1-\eta s)^{t }  = \frac{\lambda_\star}{\lambda_{\star} + s} - (1-\eta s)^{t} \leq \frac{\lambda_{\star}}{\lambda_{\star} + s}.
\end{align*}
For the case $s = 1/\eta$ we have $\frac{\lambda_\star}{\lambda_{\star} + s} - (1-\eta s)^{t} = \frac{\lambda_\star}{\lambda_{\star} + s}$. Bringing together these cases we find, for $s \in (u,1/\eta]$,
\begin{align*}
    \left( \frac{1}{s} + \frac{1}{\lambda_{\star}}\right) \left( (1-\eta s)^{t} - \frac{\lambda_{\star}}{\lambda_{\star} + s}\right)^2
    \leq 
    \left( \frac{1}{s}  + \frac{1}{\lambda_{\star}} 
    \right) \left( \frac{\lambda_{\star}}{\lambda_{\star} + s}\right)^2 
    \leq  \frac{\lambda_{\star}}{s^2}.
\end{align*}

\textbf{Case 2} Consider $s \in (\sqrt{\lambda_{\star} u},u]$. Note that for
\begin{align*}
    (1- \eta s)^{k} \leq \frac{\lambda_{\star} u}{\lambda_{\star} u + s^2} 
\end{align*}
it is sufficient to have  $k \geq \frac{\log(1 + s^2/(\lambda_{\star}u))}{-\log(1-\eta s)}$. This is satisfied with $k=t$ since
\begin{align*}
    t= \frac{1}{\eta \sqrt{\lambda_{\star} u }} \geq  \frac{1}{\eta s} \frac{s^2/(\lambda_{\star} u )}{\sqrt{1 + s^2/(\lambda_{\star} u)}} \geq  \frac{\log(1 + s^2/(\lambda_{\star}u))}{-\log(1-\eta s)},
\end{align*}
where we used the inequalities $\log(1+x) \leq \frac{x}{\sqrt{1+x}}$ and $\log(1/(1-x)) \geq x$. Therefore if $t < t^{\star}(s)$ we  find 
\begin{align*}
    0 \leq (1- \eta s )^{t} - (1-\eta s )^{t^{\star}(s)}
    \leq 
    \frac{\lambda_{\star} u}{\lambda_{\star} u + s^2}.
\end{align*}
Meanwhile if $t \geq t^{\star}(s)$ we have  
\begin{align*}
    0 \leq (1-\eta s)^{t^{\star}(s)} - (1-\eta s)^{t}
    \leq (1 - \eta s)^{t^{\star}(s)}  = \frac{\lambda_{\star}}{\lambda_{\star} + s}
    \leq \frac{\lambda_{\star}u}{\lambda_{\star} u + s^2},
\end{align*}
where the last inequality arises from $x \rightarrow x/(1+x)$ being increasing, and by using $s\le u$. Therefore we have for $s \in (\sqrt{\lambda_{\star}u},u]$ (so that $s \geq \lambda_{\star}$)
\begin{align*}
    \left( \frac{1}{s} + \frac{1}{\lambda_{\star}}\right) \left( (1-\eta s)^{t} - \frac{\lambda_{\star}}{\lambda_{\star} + s}\right)^2
    \leq
    \frac{2}{\lambda_{\star}} \left( \frac{\lambda_{\star} u}{\lambda_{\star} u + s^2} 
    \right)^2
    \leq 2 \frac{u^2 \lambda_{\star}}{s^4}.
\end{align*}

\textbf{Case 3} 
Let us suppose that $s \in (\lambda_{\star},\sqrt{\lambda_{\star}u}]$. In this case, if $t < t^{\star}(s)$ we find 
\begin{align*}
    0 \leq 
    (1 - \eta s)^{t} - (1-\eta s)^{t^{\star}(s)}
    \leq 
    1 - \frac{\lambda_{\star}}{\lambda_{\star}  + s}
    = \frac{s}{\lambda_{\star} +s}.
\end{align*}
On the other hand, if $t^{\star}(s) \le t$  we have 
\begin{align*}
    0 < 
    (1- \eta s)^{t^{\star}(s)}
    - 
    (1 - \eta s)^{t} 
    \leq 
    \frac{\lambda_{\star}}{\lambda_{\star} + s} 
    \leq 
    \frac{s}{\lambda_{\star} + s},
\end{align*}
where we note that $s > \lambda_{\star}$ and therefore $\frac{\lambda_{\star}}{\lambda_{\star} + s} \leq \frac{s}{\lambda_{\star} + s}$. Plugging in we find for $s \in (\sqrt{\lambda_{\star}u},u]$
\begin{align*}
    \left( \frac{1}{s} + \frac{1}{\lambda_{\star}}\right) \left( (1-\eta s)^{t} - \frac{\lambda_{\star}}{\lambda_{\star} + s}\right)^2
    \leq 
    \frac{1}{\lambda_{\star}}.
\end{align*}

\textbf{Case 4}
Let us suppose that $s \in (0,\lambda_{\star}]$. Now, note that 
\begin{align*}
    t^{\star}(s) = \frac{\log(1 + s/\lambda_{\star})}{-\log(1-\eta s)}
    \geq
    \frac{1-\eta s}{\eta s} \frac{s/\lambda_{\star}}{1+s/\lambda_{\star}}
    = 
    \frac{1-\eta s}{\eta} \frac{1}{s + \lambda_{\star}}
    \geq 
    \frac{1-\eta s}{\eta} \frac{1}{2 \lambda_{\star}}
    \geq 
    \frac{1 - \eta \lambda_{\star}}{2 \eta \lambda_{\star}} 
    \geq t.
\end{align*}
Therefore
\begin{align*}
    0 \leq (1 - \eta s)^{t} - (1-\eta s)^{t^{\star}(s)}
    \leq 
    1 - \frac{\lambda_{\star}}{\lambda_{\star} + s}
    = 
    \frac{s}{\lambda_{\star} + s},
\end{align*}
which when plugging in yields 
\begin{align*}
    \left( \frac{1}{s} + \frac{1}{\lambda_{\star}}\right) \left( (1-\eta s)^{t} - \frac{\lambda_{\star}}{\lambda_{\star} + s}\right)^2
    \leq 
    \frac{s}{\lambda_{\star}(\lambda_{\star} + s)}\
    \leq 
    \frac{s}{\lambda_{\star}^2}. 
\end{align*}
Combining the bounds for each of the cases, plugging in $u = (\eta t)^2/\lambda_{\star}$, and recalling \eqref{coarsegd}  then yields the required result
\begin{align*}
    \left( \frac{1}{s} + \frac{1}{\lambda_{\star}}\right) \left( (1-\eta s)^{t} - \frac{\lambda_{\star}}{\lambda_{\star} + s}\right)^2
    \leq 
    2 G^{\text{Coarse GD}}(s). 
\end{align*}

\subsubsection{Proof of Proposition \ref{prop:IterationsLimit}}
\label{sec:proof:technicalProp}
We only need to consider the case when $s \eta < 1$. We proceed to prove the following more general inequality for $u,s$ such that $1/\eta > u\geq s > 0$: 
\begin{align}
\label{equ:GeneralDifference}
      t^{\star}(s) - t^{\star}(u) \leq 
       \frac{ \sqrt{1 - \eta u} }{2 \eta \lambda_{\star}} 
    \int_{s}^{u} \frac{1}{z + \lambda_{\star}} dz
    + 
    \log\left(1 + \frac{s}{\lambda_{\star}}\right) 
    \frac{u-s + \eta s^2/12}{s \sqrt{1- \eta s}}.
\end{align}
Proposition \ref{prop:IterationsLimit} then arises from taking $s \rightarrow 0$ in the above bound. Indeed, we have $\lim_{s\rightarrow 0}t^{\star}(s) = \frac{1}{\eta \lambda_{\star}} $. We can bound the integral $\int_{s}^{u} \frac{1}{z + \lambda_{\star}} dz \leq \frac{u}{\lambda_{\star}}$ and note that $\lim_{s \rightarrow 0}\log(1 +s/\lambda_{\star})/s = 1/\lambda_{\star}$, leading to Proposition \ref{prop:IterationsLimit}.

To prove \eqref{equ:GeneralDifference}, we use the inequality \cite{topsok2006some}\footnote{
The lower bound follows from $\log(1+x)/x \geq \frac{3(2+x)}{6 + 6x + x^2}$ \cite{topsok2006some}, since we have we have that $\frac{3(2+x)}{\sqrt{6}\sqrt{6 + 6x + x^2}} \geq 1$. 
}
\begin{align}
    \label{equ:LogInqSharp} 
    \frac{1}{\sqrt{1 + x + x^2/6}}
    \leq 
    \frac{\log(1+x)}{x} 
    \leq 
    \frac{1}{\sqrt{1+x}}
    \quad \text{ for } x > 0.
\end{align}
Using this upper bound for $-\log(1-\eta u) = \log( \frac{\eta u}{1-\eta u} + 1)$ and the lower bound for $-\log(1-\eta s) = \log\left( \frac{\eta s}{1-\eta s} + 1 \right)$, we find   
\begin{align}
    & t^{\star}(s) - t^{\star}(u) 
    = 
    \frac{\log(1 + s/\lambda_{\star})}{-\log(1-\eta s)}
    - 
    \frac{\log(1 + u/\lambda_{\star})}{-\log(1-\eta u)}\nonumber\\
    & \leq 
    \log(1 + s/\lambda_{\star})\frac{1 - \eta s}{\eta s }
    \sqrt{1/(1-\eta s) + (1/(1-\eta s) - 1)^2/6}
    - 
    \log(1 + u/\lambda_{\star})\frac{\sqrt{1- \eta u}}{\eta u } \nonumber 
    \\
    & = 
    \sqrt{1 - \eta u} 
    \left( \frac{\log(1 + s/\lambda_{\star})}{\eta s } 
    - 
    \frac{\log(1 + u/\lambda_{\star})}{\eta u }
    \right)  \label{tstardiff1}\\
    &\quad\quad 
    + 
    \log\left(1 + \frac{s}{\lambda_{\star}}\right) 
    \frac{1}{\eta s }
    \left( \sqrt{1 - \eta s } - \sqrt{1 - \eta u } \right)  
    \label{tstardiff2}\\
    & \quad\quad 
    + 
    \log\left(1 + \frac{s}{\lambda_{\star}}\right) 
    \frac{1}{\eta s }
    \left( (1 - \eta s) \sqrt{1/(1 - \eta s) + (1/(1-\eta s) - 1)^2/6} - \sqrt{1 - \eta s } 
    \right).\label{tstardiff3}
\end{align}
By the fundamental theorem of calculus, we can bound the term in \eqref{tstardiff1} as 
\begin{align*}
    \frac{\log(1 + s/\lambda_{\star})}{\eta s } 
    - 
    \frac{\log(1 + u/\lambda_{\star})}{\eta u }
    & = 
    \frac{1}{\eta} 
    \int_{s}^{u} \left[
    \frac{1}{z^2 }\log(1+z/\lambda_{\star}) - 
    \frac{1}{z + \lambda_{\star}} 
    \frac{1}{z} 
    \right]dz\\
    & \leq 
    \frac{1}{\eta} 
    \int_{s}^{u} \left[
    \frac{1}{z^2 }\frac{z}{\lambda_{\star}} \frac{1}{\sqrt{1 + z/\lambda_{\star}}} - 
    \frac{1}{z + \lambda_{\star}} 
    \frac{1}{z} \right]dz
    \\
    & = 
    \frac{1}{\eta} 
    \int_{s}^{u} \left[
    \frac{1}{\sqrt{z+ \lambda_{\star}}} 
    \frac{1}{z}
    \left( \frac{1}{\sqrt{\lambda_{\star}}}
    - 
    \frac{1}{\sqrt{z + \lambda_{\star}}} 
    \right) \right]dz\\
    & \leq 
    \frac{1}{2 \eta \lambda_{\star}} 
    \int_{s}^{u} \frac{1}{z + \lambda_{\star}} dz,
\end{align*}
where we use the second inequality in \eqref{equ:LogInqSharp}. The final inequality arises from using  $\sqrt{1 + x} - 1 \leq x/2 $ for $x \geq 0 $ to conclude 
$\frac{1}{\sqrt{\lambda_{\star}}}
    - \frac{1}{\sqrt{z + \lambda_{\star}}} = \frac{\sqrt{1 + z/\lambda_{\star}} - 1}{\sqrt{z + \lambda_{\star}}} \leq \frac{z}{2\lambda_{\star}} \frac{1}{\sqrt{z + \lambda_{\star}}}$. 
    
To bound the term in \eqref{tstardiff2} we have, since $u\geq s$,
\begin{align*}
    \sqrt{1 - \eta s} - \sqrt{1 - \eta u} 
    = 
    \frac{\eta (u-s)}{\sqrt{ 1- \eta s} + \sqrt{ 1 - \eta u}}
    \leq 
    \frac{\eta (u-s)}{\sqrt{1- \eta s}}.
\end{align*}
Finally, for the term \eqref{tstardiff3}, we note that 
\begin{align*}
     (1 - \eta s) \sqrt{1/(1 - \eta s) + (1/(1-\eta s) - 1)^2/6} - \sqrt{1 - \eta s } \\
     = 
    \sqrt{1 - \eta s} 
    \left( 
    \sqrt{ 1 + \frac{(\eta s)^2}{6(1 - \eta s)}} - 1 
    \right)
    & \leq 
    \frac{(\eta s)^2}{12\sqrt{ 1- \eta s}}.
\end{align*}
Bringing together the last three displays then yields \eqref{equ:GeneralDifference}.
\subsection{Proofs for Relative Sub-optimality with Fast Decaying Eigenvalues}
\label{sec:proof:Subopt:FastDecay}
In this section we present the proofs of the lower bounds on the relative sub-optimality when the eigenvalues decay either exponentially or as a power law with a power greater than one. We begin with a series of technical results.  The first shows that the error of gradient descent is decreasing for a particular number of iterations. 
\begin{lemma}
\label{lem:DecreasingGradientDescent}
Let  for $i=1,\dots,r$,  $b_i = \frac{\lambda_{\star}}{\lambda_{\star} + s_i}$ 
as well as 
\begin{align}\label{umin}
    u_{\min}(\eta) := 
    \arg\inf_{u > 2 \lambda_{\star}}\left\{A(u,\eta) \leq 1 - \exp\left(-\frac{1}{64} \right)  \right\},
\end{align}
where
\beq\label{Au}
 A(u,\eta) := 
    \frac{\sum_{i : s_i > u } b_i+\eta /2\cdot\sum_{i=1}^{r}b_i s_i (1 + \frac{s_i}{\lambda_{\star}})^2}{ \sum_{i : \lambda_{\star} < s_i \leq u /2 } b_i}.
\eeq
If $\eta \leq \frac{1}{64 s_1 \log(1 + \frac{s_1}{\lambda_{\star}}) } $  and $u_{\min}(\eta) \leq s_1$, then  the expected estimation error of gradient descent is decreasing up to $t  =  \lceil t^{\star}(u_{\min}(\eta)) \rceil$  iterations so that
\begin{align*}
    \E_{\beta_{\star},\epsilon}[L_{\beta_{\star}}(\hbeta_{\eta,t+1})]  \leq \E_{\beta_{\star},\epsilon}[L_{\beta_{\star}}(\hbeta_{\eta,t})]
     \quad\quad \text{ for }  t \leq t^{\star}(u_{\min}(\eta)).
\end{align*}
\end{lemma}
With slight abuse of notation, we will often abbreviate $A(u) =  A(u,\eta)$.

A consequence of this lemma is that the optimal number of iterations is at least $t^{\star}(u_{\min}(\eta))$. This allows us to prove that gradient descent essentially overfits the large eigenvalue directions, which provides the required lower bounds. 
\begin{theorem}
\label{thm:GDLower2}
Suppose $\eta \leq \frac{1}{64 s_1 \log(1 + s_1/\lambda_{\star})}$ and $t \geq \frac{1}{ \eta \lambda_{\star}}$. If $u_{\min}(\eta) < s_1/2$, 
then the error of gradient descent is lower bounded as
\begin{align*}
    \mathcal{E}(\mathcal{C}^{\text{GD}}(\eta,t)) \geq \frac{1}{4} (1 - \exp(- 32^{-1}))^2 \lambda_\star \int^{s_1}_{2 u_{\min}(\eta)} G^{\text{GD}}(s) d \widehat{H}(s),
\end{align*}
where $G^{\text{GD}}(s) $ is defined in Lemma \ref{lem:GDBounds:Fine}. 
\end{theorem}

Given this, it is clear that an upper bound on $u_{\min}(\eta)$ is sufficient, which then depends upon the stepsize and eigenvalue spectrum. We now present two propositions, which provide upper bounds when the eigenvalues decay exponentially or as a power law.   
\begin{proposition}
\label{prop:Expdecay}
Suppose $s_i = \exp(-\rho(i-1))$ for $i=1,\dots,r$ and $\rho > 0 $. If
\begin{align}\label{etaexp}
    \eta \leq \frac{\lambda_{\star} }{8e} \frac{e^{\rho}}{1 +  e^{\rho}}
    (1 - e^{-\frac{1}{64}}),
\end{align}
then $u_{\min}(\eta) \leq 2 e \lambda_{\star}  \left( 1 + 2 e^{\rho} (1 - e^{-\frac{1}{64}})^{-1} \right).$

\end{proposition}
We now consider the case of eigenvalues decaying at a fast power law rate. 
\begin{proposition}
\label{prop:fastPolynomial}
Suppose $s_i = i^{-\alpha}$ for $i=1,2,3,\dots,r$ and $\alpha > 1$. If $\lambda_{\star}\leq 2^{-\alpha}$ and
\begin{align}\label{etafp}
      \eta \leq 
      \frac{\lambda_{\star}^{1-1/\alpha}}{17\cdot 2^{\alpha}}
    \frac{\alpha - 1}{1+\alpha}
     (1-e^{-\frac{1}{64}}).
\end{align}
Then $u_{\min}(\eta)
    \leq 
    2^{1+\alpha}\lambda_{\star}\left[ 1 + 2^{4+3\alpha}(1 - e^{-\frac{1}{64}})^{-1} \right]^{\frac{\alpha}{1+\alpha}}.$

\end{proposition}
Given these technical results we are ready to the present the proofs, with the remainder of this section structured as follows. 
Section \ref{sec:Proof:thm:GDVsRidge:exp} presents the proof of Theorem \ref{thm:mainGDvsRidge_Upper:exp}, which gives the lower bound on the estimation error when the eigenvalues decay exponentially.  
Section \ref{sec:Proof:thm:GDVsRidge:poly} presents the proof of Theorem \ref{thm:mainGDvsRidge_Upper:poly}, which gives the lower bound on the estimation error when the eigenvalues decay as a power law. The remaining sections then give the proofs for the aforementioned technical results. 
Section \ref{sec:proof:DecreasingGradientDescent} presents the proof of Lemma  \ref{lem:DecreasingGradientDescent}.
Section \ref{sec:proof:GDLower2} presents the proof of Theorem \ref{thm:GDLower2}.
Section \ref{sec:proof:prop:Expdecay} presents the proof of Proposition \ref{prop:Expdecay}. 
Section  \ref{sec:proof:prop:fastPolynomial} presents the proof of Proposition \ref{prop:fastPolynomial}.

\subsubsection{Proof of Theorem \ref{thm:mainGDvsRidge_Upper:exp}}
\label{sec:Proof:thm:GDVsRidge:exp}
We begin by checking that the conditions of  Theorem \ref{thm:GDLower2} and Proposition \ref{prop:Expdecay} are satisfied. 
We recall $\eta = \frac{1}{k \lambda_{\min}}$ and condition \eqref{expgdc}. 
Considering the conditions in Proposition \ref{prop:Expdecay} we see from $k \geq \frac{8e(1 + e^{-\rho})}{\lambda_{\star}\lambda_{\min} (1-e^{-\frac{1}{64}})}$ that \eqref{etaexp} holds. As a consequence of Proposition \ref{prop:Expdecay} alongside that $\lambda_{\star} < \frac{s_1}{2} (2e + 4 e^{\rho+1}(1 - e^{-\frac{1}{64}})^{-1})^{-1}$  we have $u_{\min}(\eta) < s_1/2$. Moreover, the conditions of Theorem \ref{thm:GDLower2} are satisfied since  since $k \geq \frac{ 64 }{\lambda_{\min}} s_1 \log(1 + \frac{s_1}{\lambda_\star})$ implies $\eta \leq \frac{1}{64 s_1 \log(1+\frac{s_1}{\lambda_{\star}})}$.  

Applying Theorem \ref{thm:GDLower2} to gradient descent alongside the upper bound from Proposition \ref{lem:RidgeBounds} for ridge regression with $\delta = (1-\lambda_{\min})/(k-1)$  yields  
\begin{align*}
    \mathcal{S}(\mathcal{C}^{\text{GD}}(\eta,t),\mathcal{C}^{\text{Ridge}}(\lambda_{\min},k))
    & = 
    \frac{\mathcal{E}(\mathcal{C}^{\text{GD}}(\eta,t))}{\mathcal{E}(\mathcal{C}^{\text{Ridge}}(\lambda_{\min},k))}\\
    & \geq 
    \frac{1}{16} (1 - \exp(-\frac{1}{32}))^2
    \frac{1}{\min\{1-\ep,\ep\}^2}
    \frac{\int^{s_1}_{2 u_{\min}(\eta)} G^{\text{GD}}(s) d \widehat{H}(s)}{\int G^{\text{Ridge}}(s) d \widehat{H}(s)}.
\end{align*}
We now need to bound the rightmost integrals. From Proposition \ref{prop:Expdecay} we recall that $u_{\min}(\eta) \leq  \lambda_{\star} c_{\rho}$ where $c_{\rho} = 2e + 4 e^{1+\rho}(1- e^{-\frac{1}{64}})^{-1} > 1$. This allows us to lower bound
\begin{align*}
    \int^{s_1}_{2 u_{\min}(\eta)} G^{\text{GD}}(s) d \widehat{H}(s)
    \geq 
    \int^{s_1}_{2 c_{\rho} \lambda_{\star}} G^{\text{GD}}(s) d \widehat{H}(s)
    = \int^{s_1}_{2 c_{\rho} \lambda_{\star}} \frac{\lambda_{\star}}{s^2} d \widehat{H}(s).
\end{align*}
For ridge regression we have 
\begin{align*}
    \int G^{\text{Ridge}}(s) d \widehat{H}(s)
    & = 
    \left( \frac{1-\lambda_{\min}}{ \lambda_{\star} k}\right)^2
    \left[ 
    \int^{s_1}_{c_{\rho} \lambda_{\star}} \frac{\lambda_{\star}}{ s^2} d\widehat{H}(s) 
    + 
     \int^{c_{\rho}\lambda_{\star}}_{\lambda_{\star}} \frac{\lambda_{\star}}{s^2} d\widehat{H}(s)
    + 
    \int^{\lambda_{\star}}_{0}\frac{s}{\lambda_{\star}^{2}} d \widehat{H}(s)
    \right].
\end{align*}
Bringing everything together yields 
\begin{align}\label{thm5lb}
    \mathcal{S}(\mathcal{C}^{\text{GD}}(\eta,t),\mathcal{C}^{Ridge}(\lambda_{\min},k))
    &
    \geq 
    \frac{1}{16} (1 - \exp(-\frac{1}{32}))^2
    \frac{1}{\min\{1-\ep,\ep\}^2}
    \left( \frac{ \lambda_{\star} k}{1-\lambda_{\min}}\right)^2\nonumber
    \\
    &\quad\quad\times
    \frac{\int^{s_1}_{c_{\rho} \lambda_{\star}} \frac{\lambda_{\star}}{ s^2} d\widehat{H}(s)  }{\int^{s_1}_{c_{\rho} \lambda_{\star}} \frac{\lambda_{\star}}{ s^2} d\widehat{H}(s) 
    + 
    \int^{c_{\rho}\lambda_{\star}}_{\lambda_{\star}} \frac{\lambda_{\star}}{s^2} d\widehat{H}(s)
    + 
    \int^{\lambda_{\star}}_{0}\frac{s}{\lambda_{\star}^{2}} d \widehat{H}(s) }.
\end{align}
We must now bound the ratio of integrals. From $\exp(2\rho(i-1)) \leq \int_{i-1}^{i} \exp(2\rho x)dx $ we have 
\begin{align*}
    \int^{c_{\rho}\lambda_{\star}}_{\lambda_{\star}} \frac{\lambda_{\star}}{s^2} d\widehat{H}(s)
    & = \frac{1}{d} \sum_{i : \lambda_{\star} < s_i \leq c_{p} \lambda_{\star} } \frac{\lambda_{\star}}{s_i^2}
     = 
    \frac{\lambda_{\star}}{d} \sum_{i = \lceil 1 + \frac{1}{\rho} \log(1/(c_{p} \lambda_{\star})) \rceil }^{\lceil \frac{1}{\rho} \log(1/\lambda_{\star}) \rceil} \exp(2\rho(i-1))\\
    & \leq 
    \frac{\lambda_{\star}}{d}  \int_{\frac{1}{\rho} \log(1/(c_{p} \lambda_{\star}))}^{1+\frac{1}{\rho} \log(1/\lambda_{\star} )} \exp(2\rho x)dx 
     = \frac{\lambda_{\star}}{2\rho d} \left( \frac{e^{2\rho}}{\lambda_{\star}^2} - \frac{1}{c_p^2 \lambda_{\star}^2} \right)
    \leq 
\frac{e^{2\rho}}{2\rho d \lambda_{\star}}.
\end{align*}
Similarly we have 
\begin{align*}
\int^{\lambda_{\star}}_{0}\frac{s}{\lambda_{\star}^{2}} d \widehat{H}(s)
 \! = \!
\frac{1}{d\lambda_{\star}^2} \sum_{i=\lceil 1+\frac{1}{\rho} \log(1/\lambda_{\star})\rceil }^{r} \exp(-\rho(i-1))
 \leq \frac{1}{d\lambda_{\star}^2} \int_{\frac{1}{\rho} \log(1/\lambda_{\star})}^{r}\!\exp(-\rho x) dx 
\! \leq \! \frac{1}{d\rho \lambda_{\star}}.
\end{align*}
Finally, we have the lower bound 
\begin{align*}
    \int_{c_{\rho} \lambda_{\star}}^{s_1} \frac{\lambda_{\star}}{s^2} d\widehat{H}(s) 
    & = 
    \frac{\lambda_{\star}}{d} \sum_{i=1}^{\lceil  \frac{1}{\rho} \log(1/(c_{\rho} \lambda_{\star}))\rceil} 
    \exp(2\rho(i-1)) \\
    & \geq \frac{\lambda_{\star}}{d} 
    \int_0^{\frac{1}{\rho} \log(1/(c_{\rho} \lambda_{\star}))}
    \exp(2\rho x) dx 
    = 
    \frac{\lambda_{\star}}{2\rho d}
    \left( \frac{1}{c_\rho^2 \lambda_{\star}^2 } - 1\right)
     \geq \frac{1}{4d \rho c_{\rho}^2 \lambda_{\star}},
\end{align*}
where in the last inequality we used that $\lambda_{\star} \leq \frac{s_1}{2c_{\rho}}  \leq \frac{1}{\sqrt{2}c_{\rho}}$. 
Combining these bounds we have 
\begin{align*}
    \frac{\int^{s_1}_{c_{\rho} \lambda_{\star}} \frac{\lambda_{\star}}{ s^2} d\widehat{H}(s)  }{\int^{s_1}_{c_{\rho} \lambda_{\star}} \frac{\lambda_{\star}}{ s^2} d\widehat{H}(s) 
    + 
    \int^{c_{\rho}\lambda_{\star}}_{\lambda_{\star}} \frac{\lambda_{\star}}{s^2} d\widehat{H}(s)
    + 
    \int^{\lambda_{\star}}_{0}\frac{s}{\lambda_{\star}^{2}} d \widehat{H}(s) }
    \geq \frac{1}{1 + 4 c_{\rho}^2(1 + e^{2\rho}) }
    \geq 
    \frac{1}{c_{\rho}^2} \frac{1}{8 (1 + e^{2\rho})}.
\end{align*}
Plugging this into \eqref{thm5lb} yields the result.

\subsubsection{Proof of Theorem \ref{thm:mainGDvsRidge_Upper:poly}}
\label{sec:Proof:thm:GDVsRidge:poly}
Similar to the proof of Theorem \ref{thm:mainGDvsRidge_Upper:exp} we begin by checking the conditions of Theorem \ref{thm:GDLower2} and Proposition \ref{prop:fastPolynomial}. Since $\eta = \frac{1}{k \lambda_{\min}}$ and \eqref{thm6cond1}, \eqref{thm6cond2} hold, \eqref{etafp} holds. Since $\lambda_\star < \frac{s_1}{2} \frac{1}{2^{1+\alpha}} [ 1 + 2^{4+3\alpha}(1-e^{-\frac{1}{64}})^{-1}]^{\frac{-\alpha}{\alpha+1}}$ we have from Proposition \ref{prop:fastPolynomial} that $u_{\min}(\eta) < s_1/2$. Moreover, since $k \geq \frac{64}{\lambda_{\min}} s_1 \log(1 + \frac{s_1}{\lambda_{\star}})$, we have $\eta \leq \frac{1}{64 s_1 \log(1+\frac{s_1}{\lambda_{\star}})}$, and thus the conditions for Theorem  \ref{thm:GDLower2} are met.

Writing $u_{\min}(\eta) \leq \lambda_{\star} b_{\rho}$ where $b_{\rho} =  2^{1+\alpha}(1+2^{4+3\alpha}(1-e^{-\frac{1}{64}})^{-1})^{\frac{\alpha}{\alpha+1}} > 1$,
we have
\begin{align}
    \mathcal{S}(\mathcal{C}^{\text{GD}}(\eta,t),\mathcal{C}^{\text{Ridge}}(\lambda_{\min},k))
    &
    \geq 
    \frac{1}{16} (1 - \exp(-\frac{1}{32}))^2
    \frac{1}{\min\{1-\ep,\ep\}^2}
    \left( \frac{ \lambda_{\star} k}{1-\lambda_{\min}}\right)^2
    \nonumber
    \\
    &\quad\quad\times
    \frac{\int^{s_1}_{b_{\rho} \lambda_{\star}} \frac{\lambda_{\star}}{ s^2} d\widehat{H}(s)  }{\int^{s_1}_{b_{\rho} \lambda_{\star}} \frac{\lambda_{\star}}{ s^2} d\widehat{H}(s) 
    + 
    \int^{b_{\rho}\lambda_{\star}}_{\lambda_{\star}} \frac{\lambda_{\star}}{s^2} d\widehat{H}(s)
    + 
    \int^{\lambda_{\star}}_{0}\frac{s}{\lambda_{\star}^{2}} d \widehat{H}(s) }.\label{thm6lb}
\end{align}
We now bound the ratio of integrals in a similar way to the argument in the proof of Theorem \ref{thm:mainGDvsRidge_Upper:exp}. First,  
\begin{align*}
    &\int^{s_1}_{b_{\rho} \lambda_{\star}} \frac{\lambda_{\star}}{ s^2} d\widehat{H}(s)
     = 
    \frac{ \lambda_{\star}}{d} \sum_{i : b_{\rho} \lambda_{\star} < s_i \leq s_1} 
    i^{2\alpha} 
     = 
    \frac{ \lambda_{\star}}{d} \sum_{i=1}^{\lceil (b_{\rho} \lambda_{\star})^{-1/\alpha} \rceil-1 } 
    i^{2\alpha} 
     \geq 
    \frac{ \lambda_{\star}}{d} \int_{0}^{\lceil (b_{\rho} \lambda_{\star})^{-1/\alpha} \rceil-1 }
    x^{2\alpha}dx\\
    & \geq 
    \! \frac{ \lambda_{\star}}{d} \! \int_{0}^{ (b_{\rho} \lambda_{\star})^{-1/\alpha} - 1  }
    \! x^{2\alpha}dx 
     \! = \!
     \frac{1}{1 \! + \! 2\alpha} 
    \frac{\lambda_{\star}}{d} 
    ((b_{\rho} \lambda_{\star})^{-1/\alpha} \! - \! 1)^{1+2\alpha}
    \! \geq \!
     \frac{ 1 }{1+2\alpha} \frac{1}{2^{1+2\alpha}} \frac{\lambda_{\star}}{d} (b_{\rho} \lambda_{\star})^{-(2+1/\alpha)},
\end{align*}
where in the last inequality we used that $\lambda_{\star} \leq \frac{1}{2^{\alpha} b_{\rho}}$. 
For the second integral we have 
\begin{align*}
    \int^{b_{\rho}\lambda_{\star}}_{\lambda_{\star}} \frac{\lambda_{\star}}{s^2} d\widehat{H}(s)
    & = 
    \frac{ \lambda_{\star}}{d} \sum_{i : \lambda_{\star}< s_i \leq b_{\rho}  \lambda_{\star} } 
    i^{2\alpha} 
     = 
    \frac{ \lambda_{\star}}{d}  \sum_{i= \lceil (b_{\rho} \lambda_{\star})^{-1/\alpha} \rceil}^{\lceil  \lambda_{\star}^{-1/\alpha}\rceil-1 }
    i^{2\alpha}
     \leq 
    \frac{ \lambda_{\star}}{d} 
    \int_{(b_{\rho} \lambda_{\star})^{-1/\alpha}}^{1 + \lambda_{\star}^{-1/\alpha}}
    x^{2\alpha} dx \\
    & \leq \frac{ \lambda_{\star}}{d} 
    \frac{1}{1 + 2\alpha} 
    2^{1+2\alpha}\left(\lambda_{\star}^{-(1+2\alpha)/\alpha} - (b_{\rho}\lambda_{\star} )^{-(2+1/\alpha)} \right)\\
    & \leq  
    \frac{ \lambda_{\star}}{d} 2^{1+2\alpha}
    \lambda_{\star}^{-(2+1/\alpha)}.
\end{align*}
Finally, since $\alpha > 1$, we have
\begin{align*}
    \int^{\lambda_{\star}}_{0}\frac{s}{\lambda_{\star}^{2}} d \widehat{H}(s)
    &= \frac{1}{d\lambda_{\star}^2} \sum_{i: s_r \leq s_i \leq \lambda_{\star}} i^{-\alpha} 
     = 
    \frac{1}{d\lambda_{\star}^2}
    \sum_{i=\lceil \lambda_{\star}^{-1/\alpha} \rceil}^{r}
    i^{-\alpha}
     \leq 
    \frac{1}{d\lambda_{\star}^2} \int_{\lambda_{\star}^{-1/\alpha}}^{ r+ 1} x^{-\alpha}dx \\
    & = 
    \frac{1}{(\alpha - 1)d\lambda_{\star}^2} \left( \lambda_{\star}^{1-1/\alpha} -(r+1)^{1-\alpha} \right) 
    \leq 
     \frac{1}{(\alpha - 1)d\lambda_{\star}^{1+1/\alpha}}.
\end{align*}
Bringing everything together we have 
\begin{align*}
    \frac{\int^{s_1}_{b_{\rho} \lambda_{\star}} \frac{\lambda_{\star}}{ s^2} d\widehat{H}(s)  }{\int^{s_1}_{b_{\rho} \lambda_{\star}} \frac{\lambda_{\star}}{ s^2} d\widehat{H}(s) 
    + 
    \int^{b_{\rho}\lambda_{\star}}_{\lambda_{\star}} \frac{\lambda_{\star}}{s^2} d\widehat{H}(s)
    + 
    \int^{\lambda_{\star}}_{0}\frac{s}{\lambda_{\star}^{2}} d \widehat{H}(s) }
    & \geq 
    \frac{1}{1 + 2^{2(1+2\alpha)}(1+2\alpha)b_{\rho}^{1/\alpha + 2}\frac{\alpha}{\alpha - 1}}\\
    & \geq 
    \frac{1}{2} \frac{1}{ 2^{2(1+2\alpha)}(1+2\alpha)b_{\rho}^{1/\alpha + 2}\frac{\alpha}{\alpha - 1}}.
\end{align*}
Plugging in the above bound into \eqref{thm6lb} yields the result.

\subsubsection{Proof of Lemma \ref{lem:DecreasingGradientDescent}}
\label{sec:proof:DecreasingGradientDescent}
Fix some number of iterations $t \geq 0$. The error decreases with an additional iteration if and only if the following quantity is negative:
\begin{align*}
    & \frac{n}{\sigma^2}\left(  \E_{\beta_{\star},\epsilon}[L_{\beta_{\star}}(\hbeta_{\eta,t+1})] - \E_{\beta_{\star},\epsilon}[L_{\beta_{\star}}(\hbeta_{\eta,t})]\right)
    \\
    & =
    \sum_{i=1}^{r}\left(\frac{1}{s_i} + \frac{1}{\lambda_{\star}}\right) 
    \left( (1-\eta s_i)^{t+1} - \frac{\lambda_{\star}}{\lambda_{\star} + s_i} \right)^2 
    - 
    \sum_{i=1}^{r}\left(\frac{1}{s_i} + \frac{1}{\lambda_{\star}}\right) 
    \left( (1-\eta s_i)^{t} - \frac{\lambda_{\star}}{\lambda_{\star} + s_i} \right)^2 
    \\
    & = 
    \sum_{i=1}^{n}\left(\frac{1}{s_i}+ \frac{1}{\lambda_{\star}}\right)
    (1-\eta s_i)^{t}(\eta s_i)
    \left[ 2 \frac{\lambda_{\star}}{\lambda_{\star} + s_i} - (1-\eta s_i)^{t}(2-\eta s_i)
    \right] =:D.
\end{align*}
For a $t = t^{\star}(u)$ for some $s_1 > u > 0$, we now bound the term within the square brackets by considering three different cases for the eigenvalues : (A) $\lambda_{\star} < s_i \leq u/2$; (B) $u/2 < s_i \leq u$ and (C)  $u< s_i$. 

Recall the function $t^{\star}$ from \eqref{tstar}.
For case (A) observe that we have the lower bound  for $\lambda_{\star} < s_i \leq u/2$
\begin{align}
\label{equ:LowerBound}
    \eta s_i(t^{\star}(s_i) - t^{\star}(u)) \geq 
    \frac{1}{64},
\end{align}
the proof of which is given at the end of this subsection, starting with display \eqref{t-star-bound}. 
Using this we have, for $\lambda_{\star} < s_i \leq u/2$,
\begin{align*}
    \frac{\lambda_{\star}}{\lambda_{\star} + s_i} - (1-\eta s_i)^{t^{\star}(u)}
    & = 
    (1- \eta s_i)^{t^{\star}(s_i)} - (1-\eta s_i)^{t^{\star}(u)} \\
    & = 
    (1-\eta s_i)^{t^{\star}(u)}
    \left( (1- \eta s_i)^{t^{\star}(s_i)- t^{\star}(u)} - 1 \right) \\
    & \leq 
    (1-\eta s_i)^{t^{\star}(u)} \left( \exp(- \eta s_i( t^{\star}(s_i)- t^{\star}(u)) ) - 1 \right)\\
    & \leq 
    \frac{\lambda_{\star}}{\lambda_{\star} + s_i}\left(\exp\left(-\frac{1}{64}\right) - 1\right).
\end{align*}
Meanwhile for case (B), i.e.,  $u \geq s_i > u/2 $, we immediately have from $s \rightarrow t^{\star}(s)$ being decreasing that  $\frac{\lambda_{\star}}{\lambda_{\star} + s_i} - (1-\eta s_i)^{t^{\star}(u)} = (1-\eta s_i)^{t^{\star}(s_i)} - (1-\eta s_i)^{t^{\star}(u)}\leq 0$. Finally for case (C), i.e., $s_i > u$, we have 
\begin{align*}
     \frac{\lambda_{\star}}{\lambda_{\star} + s_i } - (1-\eta s_i)^{t^{\star}(u)} \leq 
    \frac{\lambda_{\star}}{\lambda_{\star} + s_i }.
\end{align*}
Combining these three bounds and denoting $a_i(u) = (1-\eta s_i)^{t^{\star}(u)}$ yields
\begin{align*}
   D \leq 
    2 \sum_{i : s_i > u } \eta a_i(u) 
    + 
    2 (\exp(-\frac{1}{64}) - 1) \sum_{i : \lambda_{\star} <s_i \leq u/2 }\eta a_i(u) 
    + 
    \sum_{i =1}^{r} \eta^2 a_i(u) s_i (1+s_i/\lambda_{\star}).
\end{align*}
The above is then negative when 
\begin{align}\label{sum2frac}
    \frac{\sum_{i : s_i > u } a_i(u)}{ \sum_{i : \lambda_{\star} < s_i \leq u /2 } a_i(u)} 
    + 
    \eta \frac{\sum_{i=1}^{r}a_i(u) s_i (1 + \frac{s_i}{\lambda_{\star}})}{ 2\sum_{i : \lambda_{\star} < s_i \leq u /2 } a_i(u) }
    \leq 
    1 - \exp\left(-\frac{1}{64}\right).
\end{align}
We can further upper bound the terms $a_i(u)$. In particular for $s_i > u$ we have $(1-\eta s_i)^{t^{\star}(u)} \leq \frac{\lambda_{\star}}{\lambda_{\star} +s_i}$, and for $s_i \leq u$ we have $(1-\eta s_i)^{t^{\star}(u)} \geq \frac{\lambda_{\star}}{\lambda_{\star} + s_i}$. With 
$A(u)$ from \eqref{Au},
it is sufficient that
$
    A(u)
    \leq 
    1 - \exp(-\frac{1}{64} ),
$
where for the numerator of the second term in \eqref{sum2frac} we used that $(1-\eta s_i)^{t^{\star}(u)} \leq 1$ to  bound $\sum_{i=1}^{r} a_i(u) s_i (1 + \frac{s_i}{\lambda_{\star}}) \leq \sum_{i=1}^{r} b_i s_i (1 + \frac{s_i}{\lambda_{\star}})^2$. 

Since  $t^{\star}$ is decreasing in $u$, for any $t \leq t^{\star}(u_{\min}(\eta))$ we have $u \geq u_{\min}(\eta)$. Since $A(u)$ is decreasing, from \eqref{umin} it follows that the required condition $
    A(u)
    \leq 
    1 - \exp(-\frac{1}{64} )
$ holds, concluding the proof.

{\bf Proof of \eqref{equ:LowerBound}.} We now prove the inequality \eqref{equ:LowerBound}, writing $s=s_i$ for simplicity.  Using $\frac{x}{1+x} \leq \log(1+x) \leq x$ we immediately have
\begin{align}
    & t^{\star}(s) - t^{\star}(u)
    = 
    \frac{\log\left(1 +\frac{s}{\lambda_{\star}}\right)}{-\log(1-\eta s)}
    - 
    \frac{\log\left(1 +\frac{u}{\lambda_{\star}}\right)}{-\log(1-\eta u)}\label{t-star-bound}\\
    & \geq 
    \frac{1-\eta s}{\eta s} \log\big(1 \! + \!\frac{s}{\lambda_{\star}}\big)
    - \frac{1}{\eta u } \log\big(1 + \frac{u}{\lambda_{\star}}\big)
    = 
    \frac{1}{\eta} \left[ \frac{\log\big(1 \! + \! \frac{s}{\lambda_{\star}}\big)}{s} \! - \! \frac{\log\big(1 \!+ \! \frac{u}{\lambda_{\star}}\big)}{u}
    \right]
    \! - \!
    \log\big(1 \! + \!\frac{s}{\lambda_{\star}}\big)\label{t-star-bound2}.
\end{align}
We then have, using the first inequality in \eqref{equ:LogInqSharp}, and with $S(z) = \sqrt{1 + \frac{z}{\lambda_{\star}} + \frac{1}{6} \left( \frac{z}{\lambda_{\star}}\right)^2}$  
\begin{align}
    &\frac{\log\left(1 +\frac{s}{\lambda_{\star}}\right)}{s} - \frac{\log\left(1 +\frac{u}{\lambda_{\star}}\right)}{u}
     = 
    \int_{s}^{u}\frac1z\left[\frac{1}{z} \log(1 + \frac{z}{\lambda_{\star}}) - \frac{1}{z + \lambda_{\star}} \right]dz \nonumber\\
    & \geq 
    \int_{s}^{u}\frac1z\left[\frac{1}{\lambda_{\star} S(z)}  - \frac{1}{z + \lambda_{\star}} \right]dz
     = 
    \int_{s}^{u}\frac{dz}{z(z+\lambda_{\star})}\left[ 
    \frac{1 + \frac{z}{\lambda_{\star}}}{S(z)} - 1 \right].
    \label{log-diff-bd}
\end{align}
However,
\begin{align}\label{simple-ineq}
    \frac{1+x}{\sqrt{1 + x + x^2/6}} 
    \geq 
    1 + \frac{1}{8} \min(x,1).
\end{align}
To see this, split into cases $x \in [0,1]$ and $x \in (1,\infty)$. For the case $x \in [0,1]$ consider the function $f(x) = \frac{1+x}{\sqrt{1+x+x^2/6}} - 1 - \frac{1}{8} x$. Observe that $f(0) = 0$. Computing the derivative we have $f^{\prime}(x) = \frac{1}{2} \frac{1+x/6}{(1+x+x^2/6)^{3/2}} - \frac{1}{8}$. For $0 \leq x \leq 1$ we then see that $f^{\prime}(x) \geq \frac{1}{2} \frac{1}{(2+1/6)^{3/2}} - \frac{1}{8} \geq 0$. Therefore $f(x)$ is non-negative on $[0,1]$. Meanwhile for $x \in [1,\infty)$ consider $g(x) = \frac{1+x}{\sqrt{1+x+x^2/6}} - 1 - \frac{1}{8} $. We then see that $g(1) = \frac{2}{\sqrt{2 + 1/6}} - 1 - \frac{1}{8} \geq 0$, while from the calculation just above this, $g^{\prime}(x)=\frac{1}{2} \frac{1+x/6}{(1+x+x^2/6)^{3/2}}$, which
is non-negative. Plugging in \eqref{simple-ineq} into \eqref{log-diff-bd} we have
\begin{align*}
    \frac{\log\left(1 +\frac{s}{\lambda_{\star}}\right)}{s} - \frac{\log\left(1 +\frac{u}{\lambda_{\star}}\right)}{u}
    & \geq 
    \frac{1}{8} 
    \int_{s}^{u} \frac{1}{z(z+\lambda_{\star})} \min(1,\frac{z}{\lambda_{\star}})dz.
\end{align*}
Recalling that $s > \lambda_{\star}$ in the case considered here then yields 
\begin{align*}
    \frac{\log\left(1 +\frac{s}{\lambda_{\star}}\right)}{s} - \frac{\log\left(1 +\frac{u}{\lambda_{\star}}\right)}{u}
    & \geq 
    \frac{1}{16} \int_{s}^{u} \frac{1}{ z^2} dz
    = 
    \frac{1}{16} \left[ \frac{1}{s} - \frac{1}{u} \right].
\end{align*}
Plugging this into \eqref{t-star-bound2} and using that in the case considered here $u \geq 2s$ and by assumption $\eta s \log(1 +\frac{s}{\lambda_{\star}}) \leq \eta s_1 \log(1 +\frac{s_1}{\lambda_{\star}}) \leq \frac{1}{64}$ yields 
\begin{align*}
    \eta s 
    \left(t^{\star}(s) - t^{\star}(u)\right) 
    \geq 
    \frac{1}{16} \left[ \frac{s}{s} - \frac{s}{u}\right] 
    - 
    \eta s \log\left(1 + \frac{s}{\lambda_{\star}}\right)
    & \geq 
    \frac{1}{32} - \eta s \log\left(1 + \frac{s}{\lambda_{\star}}\right)
     \geq 
    \frac{1}{64},
\end{align*}
as required.

\subsubsection{Proof of Theorem \ref{thm:GDLower2}}
\label{sec:proof:GDLower2}
From Lemma \ref{lem:DecreasingGradientDescent}, and recalling  \eqref{tstar}, we know the error is decreasing for $t \leq t^{\star}(u_{\min})$. Therefore the optimal number of iterations is greater than or equal to $t^{\star}(u_{\min})$. Thus  we have
\begin{align}
    \mathcal{E}(\mathcal{C}^{\text{GD}}(\eta,t))
    & = 
    \min_{ j \in \{1,\dots,t\}} 
    \left( \E_{\beta_{\star},\epsilon}[L_{\beta_{\star}}(\hbeta_{\eta,j})] - \E_{\beta_{\star},\epsilon}[L_{\beta_{\star}}(\hbeta_{\lambda_\star})] \right) \nonumber 
    \\
    & = 
    \min_{j \in \{\lceil t^{\star}(u_{\min}(\eta) )\rceil,\dots,t \} } 
    \left( 
    \E_{\beta_{\star},\epsilon}[L_{\beta_{\star}}(\hbeta_{\eta,j})]
    - 
    \E_{\beta_{\star},\epsilon}[L_{\beta_{\star}}(\hbeta_{\lambda_\star})] \right) \nonumber 
    \\
    & \geq 
    \min_{j \in \{\lceil t^{\star}(u_{\min}(\eta) )\rceil,\dots,t \} }
    \frac{\sigma^2}{n}  \sum_{i : s_i > u_{\min}(\eta) }
    \frac{1}{\lambda_{\star}} \left( (1-\eta s_i)^{j} - \frac{\lambda_{\star}}{\lambda_{\star} + s_i}\right)^2\nonumber\\
    & \geq 
    \frac{\sigma^2}{n} \sum_{i : s_i > u_{\min}(\eta) }
    \frac{1}{\lambda_{\star}} \left( \frac{\lambda_{\star}}{\lambda_{\star} + s_i} - (1 - \eta s_i)^{t^{\star}(u_{\min}(\eta))}\right)^2\nonumber
    \\
    & \geq 
    \frac{\sigma^2}{n} \sum_{i : s_i > 2 u_{\min}(\eta) }
    \frac{1}{\lambda_{\star}} \left( \frac{\lambda_{\star}}{\lambda_{\star} + s_i} - (1 - \eta s_i)^{t^{\star}(u_{\min}(\eta))}\right)^2,\label{gdlb7}
\end{align}
where we use that $\frac{\lambda_{\star}}{\lambda_{\star} +s_i} = (1-\eta s_i)^{t^{\star}(s_i)} \geq (1-\eta s_i)^{t^{\star}(u_{\min}(\eta))} \geq (1-\eta s_i)^{j}$ for $s_i > u_{\min}(\eta)$ and $j \geq \lceil t^{\star}(u_{\min}(\eta) )\rceil$. Since $\lambda_{\star} < u_{\min}(\eta) < s_i/2 $ we can then lower bound using inequality \eqref{equ:LowerBound}, which implies $    \eta u_{\min}(\eta)(t^{\star}(u_{\min}(\eta)) - t^{\star}(s_i)) \geq \frac{1}{64}$ and thus also  $\eta s_i (t^{\star}(u_{\min}(\eta)) - t^{\star}(s_i)) \geq \frac{1}{32}$
\begin{align*}
    \frac{\lambda_{\star}}{\lambda_{\star} + s_i } - (1-\eta s_i)^{t^{\star}(u_{\min}(\eta))}
    & = 
    (1- \eta s_i)^{t^{\star}(s_i)}(1 - (1-\eta s_i)^{t^{\star}(u_{\min}(\eta)) - t^{\star}(s_i)})\\
    & \geq 
    (1- \eta s_i)^{t^{\star}(s_i)}(1 - \exp(- \eta s_i ( t^{\star}(u_{\min}(\eta)) - t^{\star}(s_i))))\\
    & \geq 
    (1- \eta s_i)^{t^{\star}(s_i)}\left(1 - \exp\left(- \frac{1}{32}\right)\right)\\
    & = 
    \frac{\lambda_{\star}}{\lambda_{\star} + s_i}\left(1 - \exp\left(- \frac{1}{32}\right)\right)
     \geq 
    \frac{1}{2} \frac{\lambda_{\star}}{s_i}\left(1 - \exp\left(- \frac{1}{32}\right)\right).
\end{align*}
Plugging in this lower bound into \eqref{gdlb7} then yields the result.

\subsubsection{Proof of Proposition \ref{prop:Expdecay}}
\label{sec:proof:prop:Expdecay}

Since $s_i = \exp(-\rho(i-1))$, we have the following upper bound, using that $\exp(\rho(i-1)) \leq \int_{i-1}^{i} \exp(\rho x ) dx$ and recalling the definition of $b_i$ and that $b_i \le  \frac{\lambda_{\star}}{s_i}$: 
\begin{align*}
    \sum_{i: s_i > u} b_i 
    & \leq  \lambda_{\star} \sum_{i=1}^{\lceil \frac{1}{\rho} \log( 1/u)\rceil } \exp( \rho(i-1))\\
    & \leq 
    \lambda_{\star} \int_0^{1 + \frac{1}{\rho} \log( 1/u)}\exp(\rho x)dx 
    =  \frac{\lambda_{\star}}{\rho}\left( \frac{e^{\rho}}{u} - 1\right)
     \leq \frac{ \lambda_{\star} e^{\rho}}{\rho u}.
\end{align*}

Using that $\exp(\rho(i-1))\geq \int_{i-2}^{i-1} \exp(\rho x) dx$, as well as that $s_i > \lambda_{\star}$---which implies $b_i \geq \frac{1}{2} \frac{\lambda_{\star}}{s_i}$---we have
\begin{align*}
    \sum_{i : \lambda_{\star} < s_i \leq u/2} b_i 
    & \geq  \frac{\lambda_{\star}}{2} \sum_{i = \lceil 1 + \frac{1}{\rho} \log(2/u)\rceil }^{\lceil  \frac{1}{\rho} \log(1/\lambda_{\star})\rceil } \exp(\rho(i-1)) 
    \geq  \frac{\lambda_{\star}}{2} \int_{ \lceil 1 + \frac{1}{\rho} \log(2/u)\rceil -2 }^{\lceil  \frac{1}{\rho} \log(1/\lambda_{\star})\rceil -1 }\exp(x\rho) dx \\
    & \geq 
     \frac{\lambda_{\star}}{2} 
    \int_{\frac{1}{\rho} \log(2/u) }^{\frac{1}{\rho} \log(1/\lambda_\star) - 1}
    \exp(x\rho) dx
    =  \frac{\lambda_{\star}}{2}\frac{1}{\rho}\left( \frac{1}{e \lambda_{\star}} - \frac{2}{u } \right) = 
    \frac{1}{2\rho}\left( e^{-1} - \frac{2 \lambda_{\star}}{u } \right).
\end{align*}
Taking the ratio of the bounds in the above two displays we then have 
\begin{align*}
    \frac{\sum_{i: s_i > u} b_i }{\sum_{i : \lambda_{\star} < s_i \leq u/2} b_i }
    \leq 
    2 e^{\rho +1} \frac{\lambda_{\star}}{u - 2 e\lambda_{\star}}.
\end{align*}

For the second term in the numerator of \eqref{Au} we have, since $s_i \leq 1$ and $\lambda_{\star} \leq 1$,
\begin{align*}
    \sum_{i=1}^{r} b_i s_i (1 + \frac{s_i}{\lambda_{\star}})^2
    & = \sum_{i=1}^{r} s_i^2 \left(\frac{1}{s_i} + \frac{1}{\lambda_\star}\right)
     \leq 
     \frac{2}{\lambda_{\star}}\sum_{i =1}^{r} s_i
    = 
    \frac{2}{\lambda_{\star}} \sum_{1=1}^{r} \exp(-\rho(i-1)) \\
    & \leq \frac{2}{\lambda_{\star}} \left(1 + \int_0^{r} \exp(-\rho x) dx \right)
    < \frac{2}{\lambda_{\star}} \left(1 + \int_0^{\infty} \exp(-\rho x) dx \right)
    = \frac{2}{\lambda_{\star}\rho}.
\end{align*}
Using this bound and the lower bound for $\sum_{i : \lambda_{\star} < s_i \leq u/2} b_i $ from above we then have
\begin{align*}
    \frac{ \sum_{i=1}^{r} b_i s_i (1 + \frac{s_i}{\lambda_{\star}})^2  }{\sum_{i : \lambda_{\star} < s_i \leq u/2} b_i}
    \leq 
    \frac{ \frac{2}{\lambda_{\star} \rho}}{\frac{1}{2\rho}(\frac{1}{e} - \frac{2 \lambda_{\star}}{u})}
    = 
     \frac{1}{\lambda_{\star}}\frac{4 e}{1 - \frac{2e \lambda_\star}{u}}
\end{align*}
and bringing everything together yields 
\begin{align}\label{uubfp}
    A(u)
    \leq  2 e^{\rho+1} \frac{\lambda_{\star}}{u - 2 e\lambda_{\star}}  + \frac{\eta}{\lambda_{\star}} \frac{4e}{1 - \frac{2e\lambda_\star}{u}}.
\end{align}

For $u = 2 e \lambda_{\star}  \left( 1 + 2 e^{\rho} (1 - e^{-\frac{1}{64}})^{-1} \right)$, we have 
\begin{align*}
    2 e^{\rho+1} \frac{\lambda_{\star}}{u - 2e \lambda_{\star}}
    = 
    \frac{1}{2} \left(1 - \exp\left(-\frac{1}{64}\right) \right).
\end{align*}
Moreover,
by \eqref{etaexp}, 
\begin{align*}
\frac{\eta}{\lambda_{\star}} \frac{4e}{1 - \frac{2e\lambda_\star}{u}}
\leq
\frac{1}{2} (1 - e^{-\frac{1}{64}})
\frac{e^{\rho}}{1 +  e^{\rho}} 
   =
\frac{1}{2} \left(1 - \exp\left(-\frac{1}{64}\right) \right).
\end{align*}
Thus, from \eqref{uubfp} and \eqref{umin}, the conclusion follows.

\subsubsection{Proof of Proposition \ref{prop:fastPolynomial}}
\label{sec:proof:prop:fastPolynomial}
This argument follows the proof of Proposition \ref{prop:Expdecay}. 
Since $s_i = i^{-\alpha}$, we have, using that $i^{\alpha} \leq \int_i^{i+1} x^{\alpha} dx$
\begin{align*}
    \sum_{i : s_i > u} b_i 
    &\leq 
    \lambda_{\star} \sum_{i=1}^{\lceil  u^{-1/\alpha}\rceil -1 } i^{\alpha}
    \leq 
     \lambda_{\star} 
    \int_{1}^{1+u^{-1/\alpha}} x^{\alpha}dx \\
    & = 
    \frac{ \lambda_{\star}}{\alpha+1}
    \left( ( 1+ u^{-1/\alpha} )^{1+\alpha} - 1 \right)
     \leq 
    \frac{2^{1+\alpha} \lambda_{\star} }{\alpha+1} u^{-(1+1/\alpha)}.
\end{align*}
Meanwhile using the lower bounds $i^{\alpha} \geq \int_{i-1}^{i} x^{\alpha} dx$ and $b_i \geq \frac{1}{2} \frac{\lambda_{\star}}{s_i}$ 
we have 
\begin{align*}
     &\sum_{i : \lambda_{\star} < s_i \leq u/2} b_i \geq 
    \frac{\lambda_{\star}}{2} \sum_{i = \lceil (u/2)^{-1/\alpha}\rceil }^{\lceil \lambda_{\star}^{-1/\alpha}\rceil -1} 
    i^{\alpha}
    \geq 
    \frac{\lambda_{\star}}{2} 
    \int_{\lceil (u/2)^{-1/\alpha}\rceil -1 }^{\lceil \lambda_{\star}^{-1/\alpha}\rceil -1 }
    x^{\alpha}dx \geq 
    \frac{\lambda_{\star}}{2} 
    \int_{(u/2)^{-1/\alpha} }^{  \lambda_{\star}^{-1/\alpha} -1 }
    x^{\alpha}dx\\
    & \geq 
    \frac{\lambda_{\star}}{2(\alpha\!+\!1)}
    \!
    \left( (\lambda_{\star}^{-1/\alpha} \! - \! 1)^{\alpha+1}
    \! - \! (u/2)^{-(1+1/\alpha)}
    \right)
    \! \geq \! 
    \frac{\lambda_{\star}}{2^{2+\alpha}(\alpha \! + \! 1)}
    ( \lambda_{\star}^{-(1+1/\alpha)} \! - \! 2^{2+\alpha + 1/\alpha}u^{-(1+1/\alpha)}),
\end{align*}
where in the last inequality we have used that $\lambda_\star^{-1/\alpha} -1 \geq \frac{1}{2} \lambda_{\star}^{-1/\alpha}$ since $\lambda_\star \leq  2^{-\alpha} $. Since $s_i \leq 1$ and $\lambda_{\star} \leq 1$
\begin{align*}
    \sum_{i=1}^{r} b_i s_i (1 + \frac{s_i}{\lambda_{\star}})^2 
    & \leq 
    \frac{2}{\lambda_{\star}} \sum_{i=1}^{r} s_i
     \leq \frac{2}{\lambda_{\star}(\alpha - 1)} (1 - r^{1-\alpha})
     \leq 
    \frac{2}{\lambda_{\star}(\alpha - 1)}.
\end{align*}
Combining these bounds we have 
\begin{align}
   A(u)& \leq 
    \frac{2^{3+2\alpha} \lambda_{\star}^{1+1/\alpha} + \frac{\eta}{\lambda_{\star}} 2^{3+\alpha} \frac{\alpha+1}{\alpha - 1} \lambda_{\star}^{1/\alpha} u^{1+1/\alpha}}
    {u^{1+1/\alpha} -  2^{2+\alpha+1/\alpha}\lambda_{\star}^{1+1/\alpha} }.
    \label{uub2}
\end{align}
Now for $u = 2^{1+\alpha}\lambda_{\star}\left( 1 + 2^{4+3\alpha}(1 - e^{-\frac{1}{64}})^{-1} \right)^{\alpha/(1+\alpha)}$ we have 
\begin{align*}
    2^{3+2\alpha}  \frac{\lambda_{\star}^{1+1/\alpha}}{u^{1+1/\alpha} -  2^{2+\alpha+1/\alpha}\lambda_{\star}^{1+1/\alpha} } = 
    \frac{1}{2} \big(1 - e^{-\frac{1}{64}}\big).
\end{align*}
From \eqref{etafp},  bounding $\frac{1}{17 \times 2^{\alpha}} = \frac{1}{2^{4+\alpha}} \frac{2^{4}}{1 + 2^{4}}  \leq \frac{1}{2^{4+\alpha}} \frac{2^{4+3\alpha}}{1 - e^{-\frac{1}{64}} + 2^{4+3\alpha}}
$,
\begin{align*}
    &\frac{\frac{\eta}{\lambda_{\star}} 2^{3+\alpha} \frac{\alpha+1}{\alpha - 1} \lambda_{\star}^{1/\alpha} u^{1+1/\alpha}}
    {u^{1+1/\alpha} -  2^{2+\alpha+1/\alpha}\lambda_{\star}^{1+1/\alpha} } 
    \leq 
    \frac12 \frac{\frac{2^{4+3\alpha}}{1 - e^{-\frac{1}{64}} + 2^{4+3\alpha}} u^{1+1/\alpha}}
    {u^{1+1/\alpha} -  2^{2+\alpha+1/\alpha}\lambda_{\star}^{1+1/\alpha} } \\
    &= 
   \big(1 - e^{-\frac{1}{64}}\big) 
   \frac{(u/\lambda_\star)^{1+1/\alpha}}{1 - e^{-\frac{1}{64}} + 2^{4+3\alpha}}
    \leq \frac{1}{2}\big(1 - e^{-\frac{1}{64}}\big).
\end{align*}
Thus, from \eqref{uub2}, we find that 
$
    A(u)
    \leq 
    1 - \exp(-\frac{1}{64})
$
holds for this value of $u$, hence we have $u_{\min} \leq u$.

\subsection{Ridge Regression with a Logarithmic Grid}
\label{sec:LogGridRidge}
In this section we provide proofs for the results associated to tuning ridge regression with a logarithmic grid. Precisely, let us consider the ridge regression estimator with the sequence of regularization parameters on a logarithmic grid. Thus, for $k \geq j \geq 1$, 
\begin{align*}
    \lambda_j = e^{ \frac{j-1}{k-1} \log (\lambda_{\min})}.
\end{align*}
We have $\lambda_1 = 1$ and $\lambda_k = \lambda_{\min}$ i.e., the sequence is decreasing as a function of the index $j$, and covers the range $[\lambda_{\min},1]$. Naturally, we let the set of estimators then be denoted 
\begin{align*}
    \mathcal{C}^{\text{$\log$-Ridge}}(\lambda_{\min},k) := 
    \{ \widehat{\beta}_{\lambda} : \lambda = \lambda_j, j \in \{1,\dots,k\} \}.
\end{align*}
Given this, we can then bound the ratio of the excess risks of ridge with a logarithmic grid and gradient descent. This is summarized within the proof of Theorem \ref{thm:log-ridgeVsGD} which is provided in the following Section \ref{sec:proof:log-ridge}. This is followed by the proof of Proposition \ref{lem:log_grid} presented within Section \ref{sec:proof:lem:log_grid}.

\subsubsection{Proof of Theorem \ref{thm:log-ridgeVsGD}}
\label{sec:proof:log-ridge}
We split the proof into two subsections: one for the upper bound and one for the lower bound. We begin with the upper bound.

\textbf{Upper Bound}
Consider the lower bound from Proposition  \ref{lem:log_grid} for ridge regression with a logarithmic grid and the upper bound of gradient descent from Proposition \ref{lem:GDBounds:Fine}. Using these results, the ratio of errors is then  upper bounded as
\begin{align*}
    \mathcal{S}\left( \mathcal{C}^{\text{GD}}(\eta,k), 
    \mathcal{C}^{\text{$\log$-Ridge}}(\lambda_{\min},k)
    \right)  
    & \leq 216 \frac{\int G^{\text{GD}}(s) d \widehat{H}(s)}{\int G^{\text{$\log$-Ridge}}(s) d \widehat{H}(s)}.
\end{align*}
Considering the lower bound  
\begin{align*}
    \int G^{\text{$\log$-Ridge}}(s) d \widehat{H}(s) 
    \geq \frac{1}{d} \frac{1}{8} \frac{\log^2(1/\lambda_{\min})}{(k-1)^2} 
    \Big( 
    \lambda_\star \sum_{i : s_i > \lambda_{\star}} \frac{1}{s_i^2} 
    + 
    \frac{1}{\lambda_{\star}^2} \sum_{i : s_i \leq \lambda_{\star}} s_i
    \Big) 
\end{align*}
alongside the equality \eqref{equ:GDSpectralInt} for $\int G^{\text{GD}}(s) d \widehat{H}(s)$ leads to the upper bound 
\begin{align*}
    & \mathcal{S}\left( \mathcal{C}^{\text{GD}}(\eta,k), 
    \mathcal{C}^{\text{$\log$-Ridge}}(\lambda_{\min},k)
    \right) \\ 
    & \leq 
    \frac{1728}{\log^2(1/\lambda_{\min})}
    \Bigg[ 
    k^2 \frac{ \lambda_{\star} \sum_{i: s_i > \lambda_\star} \frac{1}{s_i^2}}{\Big( 
    \lambda_\star \sum_{i : s_i > \lambda_{\star}} \frac{1}{s_i^2} 
    + 
    \frac{1}{\lambda_{\star}^2} \sum_{i : s_i \leq \lambda_{\star}} s_i
    \Big)} \\
    & \quad\quad\quad\quad + 
    k^2 \frac{ \frac{1}{\lambda_{\star}^4} \sum_{i : \lambda_{\star}^2/(k \lambda_{\min}) < s_i \leq \lambda_{\star} } s_i^3}{
    \Big( 
    \lambda_\star \sum_{i : s_i > \lambda_{\star}} \frac{1}{s_i^2} 
    + 
    \frac{1}{\lambda_{\star}^2} \sum_{i : s_i \leq \lambda_{\star}} s_i
    \Big)} 
    +\frac{1}{\lambda_{\min}^2} 
    \frac{\sum_{i: s_i \leq \lambda_{\star}^2/(k \lambda_{\min})}
    s_i}{    \Big( 
    \lambda_\star \sum_{i : s_i > \lambda_{\star}} \frac{1}{s_i^2} 
    + 
    \frac{1}{\lambda_{\star}^2} \sum_{i : s_i \leq \lambda_{\star}} s_i
    \Big) }
    \Bigg]\\
    & \leq \frac{1728}{\log^2(1/\lambda_{\min})}
    \Bigg[ 
    k^2\lambda_\star^3 \frac{  \sum_{i : s_i > \lambda_{\star}} \frac{1}{s_i^2}}{\sum_{i : s_i \leq \lambda_{\star}} s_i}
     + 
    \frac{k^2}{\lambda_\star^2 } \frac{\sum_{i : \lambda_{\star}^2/(k \lambda_{\min}) < s_i \leq \lambda_{\star} } s_i^3}{
    \sum_{i : s_i \leq \lambda_{\star}} s_i
    }
    + \frac{\lambda_{\star}^2}{\lambda_{\min}^2}
    \Bigg].
\end{align*}
Using the bounds \eqref{equ:PartialSumBounds1}, \eqref{equ:PartialSumBounds2} and \eqref{equ:PartialSumBounds3} for $\sum_{i : s_i > \lambda_{\star}} \frac{1}{s_i^2}$, $\sum_{i : s_i \leq \lambda_{\star}} s_i$ and $\sum_{i : \lambda_{\star}^2/(k \lambda_{\min}) < s_i \leq \lambda_{\star} } s_i^3$ respectively, we find
\begin{align*}
    \sum_{i : s_i > \lambda_{\star}} \frac{1}{s_i^2} 
     \leq \lambda_{\star}^{-(2+1/\alpha)},\qquad
    \sum_{i : s_i \leq \lambda_{\star}} s_i  \geq 
    \frac{r^{1-\alpha}}{2},\qquad
    \sum_{i : \lambda_{\star}^2/(k \lambda_{\min}) < s_i \leq \lambda_{\star} } 
    s_i^3 \leq 
    \mathcal{J}_{\alpha,\lambda_\star}(k,\lambda_{\min}).
\end{align*}
Hence
\begin{align*}
     \mathcal{S}\left( \mathcal{C}^{\text{GD}}(\eta,k), 
    \mathcal{C}^{\text{$\log$-Ridge}}(\lambda_{\min},k)
    \right)  
     & \leq
    \frac{3456}{\log^2(1/\lambda_{\min})}\\
    & \quad\quad\quad\quad 
    \times 
    \Big[ \frac{k^2}{r^{1-\alpha}} \Big( \frac{1}{\lambda_{\star}^{1/\alpha-1}} + \frac{1}{\lambda_{\star}^2} \mathcal{J}_{\alpha,\lambda_{\star}}(k,\lambda_{\min}) \Big) + \frac{\lambda_{\star}^2}{\lambda_{\min}^2} 
    \Big].
\end{align*}
Following the proof of Theorem \ref{thm:Informal}, choosing 
\begin{align*}
    \overline{r}_{\alpha, \lambda_{\star},k,\lambda_{\min}}
    = 
    \Big[ \frac{\lambda_{\min}^2k^2 }{\lambda_{\star}^2} \Big( \frac{1}{\lambda_{\star}^{1/\alpha -1}} + \frac{1}{\lambda_{\star}^2} \mathcal{J}_{\alpha,\lambda_{\star}}(k,\lambda_{\min}) \Big) 
     \Big]^{1/(1-\alpha)},
\end{align*}
we see that once $r \geq \overline{r}_{\alpha, \lambda_{\star},k,\lambda_{\min}}$ we find 
\begin{align*}
    \mathcal{S}\left( \mathcal{C}^{\text{GD}}(\eta,k), 
    \mathcal{C}^{\text{$\log$-Ridge}}(\lambda_{\min},k)
    \right)  
     \leq
     \frac{6912}{\log^2(1/\lambda_{\min})} 
     \frac{\lambda_{\star}^2}{\lambda_{\min}^2},
\end{align*}
as required. Plugging in the definition of $\lambda_{\min} = d\sigma_{\min}^2/n$ and $\lambda_{\star} = d \sigma^2/n$ then yields the result. 

\textbf{Lower Bound}
For ridge regression with a logarithmic grid, considering the upper bound from Proposition \ref{lem:log_grid} to get 
\begin{align*}
    \mathcal{E}(\mathcal{C}^{\text{$\log$-Ridge}}(\lambda_{\min},k))
    & \leq  \lambda_{\star} \int G^{\text{$\log$-Ridge}}(s) d \widehat{H}(s)\\
    & \leq 
    \frac{\lambda_{\star}}{d} \frac{\log^2(1/\lambda_{\min})}{(k-1)^2} 
    \Big( 
    \lambda_\star \sum_{i : s_i > \lambda_{\star}} \frac{1}{s_i^2} 
    + 
    \frac{1}{\lambda_{\star}^2} \sum_{i : s_i \leq \lambda_{\star}} s_i
    \Big)\\ 
    & \leq 
    2 \frac{\lambda_{\star}}{d} \frac{\log^2(1/\lambda_{\min})}{(k-1)^2}
    \frac{1}{\lambda_{\star}^2} \sum_{i : s_i \leq \lambda_{\star}} s_i
\end{align*}
where the second equality comes from $\sum_{i: s_i < \lambda_{\star}} s_i \geq \lambda_{\star}^{3} \sum_{i : s_i > \lambda_{\star}} \frac{1}{s_i^2}$ when $r \geq 2^{1/(1-\alpha)}(1+\lambda_{\star}^{-1/\alpha})$, see equation  \eqref{equ:eigenval_tail_lower_bound}.

Meanwhile, for gradient descent follow the proof of Theorem \ref{thm:mainGDvsRidge_Lower}, precisely, combine equations \eqref{equ:GD:lower_calc_1} and \eqref{equ:GD:lower_calc_2} which arise from applying Proposition \ref{lem:GDBounds:Fine}, to get the lower bound
\begin{align*}
    \mathcal{E}(\mathcal{C}^{\text{GD}}(\eta,k))
    & \geq 
        \frac{1}{16} 
    \min\{1-\kappa,\kappa\}^2
    \lambda_{\star}
    \left( \frac{1}{\lambda_{\min} k} \right)^2
    \frac{1}{d}
    \sum_{i : s_i \leq \frac{\lambda_{\star}^2 \kappa }{4 k \lambda_{\min}} }
    s_i
\end{align*}

Taking the ratio of the risks and applying the upper and lower bounds yields 
\begin{align*}
    \mathcal{S}\left( \mathcal{C}^{\text{GD}}(\eta,k), 
    \mathcal{C}^{\text{$\log$-Ridge}}(\lambda_{\min},k)
    \right)  
    & \geq \frac{\min\{1-\kappa,\kappa\}^2}{32\log^2(1/\lambda_{\min})}  \Big( 1 - \frac{1}{k}\Big)^2
    \frac{\lambda_{\star}^2}{\lambda_{\min}^2} 
    \frac{\sum_{i : s_i \leq \frac{\lambda_{\star}^2 \kappa }{4 k \lambda_{\min}} }
    s_i}{
    \sum_{i : s_i \leq \lambda_{\star}} s_i
    } \\
    & \geq 
    \frac{\min\{1-\kappa,\kappa\}^2}{128\log^2(1/\lambda_{\min})} 
    \frac{\lambda_{\star}^2}{\lambda_{\min}^2}  
    \Big( 
    1 \! - \!
    \frac{2}{1 \! - \! \alpha} \left( \frac{4 k \lambda_{\min}}{\lambda_{\star}^2 \kappa }\right)^{1/\alpha - 1}
    \!\! \frac{1}{r^{1-\alpha}}
    \Big)
\end{align*}
where the second inequality arises from $k \geq 2$ and the lower bound on the ratio of series \eqref{equ:series_ratio_lower_bound}. 
If we then let 
\begin{align*}
    \underline{r}_{\alpha,\lambda_{\star},k,\lambda_{\min}} 
    = \max\Big\{ 2^{1/(1-\alpha)}(1+\lambda_{\star}^{-1/\alpha}), \Big(\frac{4}{1-\alpha}\Big)^{1/(1-\alpha)} \left( \frac{4 k \lambda_{\min}}{\lambda_{\star}^2 \kappa }\right)^{1/\alpha}
    \Big\}
\end{align*}
then we have the desired lower bound once $r \geq \underline{r}_{\alpha,\lambda_{\star},k,\lambda_{\min}}$.

\subsubsection{Proof of Proposition \ref{lem:log_grid}}
\label{sec:proof:lem:log_grid}
Since $\lambda_\star \in (\lambda_{\min},1)$ there exists a $2 \leq j_{\star} \leq k$ such that $\lambda_{j_{\star}} < \lambda_{\star} < \lambda_{j_{\star} -1 }$. Furthermore, there exists an $\epsilon \in (0,1)$ such that 
$\lambda_{\star} = \exp\big( \frac{j_{\star} - 2 + \epsilon}{k-1} \log(\lambda_{\min}) \big)$. Leveraging this alongside the risk being monotonic yields 
\begin{align*}
    \mathcal{E}(\mathcal{C}^{\text{$\log$-Ridge}}(\lambda_{\min},k)) 
    & = \min_{j=1,\dots,k}
    \E_{\beta_{\star},\epsilon}[L_{\beta_{\star}}(\widehat{\beta}_{\lambda_j})]
    -
    \lambda_\star \int \frac{1}{s + \lambda_{\star}}
    d \widehat{H}(s) \\
    & = \min_{j=1,\dots,k}
    \frac{\sigma^2}{n} \sum_{i=1}^{r} \left( \frac{1}{s_i} + \frac{1}{\lambda_{\star}} \right) 
    \left( \frac{\lambda_j}{\lambda_j + s_i} - \frac{\lambda_{\star}}{\lambda_{\star} + s_i} \right)^2\\
    & = \min_{j =1,\dots,k} 
    \frac{\sigma^2}{n} \sum_{i=1}^{r} \left( \frac{1}{s_i} + \frac{1}{\lambda_{\star}} \right) 
    \left( \frac{\lambda_j}{\lambda_j + s_i} - \frac{\lambda_{\star}}{\lambda_{\star} + s_i} \right)^2\\
    & = 
    \min_{j = j_{\star}, j_{\star}-1}
    \frac{\sigma^2}{n} \sum_{i=1}^{r} \left( \frac{1}{s_i} + \frac{1}{\lambda_{\star}} \right) 
    \frac{s^2_i}{(\lambda_j + s_i)^2(\lambda_{\star} + s_i)^2} 
    (\lambda_j - \lambda_{\star})^2 \\
\end{align*}
We now set to bound this minimum for the two cases $j=j_{\star},j_{\star} - 1$. Begin by noting that
\begin{align*}
    |\lambda_{j_{\star}} - \lambda_{\star}| & = 
    \lambda_{\star} \Big| \frac{\lambda_{j_{\star}}}{\lambda_{\star}} - 1 \Big|\\
    & = \lambda_{\star}\big| \exp\Big( - \frac{1 - \epsilon}{k-1}\log(1/\lambda_{\min}) \Big) - 1\big|\\
    & = \lambda_{\star}\Big( 1 - \exp\Big( - \frac{1 - \epsilon}{k-1}\log(1/\lambda_{\min}) \Big) \Big), \\
    |\lambda_{j_{\star} - 1} - \lambda_{\star}| & = 
    \lambda_{j_{\star} - 1} 
    \Big|1 - \frac{\lambda_{\star}}{\lambda_{j_{\star} - 1}} \Big|
    = \lambda_{j_{\star} - 1}
    \Big(1 - \exp\big( - \frac{\epsilon}{k-1}\log(1/\lambda_{\min})\big)\Big).
\end{align*}
Moreover, since $k \geq 1 +3\log(1/\lambda_{\min}) \geq 1 + \frac{\epsilon\log(1/\lambda_{\min})}{\log(2)}$  and $\lambda_{\star} = \exp\big( -\frac{\epsilon \log(1/\lambda_{\min})}{k-1}\big)\lambda_{j_{\star} - 1}$ we have $2 \lambda_{\star} \geq \lambda_{j_{\star} - 1} \geq \lambda_{\star}$, and thus, can bound for $i=1,\dots,r$
\begin{align*}
    \frac{\lambda_{\star}}{\lambda_{\star} + s_i} 
    \leq \frac{\lambda_{j_{\star} - 1}}{\lambda_{j_{\star} - 1} + s_i} \leq 
    2 \frac{\lambda_{\star}}{\lambda_{\star} + s_i}.
\end{align*}
Similarly, $\lambda_{\star} = \exp\big( \frac{(1-\epsilon)\log(1/\lambda_{\min})}{k-1}\big) \lambda_{j_{\star}} \leq 2 \lambda_{j_{\star}}$ and thus $\frac{1}{2} \lambda_{\star} \leq \lambda_{j_{\star}} \leq \lambda_{\star}$ leading to 
\begin{align*}
    \frac{1}{\lambda_{\star}  + s_i} \leq \frac{1}{\lambda_{j_{\star}} + s_i} \leq \frac{2}{\lambda_{\star} + s_i}
\end{align*}
Utilising these two facts as well as the inequality $2 \geq u \geq 0$, $\frac{u}{3} \leq 1-e^{-u} \leq u$, since $k \geq 1 +  \log(1/\lambda_{\min})$ we get the desired upper and lower bounds. 

\section{Proofs for Orthogonal Design}

In this section we present the proofs associated to the orthogonal design setting. This section is structured as follows. Section \ref{sec:mmxgridpf} presents the proof of Theorem \ref{mmxgrid}. Section \ref{pf:err} presents the proof of  Proposition \ref{err}. Section \ref{sec:ergdpf} presents the proof of Proposition \ref{ergd}.

\subsection{Proof of Theorem \ref{mmxgrid}}
\label{sec:mmxgridpf}
By \eqref{errdec}, using $n=d$ and $\sigma=1$, the excess risk of an estimator $\Phi_j$ with $\phi_j=\Phi_j(s)$  is
\begin{align*}
    & \frac{(1 + \psi s )}{s} \left( \phi_j - \frac{1}{1 + \psi s} \right)^2
     =  \frac1s \left(\sqrt{1 + \psi s }\cdot \phi_j - \sqrt{\frac{1}{1 + \psi s }} \right)^2.
\end{align*}
Since $\psi \in [\psi_-,\psi_+]$, for $x := \sqrt{1 + \psi s }$ we have $x \in [x_-,x_+]$.
We can equivalently write the minimax excess risk objective \eqref{mmxrisk} as
$$\inf_{T_k  = \{\phi_1,\ldots,\phi_k\}\in E} \,\sup_{x\in [x_-,x_+]} \,\min_{j\in[k]}\,
 \left(x \phi_j - \frac{1}{x} \right)^2.
$$
Moreover, by monotonicity of the transform $z\to z^2$ for $z\ge 0$, we can equivalently study the objective
\begin{align}\label{redmmx}
\inf_{T_k  = \{\phi_1,\ldots,\phi_k\}\in E} \,\sup_{x\in [x_-,x_+]} \,\min_{j\in[k]}\,
\left|x \phi_j - \frac{1}{x} \right|.
\end{align}

Define the function $g: (0,\infty) \times (0,\infty) \mapsto [0,\infty)$ by  $g(x,l) = |x l - \frac{1}{x} \big|$. For any fixed $l>0$, the one-variable function  $x\to g(x,l)$ has a unique zero for $x^*(l) = 1/\sqrt{l}$. Moreover, $g(x,l) = x l - \frac{1}{x}$ for $x \ge x^*$ and $g(x,l) =  \frac{1}{x} -x l$ for $x < x^*$.
Further, the map $l\mapsto x^*(l)$ is monotonically decreasing for $l\in(0,\infty)$. 
Thus, if $x>x^*(l)$, then $g(x,l) < g(x,l')$.
Similarly, for $l <l'$ and  $x<x^*(l')$, we have $g(x,l) > g(x,l')$. 

In the interval $x\in[x^*(l'),x^*(l)]$, we have $g(x,l)  \le g(x,l')$ if and only if
\begin{align*}
    \frac{1}{x} -x l & \le x l' - \frac{1}{x},
\end{align*}
or also $\sqrt\frac{2}{l+l'} \le x$.

Thus, for a grid $\phi_1 \le \phi_2 \le\ldots \le \phi_k$, defining $q_i = \sqrt{\frac{2}{\phi_i+\phi_{i+1}}} $, with $q_k:=0$, $q_{0}=\infty$, and thus 
$$q_k\le q_{n-1} \le \ldots \le q_1 \le q_{0},$$
we have that the minimum of the functions $g(x,\phi_i)$, $i=0,\ldots,k$ has a value equal to $g(x,\phi_j)$ precisely for $x \in I_j:= [q_j,q_{j-1}]$. In addition, the functions $g(x,\phi_j)$, $g(x,\phi_{j+1})$ are equal at $q_j$, $j=1,\ldots,k-1$.

Moreover, the value at the maximum on the interval $I_j$ (for $j=1,\ldots,k-1$) can be verified to be 
$$g(q_j,\phi_j) = g(q_j,\phi_{j+1}) = \frac{\phi_{j+1} -\phi_{j}}{\sqrt{2(\phi_{j+1} + \phi_{j})}}.$$

In \eqref{redmmx}, we need to find the supremum for $x\in [x_-,x_+]$, so we can restrict to grid points such that $1/x_+^2 \le \phi_1$, $\phi_k \le 1/x_-^2$. Since we also need to consider the boundary values, problem \eqref{redmmx} simplifies to
\begin{align}
\label{kp1}
& \inf_{x_+^{-2} \le \phi_1 \le \ldots \le \phi_k\le x_-^{-2}} \,
m(\phi_1,\ldots,\phi_k) \\
& :=
 \max \left( \frac{1}{x_-} -x_- \phi_k, \max_{j=1,\ldots,k-1} 
\frac{\phi_{j+1} -\phi_{j}}{\sqrt{2(\phi_{j+1} + \phi_{j})}}, x_+ \phi_1-\frac{1}{x_+} \right).
\end{align}
Now, we provide a compactness and ``infinite descent" argument to find the minimizer and prove that the minimum is achieved.  
We claim that the solution is achieved when the $k+1$ quantities in \eqref{kp1} are equal. Otherwise, we argue that we can decrease the objective. Suppose that the $k+1$ terms are not all equal. Then there is an index $j \in \{2,\ldots, k-1\}$ such that $\max(t_{j-1},t_{j+1}) \le t_j$, with $t_j=m$ achieving the maximum, and at least one inequality is strict. Suppose without loss of generality that  $t_{j-1} < t_j$. 
Then, $t_j = t_j(\phi_{j-1},\phi_j)$ is strictly increasing in $\phi_j$, and strictly decreasing in $\phi_{j-1}$. This is readily verified by elementary calculus. Thus, if we decrease $\phi_j$, then we decrease $t_j$, while increasing $t_{j-1}$. By continuity, after decreasing $\phi_j$ sufficiently, we will have $t_{j-1} = t_j$. This procedure keeps $\phi_j \ge \phi_{j-1}$, because $t_j(\phi_{j-1},\phi_{j-1})=0$, which is less than $m$, as $m>0$ must be strictly positive (otherwise all terms must be equal, which is a contradiction).

This shows that if not all $t_j$ are equal, then we can decrease $m$. Now, by the compactness of $x_+^{-2} \le \phi_1 \le \ldots \le \phi_k\le x_-^{-2}$  and the continuity of $m$, $m$ has a minimum that is attained. Since the minimum cannot be attained when any $t_j$ are unequal, at the minimizing $\phi_j$, all $t_j$ are equal.

In conclusion, the minimax optimal grid is characterized as a sequence
$$ \frac{1}{1 + \psi s_+} \le \phi_1 \le\ldots \le \phi_k\le  \frac{1}{1 + \psi s_-}$$
such that 
$$\frac{1}{x_-} -x_- \phi_k
=  x_+ \phi_1-\frac{1}{x_+} 
= 
\frac{\phi_{j+1} -\phi_{j}}{\sqrt{2(\phi_{j+1} + \phi_{j})}}, \, j=1,\ldots,k-1.
$$

Let  $\phi_1 \ge x_+^{-2}$ and $c=x_+ \phi_1-\frac{1}{x_+}\ge 0$. If $c=0$, then all grid points $\phi_j$ must be equal, thus $\psi_-=\psi_+$, and the conclusion clearly holds. Otherwise, $c>0$ and $\phi_1 > x_+^{-2}$. 
The above display implies the quadratic equation
$$\phi_{j+1}^2 - 2\phi_{j+1}(\phi_{j} + c^2) + \phi_{j}^2-2c^2\phi_{j}=0$$
for $\phi_{j+1}$.
This has the only solution $\phi_{j+1}>\phi_j$ given by
\begin{align*}
  \phi_{j+1} 
    &= \phi_{j} + c \left( c + \sqrt{c^2 + 4\phi_{j}} \right).
\end{align*}
Then $u_j := c^2 + 4\phi_{j}$ satisfies
\begin{align*}
   u_{j+1} = u_{j} + 4c (c + \sqrt{u_j})
   = (\sqrt{u_j}+2c)^2.
\end{align*}
Thus $v_j := \sqrt{u_j}/2$ is an arithmetic progression with $v_{j+1} = v_j+c$.
From $c=  x_+ \phi_1-1/x_+$ and $v_1 = \sqrt{\phi_{1}+c^2/4}$, we find
$v_1= \frac{1}{{x_+}}+\frac{c}{2}$; and similarly $v_k= \frac{1}{{x_-}}-\frac{c}{2}$.
Using $v_k = v_1 + (k-1)c$, we find $c = \frac1k\left(\frac{1}{{x_-}}-\frac{1}{{x_+}}\right)$. This finishes the proof.

\subsection{Proof of Proposition \ref{err}}
\label{pf:err}

{\bf Minimizing over $j$, for a fixed $x$.}
With the notations from Theorem \ref{mmxgrid}, we need to evaluate
$$\sup_{x\in [x_-,x_+]} \,\min_{j\in[k]}\,
 \left(x \phi_j - \frac{1}{x} \right)^2
$$
for $\phi_j = (j/k)/[(j/k)+s]= j/(j+a)$, where $a=ks$. For a fixed $x$ the minimum over $j$ is achieved when $\phi_j$ is closest to $1/x^2$, i.e., at either $\lfloor b \rfloor$ or $\lceil b \rceil$ (where $\lceil b \rceil$ denotes the ceiling of $b$), with
$$b:=\frac{a}{x^2-1},
$$
or, the unique integer $j\in[k]$ closest to these two values if neither belongs to $[k]$. One can verify that the minimum is achieved at $\lfloor b \rfloor$ or $\lceil b \rceil$ iff $\lfloor b \rfloor \le k$, or equivalently $b<k+1$. Since $\psi_->k/(k+1)$, this condition holds for all $x\in [x_-,x_+]$.
As in the proof of Theorem \ref{mmxgrid}, we can equivalently study the objective, for $j=1,\ldots,k$
\begin{align*}
         o_j:= \left|x \phi_j - \frac{1}{x} \right|. 
\end{align*}
The minimizer is either $j_b= \lfloor b \rfloor$ or $j_b+1$. Writing $b = j_b+\tau_b$, for some $\tau_b\in[0,1)$, using that $b+a = ax^2/(x^2-1)$ the objective at $j_b$ is
\begin{align*}
         o_{j_b} &= x \left|\frac{j_b}{j_b+a} - \frac{1}{x^2} \right|
         = xa\frac{\tau_b}{(j_b+a)(b+a)}
         = \frac{\tau_b}{j_b+a} \frac{x^2-1}{x}.
\end{align*}
We can also write $x^2-1=a/b$ and $x = (1+a/b)^{1/2}$. Using $b = j_b+\tau_b$, we can express
\begin{align*}
         o_{j_b} &= \frac{\tau_b}{j_b+a} \frac{a/b}{(1+a/b)^{1/2}}
         = \frac{\tau_b}{j_b+a} \frac{a/[j_b+\tau_b]}{(1+a/[j_b+\tau_b])^{1/2}}.
\end{align*}
Similarly, one can verify that the objective at $j_b+1$ is
\begin{align*}
         o_{j_b+1} &=  \frac{1-\tau_b}{j_b+1+a}\frac{a/[j_b+\tau_b]}{(1+a/[j_b+\tau_b])^{1/2}}.
\end{align*}
The objective for a fixed $x$ (or equivalently for a fixed $b$) equals the minimum of $o_{j_b}$ and $o_{j_b+1}$. Since $x$ ranges in $x\in[x_-,x_+] = [\sqrt{1+s\psi_-},\sqrt{1+s\psi_+}]$, we have that $b=sk/(x^2-1)$ ranges in $b\in[b_-,b_+] = [k/\psi_+,k/\psi_-]$. Thus, the objective can be written as
\beq\label{Q}
\sup_{b\in [k/\psi_+,k/\psi_-]} \,
Q(j_b,\tau_b):=
\frac{a/[j_b+\tau_b]}{(1+a/[j_b+\tau_b])^{1/2}}
 \min\left(\frac{\tau_b}{j_b+a},\frac{1-\tau_b}{j_b+1+a}\right)
\eeq
where $b=j_b+\tau_b$, and $j_b = \lfloor b\rfloor$ is the integer part of $b$, while $\tau_b\in[0,1)$ is the fractional part.

{\bf Maximizing over $x$.} To maximize over $x$, it is enough to maximize over $b=j_b+\tau_b$, and we maximize over each component separately.
It is readily verified that the objective $Q$ from \eqref{Q} is $\emph{decreasing}$ in $j_b$, for any fixed $\tau_b$. Hence, it achieves its maximum for the smallest possible integer $j_b$ such that $b=j_b+\tau_b\in [k/\psi_+,k/\psi_-]$, which is $j_b = \lfloor k/\psi_+\rfloor$ or  $j_b = \lfloor k/\psi_+\rfloor+1$, depending on a condition on $\tau_b$. 
Due to the assumption $\lfloor k/\psi_+ \rfloor +1 \le k/\psi_-$, the condition is $\tau_b\in[\{ k/\psi_+ \},1)$ for $j_b = \lfloor k/\psi_+\rfloor$, and $\tau_b\in[0,k/\psi_--(\lfloor k/\psi_+ \rfloor +1)]$, and $\tau_b<1$ for $j_b = \lfloor k/\psi_+\rfloor+1$.

One can readily verify that for any fixed $j\ge 0$, with $\tau^*= \tau^*(j,a): = (j+a)/[2(j+a)+1]$, the function $Q(j,\tau)$ is strictly increasing in $\tau_b$ on the interval $[0,\tau^*)$, and strictly decreasing on the interval $(\tau^*,1)$. 

Finally, we maximize over $\tau_b$, for a fixed $j_b$. 
Then, for $j_b=\lfloor k/\psi_+\rfloor$ it is readily verified that the maximum over $\tau_b$  in the allowed set is achieved at
\begin{align*}
         \tau^*_{j_b} &= \max\left(\{ k/\psi_+ \},\frac{\lfloor k/\psi_+\rfloor+a}{2(\lfloor k/\psi_+\rfloor+a)+1}\right).
\end{align*}
Moreover, for $j_b=\lfloor k/\psi_+\rfloor+1$ the maximum over $\tau_b$ is achieved at 
\begin{align*}
         \tau^*_{j_b+1} &=
         \min\left(k/\psi_--(\lfloor k/\psi_+ \rfloor +1),\frac{\lfloor k/\psi_+\rfloor+1+a}{2(\lfloor k/\psi_+\rfloor+1+a)+1}\right)
         .
\end{align*}

The optimal objective can be checked to be
$$
O_*= \max \left \{ Q(j_b,\tau^*_{j_b}), Q(j_b+1,\tau^*_{j_b+1}) \right\}.
$$
Squaring and multiplying by $1/s$, we find the excess risk of ridge regression, finishing the proof.

\subsection{Proof of Proposition \ref{ergd}}
\label{sec:ergdpf}
Following the proof of Proposition \ref{err}, we need to evaluate
\begin{align*}
    E:=\sup_{x \in [x_{-},x_{+}]} \min_{j=1,2,\dots,k} \left( x b^{j-1} - \frac{1}{x}\right)^2,
\end{align*}
where $b= 1-s \eta$. Let $h_{x} = 2\log_{1/b}(x)$ and $\tau_x = \{ h_{x}\}$. The inner minimum is achieved at either $j(x) = \lfloor h_{x} \rfloor +1$ or $j(x) + 1$. Writing $j(x) = h_{x} - \tau_x +1$, we find 
\begin{align*}
\min_{j=1,2,\dots,k} \left| x b^{j-1} - \frac{1}{x}\right| 
= 
\min\left\{ 
x b^{ h_x - \tau_x } - \frac{1}{x}, 
\frac{1}{x} - x b^{h_x - \tau_x + 1}
\right\}
= 
\frac{1}{x} \min\big\{ b^{-\tau_x} - 1, 1 - b^{1-\tau_x}\big\}
\end{align*}
since $xb^{h_x} = \frac{1}{x}$. Factoring out $b^{-\tau_x}$ from the above and squaring we can then write
\begin{align*}
    E
    & =
    \sup_{x \in [x_{-},x_{+}]}
    \frac{1}{x^2} b^{-2\tau_x}
    \min\big\{ 1 - b^{\tau_x}, b(b^{\tau_x - 1} - 1) \big\}^2\\
    & = 
    \sup_{x \in [x_{-},x_{+}]}
    b^{j(x) - 1 - \tau_x}
   \min\big\{ 1 - b^{\tau_x}, b(b^{\tau_x - 1} - 1) \big\}^2,
\end{align*}
where for the second equality we note that $x^{-2} = b^{h_x} = b^{j(x) + \tau_x - 1}$. 
The fractional part that maximizes the minimum for $\tau_x \in [\{h_x\},1)$ occurs at $\tau^* = \log_{1/b}(2/(b+1))$. Thus, with
\begin{align*}
    Q_\star(\tau) = \min\{ 1 - b^{\tau}, b (b^{\tau - 1} - 1) \},
\end{align*}
we find  since $(1-\eta)^{k} < (1+\psi_{+})^{-1}$ that
\begin{align}\label{e-eq}
    E
    & =
    b^{\lfloor h_{x_{-}} \rfloor  - \tau^*} \max\bigg\{  \max_{\tau \in [ \{h_{x_{-}}\},1)} Q_\star(\tau)^2, 
    b\cdot \max_{\tau \in [0,\{h_{x_{-}}\}) } Q_\star(\tau)^2
    \bigg\}.
\end{align}
Where we note that $b^{j(x)}$ is decreasing in $x$, therefore, the first term within the maximum above is associated to choosing $x$ over the range $[x_{-},\lceil x_{-} \rceil)$, and the second is associated to choosing $x$ over the range $[\lceil x_{-}\rceil, \lceil x_{-} \rceil +1) $.  

Since $Q_{\star}(\cdot)$ attains its maximum at $\tau^{\star}$ and is decreasing on either side of $\tau^{\star}$, we can rewrite \eqref{e-eq}, with $a \vee b = \max\{a,b\}$ as well as $a \wedge b = \min\{a,b\}$,  as follows 
\begin{align*}
    E
    & =
    b^{\lfloor h_{x_{-}} \rfloor  - \tau^*} \max\big\{   Q_\star(\{h_{x_{-}}\} \vee \tau^{\star})^2, 
    b\cdot Q_\star(\{h_{x_{-}}\} \wedge \tau^{\star})^2
    \big\}.
\end{align*}
This quantity can take a total of four values as a function of $\{h_{x_{-}}\}$, $\tau^{\star}$ and $b$.

%\newpage
\bibliographystyle{alpha}
\bibliography{References}

\end{document}